\documentclass[11pt,openright,a4paper]{report}

\usepackage[latin1]{inputenc}
\usepackage[english]{babel}
\usepackage[margin=20pt]{caption}
\usepackage[T1]{fontenc}
\usepackage[centertags]{amsmath}
\usepackage{mathrsfs}
\usepackage{graphicx}   
\usepackage{subfig} 
\usepackage{amsmath}
\usepackage{amsthm}
\usepackage{amsfonts}
\usepackage{amssymb}
\usepackage{wasysym}
\usepackage{ae}  
\usepackage{fancyhdr}
\usepackage{dsfont}
\usepackage{indentfirst}

\usepackage{helvet}

\theoremstyle{plain}
\newtheorem{theorem}{Theorem}[section]
\newtheorem{lemma}{Lemma}[section]
\newtheorem{proposition}{Proposition}[section]
\newtheorem{problem}{Problem}[section]
\theoremstyle{definition}
\newtheorem{definition}{Definition}[section]
\newtheorem{remark}{Remark}[section]
\newtheorem{example}{Example}

\fancypagestyle{plain}{\fancyhead{}}
\begin{document}

\hyphenation{ }

\newcommand{\bb}[1]{\bold{#1}}


%
%
%
 
\begin{titlepage}
    \begin{center}
       \begin{figure}[htbp]
            \begin{center}
            \includegraphics[scale=0.2]{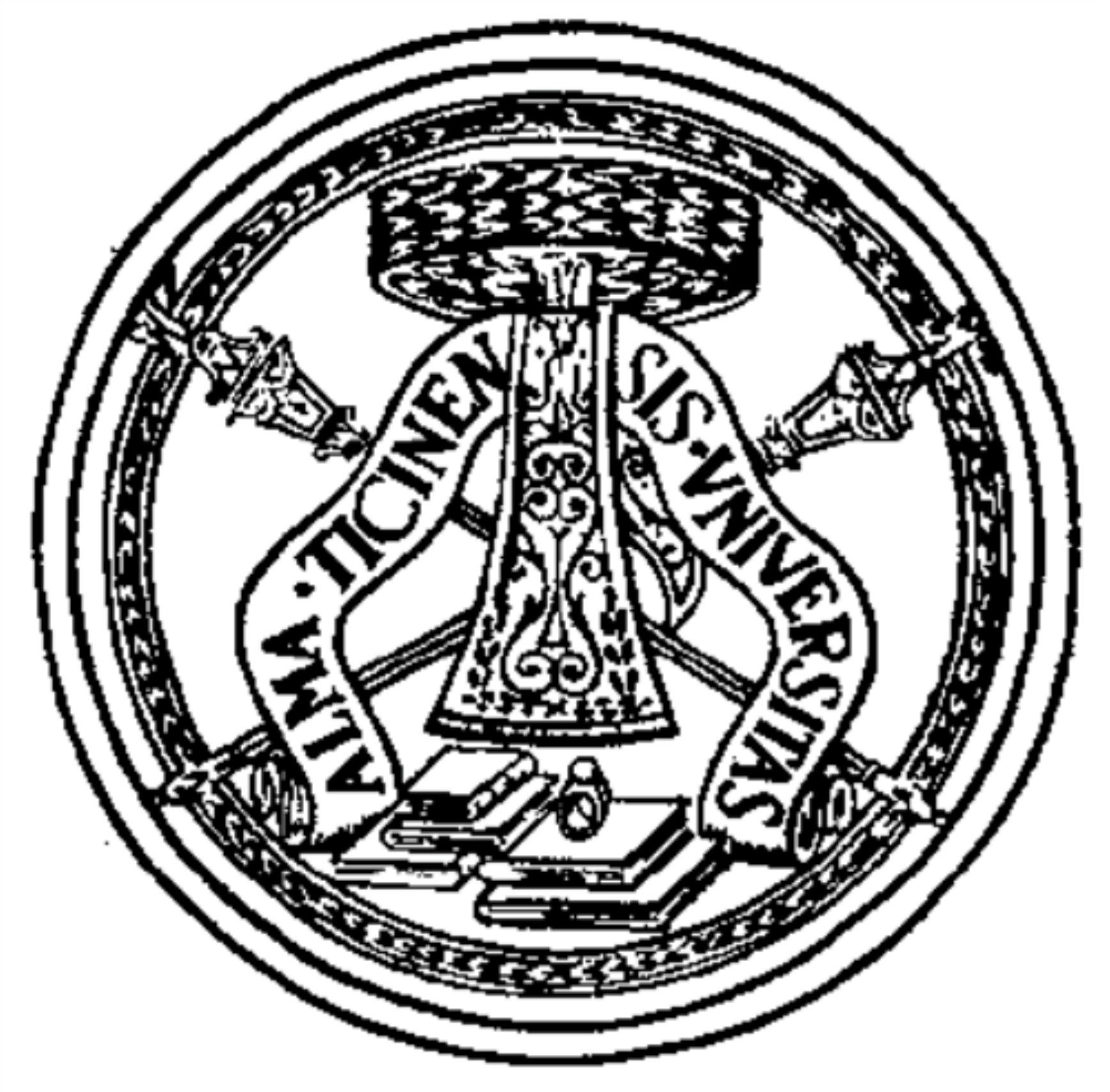}
            \end{center}
       \end{figure}
        \Large{\textbf{UNIVERSIT\`A DEGLI STUDI DI PAVIA}}\\
			\vspace*{\stretch{0.2}}
        \hrule 
			\vspace*{\stretch{0.2}}
        \large{\textbf{FACOLT\`A DI SCIENZE MM. FF. NN.}} \\
        \normalsize{\textbf{CORSO DI LAUREA SPECIALISTICA IN SCIENZE FISICHE}} \\
          \vspace{0.6cm}

       \LARGE{Flussi gradiente generati da un potenziale di interazione repulsivo non liscio}\\
       {Gradient flows driven by a non-smooth repulsive interaction potential}
        \vspace*{\stretch{1}}
    \end{center}
        \begin{flushleft}
            \large{Tesi di laurea di} \\ \Large{Giovanni Bonaschi} \\
        \end{flushleft}
        \vspace*{\stretch{0.5}}
        \begin{flushright}
            \large{Supervisore} \\ 
            \Large{Chiar.mo Prof. Giuseppe Savar\'{e}} \\
                    {Dipartimento di Matematica}
                   \end{flushright}
    \vspace*{\stretch{0.8}}
    \begin{center}
        \large{Anno Accademico 2009/2010}
    \end{center}
\end{titlepage}

\renewcommand{\thepage}{\Roman{page}} 
\setcounter{page}{0}

\clearpage{\pagestyle{empty}\cleardoublepage}
\begin{flushright}
\textit{Alla nonna}
\end{flushright}
  
\vspace{\stretch{0.8}}
  
\clearpage{\pagestyle{empty}\cleardoublepage}
        \begin{flushright}
            \textit\emph{In re mathematica ars proponendi pluris facienda est quam solvendi.} \\
 				\bigskip          
            \textit{(Georg Cantor)}
        \end{flushright}

\vspace{\stretch{0.8}}

\clearpage{\pagestyle{empty}\cleardoublepage}

\chapter*{Abstract}

Nell'ultimo secolo è emerso che una enorme varietà di fenomeni, spaziando dalla fisica alle scienze sociali, può essere descritta tramite 
 equazioni alle derivate parziali. Sviluppare metodi per ottenere una soluzione, esatta o approssimata, è di cruciale importanza. In 
questa tesi abbiamo considerato esempi appartenenti ad un classe abbastanza generale di PDE, riguardanti la derivata prima nella coordinata 
 temporale e fino al secondo ordine in quelle spaziali. Queste equazioni presentano la tipica 
 struttura dei "Flussi Gradiente": una teoria che è stata sviluppata negli ultimi decenni, partendo dall'esempio dei sistemi dinamici 
 in spazi finito dimensionali, e possiede una grande varietà di applicazioni. Tra queste, di particolare importanza, si trova un modello 
 cinetico per flussi granulari. Questo modello descrive materiali granulari (sabbia, polvere, molecole o persino asteroidi...) ed è stato 
 analizzato in questa tesi nel caso unidimensionale e spazialmente omogeneo. Esso è, sotto queste due condizioni, governato dall'equazione 
 \begin{equation*}
 \partial_t f(v,t)=-\kappa \partial_v(Ff)+ \sigma \partial_v(vf) + T_e \partial_v^2f,
 \end{equation*} 
 dove $F=f \ast \nabla W$. Nel lato destro dell'equazione troviamo termini rispettivamente di tipo collisionale , di deriva e di bagno termico.
 Partendo dal modello è stato studiato il termine con $F$, analizzando in particolare l'evoluzione di una 
 distribuzione $\mu(t)$, continua o discreta, di massa autointeragente tramite il potenziale di interazione $W$:
 \begin{equation*}
 \frac{\partial \mu}{\partial t}=\nabla \cdot (\mu (\nabla W \ast \mu)) \qquad t > 0.
 \end{equation*} 
Il potenziale che è stato analizzato presenta un punto angoloso concavo nell'origine mentre all'infinito 
è limitato ad una crescita quadratica; è quindi di tipo repulsivo a brevi distanze. 
La teoria dei flussi gradiente è ben consolidata per potenziali che sono perturbazioni regolari di funzioni convesse. Il tipo di 
potenziale analizzato nella tesi non rientra in questa classe a causa 
del comportamento nell'origine. Siamo comunque riusciti a dimostrare, almeno nel caso unidimensionale ($d=1$), l'esistenza e unicità della 
soluzione. In dimensione maggiore di uno il problema presenta maggiori difficolt\`a e verrà affrontato in un lavoro successivo. \\

La tesi si snoda in quattro capitoli. Il primo introduce la teoria del trasporto ottimo, formulata da Monge nel 1781 e che presenta tante
 profonde connessioni a problemi fisici, matematici ed economici. Grazie a questa teoria, oltre ad altri importanti risultati, possiamo 
 definire una distanza nello spazio delle misure di probabilità: la distanza di Wasserstein. \\

Il secondo capitolo introduce e riassume la teoria dei flussi gradiente. Dopo la definizione e i principali risultati, tra cui l'esistenza e 
unicità della soluzione, viene descritto lo schema dei movimenti minimizzanti. Questo schema è una procedura di approssimazione della soluzione di notevole 
importanza sia dal punto di vista teorico, sia da quello "pratico" (approssimazione numerica). \\

Il terzo capitolo si compone di due parti. Inizialmente viene introdotto e descritto il modello cinetico per 
flussi granulari. Nella seconda parte viene presentato il lavoro 
di Carrillo \emph{et.al.}, che ha ispirato, con un metodo di dimostrazione alternativo, il lavoro di questa tesi. \\

Nel quarto capitolo viene analizzato il caso del potenziale repulsivo. La dimostrazione avviene in tre passaggi. Inizialmente viene 
effettuato un parallelo tra l'equazione alle derivate parziali e un modello discreto di una famiglia di particelle, governato da un 
sistema di ODE. In questa situazione semplificata si prova l'esistenza di una soluzione e si mostra che, senza ulteriori condizioni, 
il modello ammette infinite soluzioni. Viene poi trovata una soluzione approssimata nel caso generale in dimensione arbitraria, grazie 
allo schema dei movimenti minimizzanti. Nella parte finale, restringendosi al caso unidimensionale, si dimostra sia l'esistenza sia 
l'unicità della soluzione.

\clearpage

\pagestyle{fancy}

\tableofcontents

\thispagestyle{empty}

\fancyhf{}

\renewcommand{\thepage}{\arabic{page}} 
\setcounter{page}{1}

\lhead[\thepage]{\small{\textsc{Introduction}}}
\rhead[\small{\textsc{Introduction} }] {\thepage}		

\clearpage{\pagestyle{empty}\cleardoublepage}

\addcontentsline{toc}{chapter}{Introduction} 
\markboth{Introduction}{Introduction}

\chapter*{Introduction}

In the last century there has been great interest in the study of partial differential equations both from the theoretical and from the 
practical point of view, owing the realization that an enormous variety of phenomena, covering topics from physics to social 
sciences, can be described by these equations. The search for a solution is, of course, the first problem to be addressed. A solution, 
if it exists, can be found explicitly or, in many cases, estimated with numerical methods. As a consequence, developing methods 
for finding approximated solutions to PDE problems is of key importance, 
besides the formulation of very general theorems. In the present thesis we consider some interesting examples of a quite general class 
of PDE's modeling evolution phenomena which resemble the typical structure of "Gradient Flows". A very general theory has been developed
 in the last decades, covering a large class of applications in an increasing level of generality, starting from the prototypical examples 
of dynamical systems in finite dimensional spaces. \\

Over the last years, due to industrial application and to the evolution of the trends in theoretical physics, a lot of attention was 
given to the modelling of granular material (sand, powders, heaps of cereals, grains, molecules, snow, or even asteroids...). In this 
thesis is briefly introduced a kinetic model for granular flows. It is analyzed in the one-dimensional and spatially homogeneous case. 
Under these assumptions it is driven by the equation
 \begin{equation*}
 \partial_t f(x,v,t)=-\kappa \partial_v(Ff)+ \sigma \partial_v(vf) + T_e \partial_v^2f,
 \end{equation*}
  where $F=f \ast \nabla W$. The terms in the RHS of the equation describe collisions, drift and a heat bath.
 Here we focus on a model describing the evolution of a discrete or continuous 
distribution of masses $\mu(t)$ under the influence of a dissipative force generated by an interaction potential $W$, 
described by the equation  \begin{equation*}
 \frac{\partial \mu}{\partial t}=\nabla \cdot (\mu (\nabla W \ast \mu)) \qquad x \in \mathbb{R}^d \; , \; t > 0.
 \end{equation*} 
The particular structure, 
allowing for collapsing or diffusing phenomena, can be hardly settled in the usual setting of classical function spaces, but can be 
nicely attacked by using some tools of measure theory and optimal transportation which have been recently developed (\cite{CV,AGS}) 
and have influenced many interesting contributions. \\

In the first chapter we present a quick overview of the theory of optimal transportation. This is the natural setting for the formulation 
of our evolution problem as a suitable gradient flow and is also a nice and important theory, which has an old historical motivation 
(since the pioneering papers of Monge) and many deep connections with various problems in physics, finance, probability, geometry and 
analysis. 
First is introduced the class of transportation problems (arising from the original formulation proposed by Monge). Its solution 
provides the optimal transport map between two distributions, which realizes the minimal quadratic cost, measured in a suitable integral 
way. Then  is presented the relaxed formulation of the problem introduced by Kantorovich: it allows for the greatest generality on the
 measures and in particular it is well suited to deal with concentration phenomena. Thanks to this point of view, is defined a distance 
 (the Wasserstein distance) on the space of probability measure (nowadays called Wasserstein space). This is one of the main tools 
 which characterize the gradient flow theory. \\

In the second chapter we present some useful results concerning evolution PDE's in the space of measures, starting from the classical 
continuity equation. The concepts of 
curve of maximal slope and gradient flows are then introduced following the "metric" point of view inspired by E. De Giorgi and fully 
developed in \cite{AGS}. These concepts are shown to be equivalent under certain conditions. The minimizing 
movement scheme, a variational approximation algorithm particularly useful in the theory of gradient flows is then briefly described. \\

In the third chapter is introduced the one-dimensional kinect model for granular flows, explaining the connection with the theory of gradient flows. 
An important part of the chapter is devoted to the 
study of the evolution equation driven by an interaction potential. The current available results show that the problem is well 
posed when the potential satisfy some 
regularity assumptions. The work by Carrillo \emph{et al.} (see \cite{IL}) slightly weakens the above mentioned assumptions and provides a new method 
 applied to a particular case in this thesis. \\

In the fourth chapter we focus on a particular model, characterized by a repulsive potential which lacks differentiability near the origin. 
Since this kind of singularity is of concave type, this class of potentials cannot be attacked by using the existing results and an 
ad-hoc analysis had to be studied. Our contribution can be divided in three steps: \begin{itemize}
\item By comparing the PDE to an ODE system the existence of a solution is envisaged.
\item By means of the minimizing movement scheme a trial solution is provided.
\item The trial solution is proved to be exact when the problem is one dimensional. 
A sketch of a possible proof for the problem in many dimension is also presented.
\end{itemize}

\clearpage

\lhead[\thepage]{\small{\textsc{Optimal transportation}}}
\rhead[\small{\textsc{Optimal transportation} }] {\thepage}	

\clearpage{\pagestyle{empty}\cleardoublepage}

\chapter{Optimal transportation}
\markboth{Optimal transportation}{Optimal transportation}

\section{Historical prelude}

The problem has been studied firstly by Gaspard Monge (1746-1818) in 1781 in his famous work \emph{''M\'{e}moire sur la th\'{e}orie des d\'{e}blais 
et des remblais''}. We report the very intuitive description of the problem made in \cite{CV}: assume you have a certain amount of soil
to extract from the ground and transport to places where it should be incorporated in a construction (as in figure). \begin{figure} \centering 
\includegraphics[scale=1]{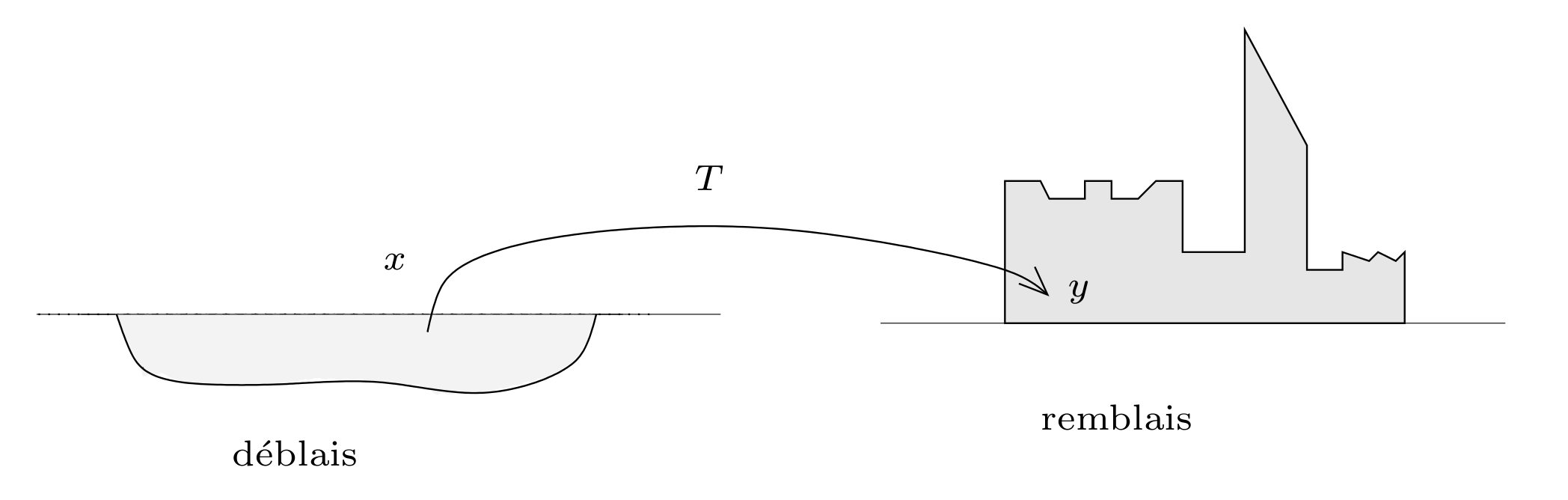}
\caption{Image taken from \cite{CV}}
\end{figure}  
Well known are the places where the material should be extracted, and the ones where it should be transported to. But the assignment has 
to be determined: to which destination should one send the material that has been extracted at a certain place? The answer does matter 
because transport is costly, and you want to minimize the total cost. Monge assumed that the transport cost of one unit of mass along a 
certain distance was given by the product of the mass by the distance. Mathematically it means that we are trying to minimize an integral like
\begin{equation*}
\int\left|T(x)-x\right|m(x)dx
\end{equation*}
Another interesting example allows us to see an economic perspective of the problem. Consider a large number of bakeries, producing loaves, 
that should be transported each morning to caf\'{e}s. The amount of bread that can be produced at each bakery, and the amount that will be consumed 
at each caf\'{e} are known in advance, and can be modeled as probability measures on a certain space. In our case we have a "density of production" 
and a "density of consumption" living in Paris equipped with the natural metric, given by the length of the shortest path between two point. 
The problem is to find in practice where each unit of bread should go in such to minimize the total transport cost. \\
Subsequently the problem, forgotten for many years (with the exception of some work about locational optimization during the nineteenth century), 
was studied again starting from the 1920s and 1930s, respectively in USSR and USA. This happened because there were found many connections 
between optimization problem, economy and war application. A deeper analysis of the historical point of view can be found in \cite{storiatrasporto}.
\\ A great breakthrough happened with the work of Leonid Kantorovich (1912-1986), that in 1942\footnote{Kantorovich like Monge published 
his work after many years because of national interests. Monge was a "warrior scientists" of the French Revolution, USSR used to keep many 
economics research secret} studied an economic problem deeply related with Monge's. In fact he discovered later this connection, but after 
that the problem of optimal transportation has been called the Monge-Kantorovich problem (the new formulation was very useful and 
circumvented the problem of indivisible masses). Kantorovich's work was remarkable and he was awarded with the Nobel Prize of economics, 
jointly with Tjalling Koopmans. His other main contributions were stating and proving a duality theorem (with Rubinstein) and the 
definition of a distance 
between probability measures. This distance (called Kantorovich-Rubinstein but nowadays renamed Wasserstein's distance) is one of the main 
tools of the transportation theory used in this thesis. \\ Other important contributions (see \cite{VILTO}) arrived from 
Ronald Dobrushin (study of particle 
systems), Hiroshi Tanaka (time-behavior of a simple variant of the Boltzmann equation), John Mather (Lagrangian dynamical system), 
Yann Brenier (incompressible fluid mechanics), Mike Cullen (semi-geostrophic equations for meteorology) and Mikhail Gromov (geometry). 
An exhaustive overview of transportation problem from a pobabilistic point of view can be found in the book of Rachev-Rushendorf.

\section{Discrete model of transportation}

\subsection*{Primal problem} 

Let us introduce here the simplest Kantorovich formulation of transportation problem at a finite discrete (but still useful in many
 economic questions) level. \\ We consider 
\begin{itemize} 
\item $m$ initial points of the configuration (indexed with $i=1,\cdots,m$);
\item $n$ final points (indexed with $j=1,\cdots,n$);
\item $c_{ij}$, the cost of the transport from $i$ to $j$;
\item $\mathbf{p}=(p_i)$, the initial distribution of masses;
\item $\mathbf{q}=(q_j)$, the final distribution;
\item $x_{ij} \ge 0$ is the quantity transported from $i$ to $j$ (unknown).
\end{itemize}

\begin{problem}[Primal]
Find $\mathbf{x}=(x_{ij})$ minimizing the cost $\mathcal{C}(\mathbf{x})=\sum_{i,j}c_{ij}x_{ij}$ with the constraints
\begin{itemize}
\item $\forall i,j \qquad p_i=\sum_jx_{ij}$ and $ q_j=\sum_ix_{ij}$, 
\item $\sum_ip_i=\sum_jq_j$.
\end{itemize}
\end{problem}
The two conditions mean that every mass is entirely transported and the total mass must be conserved.
Two important things must be noticed. The first is that no conditions on the cost are assumed; even if one usually imagines a positive 
and maybe linear cost, still negative (and also quite complicated) costs can be considered. 
And the second, easy to discover, is that adding a constant to the cost leaves the problem unchanged.

\subsection*{Dual problem}

This problem can be explained with an example. We can start thinking at the case with bakeries and caf\'{e}s. There is a delivery company that offer 
to do the transportation, payed $u_i$ for every unity of product taken from $p_i$ and $v_j$ for each delivery made to $q_j$. Obviously 
they offer a savings, so a condition must be imposed: \begin{equation*}
u_i+v_j\le c_{ij}
\end{equation*}
The delivery company wants to maximize their earnings, they aim to solve the new problem
\begin{problem}[Dual]
Find $\mathbf{w}=(u_i,v_j)$ maximizing $\mathcal{P}(\mathbf{w})=\sum_ip_iu_i+\sum_jq_jv_j$  under the constraint \begin{equation}
u_i+v_j\le c_{ij}.
\end{equation}
\end{problem}
Applying Von Neumann theorem, it would not be too difficult to show that the two problems are equivalent, i.e. 
\begin{equation*}
\text{inf}_{\mathbf{x}}\mathcal{C}(\mathbf{x})=\text{sup}_{\mathbf{w}}\mathcal{P}(\mathbf{w}),
\end{equation*}
where obviously $\mathbf{x}$ and $\mathbf{w}$ have to satisfy the constraints imposed by the respective problems.

\section{Notation and measure-theoretic results}

\subsection*{Probability measures}

Given a separable metric space $(X,d)$, we denote by $\mathcal{B}(X)$ be the Borel $\sigma$-algebra on $X$, i.e. is the $\sigma$-algebra 
generated by the open sets of $X$. A Borel probability measure is a function $\mu:\mathcal{B}(X) \to [0,1] \ \sigma$-additive. \\
We denote by $\mathcal{P}(X)$] the set of Borel probability measure $\mu:\mathcal{B}(X)\to [0,1]$ on  $X$. \\
The support of $\mu \in \mathcal{P}(X)$ is the closed set \begin{equation}
\text{supp}(\mu):=\bigl\{x \in X:\mu \bigl(B_r(x)\bigr)>0 \ \forall\text{m} r > 0\bigr\}.
\end{equation}
When $X$ is a Borel subset of an euclidean space $\mathbb{R}^d$, we set \begin{equation}
\text{m}_2(\mu):=\int_X|x|^2d\mu,
\end{equation}
and we can make the identification \begin{equation}
\mathcal{P}(X)=\bigl\{\mu \in \mathcal{P}(\mathbb{R}^d):\mu(\mathbb{R}^d \setminus X)=0\bigl\}.
\end{equation}
and we denote by $\mathcal{P}_2(X)$ the subspace of $\mathcal{P}(X)$ made by measures with finite
quadratic moment:
\begin{equation}
\mathcal{P}_2(X):=\bigl\{\mu \in \mathcal{P}(X):\text{m}_2(\mu)<\infty\bigr\}.
\end{equation}
We denote by $\mathcal{L}^d$ the Lebesgue measure in $\mathbb{R}^d$ and set 
\begin{equation}
\mathcal{P}^a_2(X):=\bigl\{\mu \in \mathcal{P}_2(X): \mu \ll \mathcal{L}^d\bigr\},
\end{equation}
whenever $X \in \mathcal{B}(\mathbb{R}^d)$.

\subsection*{Transport maps and transport plans}

\begin{definition}[Push-forward] If $\mu \in \mathcal{P}(X_1)$, and $\mathbf{r}:X_1\to X_2$ is a Borel map, we denote by $\mathbf{r}_\# \mu \in \mathcal{P}(X_2)$ 
the push-forward of $\mu$ through $\mathbf{r}$, defined by \begin{equation}
\mathbf{r}_\# \mu (B) :=\mu\bigr(\mathbf{r}^{-1}(B)\bigr) \ \forall B \in \mathcal{B}(X_2).
\end{equation}
\end{definition}
From that definition we can get the change-of-variable formula, that holds for every Borel map $f:X_2\to\mathbb{R}$: \begin{equation}
\label{changevar}
\int_{X_1}f\bigl(\mathbf{r}(x)\bigr)d\mu(x)=\int_{X_2}f(y)d\mathbf{r}_\# \mu(y)
\end{equation}
Useful are \begin{definition} $\pi^i$ are the projection operators, with $i=1,2$, working on the product space $\mathbf{X}:=X_1 \times X_2$, defined 
by \begin{equation}
\pi^1:(x_1,x_2)\mapsto x_1 \in X_1, \ \ \ \pi^2:(x_1,x_2)\mapsto x_2 \in X_2. 
\end{equation}
If $\mathbf{X}$ is endowed with the canonical product metric and the Borel $\sigma$-algebra and $\boldsymbol{\mu}\in\mathcal{P}(X)$, we define 
marginals of $\boldsymbol{\mu}$. They are the probability measures \begin{equation}
\mu^i;=\pi^i_\#\boldsymbol{\mu}\in\mathcal{P}(X_i), \ i=1,2.
\end{equation}
At last we need to define what is a transport plan: given $\mu^1\in \mathcal{P}(X_1)$ and $\mu^2\in \mathcal{P}(X_2)$, the class $\Gamma(\mu^1,\mu^2)$ 
of transport plans between $\mu^1$ and $\mu^2$ is defined by \begin{equation}
\Gamma(\mu^1,\mu^2):=\bigl\{\boldsymbol{\mu}\in \mathcal{P}(X_1\times X_2):\pi^i_\#\boldsymbol{\mu}=\mu^i,i=1,2\bigr\}.
\end{equation}
\end{definition}
It is important to enounce the following useful theorem: \begin{theorem}
Let $X_1,X_2,X_3$ be complete separable metric spaces and let $\boldsymbol{\gamma}^{12} \in \mathcal{P}(X_1 \times X_2), \ \boldsymbol{\gamma}^{13} \in 
\mathcal{P}(X_1\times X_3)$ such that $\pi^1_\#\boldsymbol{\gamma}^{12}=\pi^1_\#\boldsymbol{\gamma}^{13}=\mu^1$. Then there exists $\boldsymbol{\mu} 
\in \mathcal{P}(X_1 \times X_2 \times X_3)$ such that \begin{equation}
\pi^{1,2}_\#\boldsymbol{\mu}=\boldsymbol{\gamma}^{12}, \qquad \pi^{1,3}_\#\boldsymbol{\mu}=\boldsymbol{\gamma}^{13}.
\end{equation}
\label{JL}
\end{theorem}

\subsection*{Narrow convergence}

\begin{definition}[Narrow convergence]A sequence $(\mu_n)\subset \mathcal{P}(X)$ is narrowly convergent to $\mu \in \mathcal{P}(X)$ as $n \to \infty$ if
\begin{equation}
\lim_{n\to\infty}\int_Xf(x)d\mu_n(x)=\int_Xf(x)d\mu(x)
\label{narrowlsc}
\end{equation}
for every function $f\in C_b^0(X)$.
\end{definition}
An important result is this theorem:
\begin{theorem}[Prokhorov]If a set $\mathcal{K}\subset\mathcal{P}(X)$ is tight, i.e. \label{prok} 
\begin{equation}
\forall \epsilon > 0 \ \exists K_\epsilon \text{ compact in X such that } \mu (X \setminus K_\epsilon) \le \epsilon \ \forall \mu \in \mathcal{K},
\end{equation}
then $\mathcal{K}$ is relatively compact\footnote{$X$ is relative compact in $Y$ if his closure is a compact subset of $Y$} in $\mathcal{P}(X)$. 
Conversely, every relatively compact subset of $\mathcal{P}(X)$ is tight.
\end{theorem}

When one needs to pass to the limit in expressions like \ref{narrowlsc} w.r.t lower semicontinuous function $f$, the following 
property is quite useful:
\begin{proposition}
Given a sequence $(\mu_n) \subset \mathcal{P}(X)$ narrowly convergent to $\mu$ and a l.s.c. function $g:X\to (-\infty,+\infty]$ bounded 
from below we have that \begin{equation}
\liminf_{n \to \infty}\int_{X}g(x)d\mu_n(x)\geq \int_Xg(x)d\mu(x).
\end{equation}
\end{proposition}

\section{Formulation of the Kantorovich problem}

Let $X,Y$ be a complete and separable metric spaces and let $c:X\times Y \to [0,+\infty]$ be a Borel cost function. Given $\mu \in 
\mathcal{P}(X)$ and $\nu \in \mathcal{P}(Y)$, the optimal transport problem, in Monge's formulation, is given by 
\begin{equation}
\text{inf}\biggl\{\int_Xc\left(x,\mathbf{t}(x)\right)d\mu(x):\mathbf{t}_\# \mu = \nu \biggr\}.
\label{prob1}
\end{equation}
This problem can be ill posed because sometimes there is no transport map $\mathbf{t}$ such that $\mathbf{t}_\# \mu$ (this happens for instance when 
$\mu$ is a Dirac mass and $\nu$ is not a Dirac mass). In Kantorovich's formulation the problem become
\begin{equation}
\text{min}\biggl\{\int \int_{X \times Y}c(x,y)d \boldsymbol{\gamma}(x,y):\boldsymbol{\gamma} \in \Gamma(\mu,\nu)\biggr\}.
\label{prob2}
\end{equation}
\begin{definition}
$\Gamma_0(\mu,\nu)$ is the space of the optimal plans that realizes the minimum in (\ref{prob2}). 
\end{definition}

The dual problem has the same structure as in the discrete case, 
\begin{equation} \text{sup}\biggl\{ \int_X\phi(x)d\mu(x) \ +\int_Y\psi(x)d\nu(y) \biggr\} \end{equation}
with the condition
\begin{equation} (\phi,\psi)\in C_b^0(X)\times C_b^0(Y) \; , \; \phi(x)+\psi(y) \leq c(x,y) \end{equation}
\begin{remark}
The equivalence between the two formulations of the problem is not always true, but has been demonstrated for very general cost function. 
The most recent works are the one of Pratelli \cite{Prat}, 
Schachermayer and Teichmann \cite{SCTE} and Beiglbock \emph{et al.} \cite{AAAA}.
\end{remark} 
Detailed demonstrations about the existence can be found in \cite{AGS} and in \cite{CV}. Here is reported 
the demonstration for a general lower-semicontinuous cost function. It is needed for the subsequent results.
\begin{theorem}
Suppose that $ c:X\times Y \to [0,\infty]$ is a lower semicontinuous function. If there exists $\boldsymbol{\gamma} \in \Gamma(\mu,\nu)$ 
with a finite cost $\mathcal{C}(\boldsymbol{\gamma})$, then the problem admits minimum. 
\proof $\Gamma(\mu,\nu)$ is not empty. We need to prove that $\Gamma(\mu,\nu)$ 
is compact. Fixed $\delta$, there are two compact subset $K\subset X,L\subset Y$ such that \begin{equation*}
\mu[X\setminus K]\leq \delta \ , \ \nu[Y\setminus L] \leq \delta . 
\end{equation*} 
For every $\boldsymbol{\gamma} \in \Gamma(\mu,\nu)$ we have that \begin{equation*}
\boldsymbol{\gamma}[(X \times Y) \setminus (K \times L)] \leq \boldsymbol{\gamma}[(X \times (Y \setminus L)] + \boldsymbol{\gamma}[(X  \setminus K ) \times Y)] \leq 2\delta
\end{equation*}
We get that $\Gamma(\mu,\nu)$ is tight and so, thanks to the Prokhorov theorem (\ref{prok}) it is relatively compact. The 
conditions defining $\Gamma(\mu,\nu)$ are continuous with respect to the narrow topology, so it is weakly closed and so compact. Then, 
taking a minimizing sequence ${\boldsymbol{\gamma}_k}$, it converge (up to a subsequence) to $\gamma \in \Gamma(\mu,\nu)$. The cost function can be viewed 
as the supremum of an increasing sequence of limited functions $c_i$, thanks to its lower semi-continuity. Using the monotone convergence we have \begin{eqnarray*} 
 \int c(x,y)d\boldsymbol{\gamma} & = & \lim_{i \to \infty} \int c_i(x,y)d\boldsymbol{\gamma} \ \leq \ \lim_{i \to \infty} \limsup_{k \to \infty} \int c_i(x,y)d\boldsymbol{\gamma}_k \\
 & \leq & \limsup_{k \to \infty} \int c(x,y)d\boldsymbol{\gamma}_k \ = \ \inf_{\Gamma(\mu,\nu)}\mathcal{C}(\boldsymbol{\gamma})
\end{eqnarray*} Now we have that $\boldsymbol{\gamma}$ realize the minimum.
\endproof
\label{existence}
\end{theorem} 

\section{Wasserstein distance}

Now the two spaces are taken to be $\mathbb{R}^d$ and the cost function is the squared distance between points.
\begin{definition}
The \emph{Kantorovich-Rubinstein-Wasserstein distance} between two probability measure is defined to be \begin{equation}
W_2(\mu,\nu):=\left(\int_{X\times X}\text{d}^2(x,y)d\boldsymbol{\gamma}\right)^{1 / 2} \; , \; \boldsymbol{\gamma}\in \Gamma_0(\mu,\nu)
\end{equation}
\end{definition}
The existence of the minimum has been demonstrated in theorem \ref{existence}, but the two hypothesis must be satisfied. The l.s.c. of the cost 
function is obviously satisfied. The cost $\mathcal{C}(\boldsymbol{\gamma})$ may be infinite and so some type of constraints must be found. 
A possible request is to take probability measure with a compact support, included in a ball with finite radius. This request is too strong for 
our purpose. Taking $\boldsymbol{\gamma}=\mu\otimes\nu$ we have: \begin{equation*}
\int\text{d}^2(x,y)d\mu\otimes\nu=\int|x|^2d\mu\otimes\nu + \cdots + \int|y|^2d\mu\otimes\nu ,
\end{equation*}
and this is finite when $\mu$ and $\nu$ have finite quadratic moment. For this reason probability measure are searched in the space 
$\mathcal{P}_2(\mathbb{R}^d)$. Before proving that it is a distance we show two short examples and an existence theorem.
\begin{example}
Take $\mu=\delta_x$ and $\nu=\delta_y$. There is only one optimal plan that is $\boldsymbol{\gamma}=\delta_{(x,y)}$. An easy calculation gives 
\begin{equation*}
W_2(\delta_x,\delta_y)=|x-y|
\end{equation*}
\end{example}
\begin{theorem}[Brenier, see \cite{Bren}]
For any $\mu \in \mathcal{P}_2^a(\mathbb{R}^d)$, $\nu \in \mathcal{P}_2(\mathbb{R}^d)$ Kantorovich's optimal transport problem
(\ref{prob2}) with $c(x, y) = |x - y|^2$ has a unique solution $\boldsymbol{\gamma}$. Moreover: \begin{itemize}
\item $\nu = r_{\#}\mu$;
\item $\boldsymbol{\gamma}=(id\otimes r)_{\#}\mu$ is the unique optimal plan; 
\item $r$ is the gradient of a convex function.
\end{itemize}
\end{theorem}
\begin{example}
Thanks to the previous theorem, if $\mu \ll \mathcal{L}$ there exists a unique optimal map $r$ such that $r_\#(\mu)=\nu$. Because of this we have that  
\begin{equation*}
W_2(\mu,\nu)=\left( \int_{\mathbb{R}^d}|r(x)-x|^2d\mu(x)\right)^{1 / 2} 
\end{equation*}  
\end{example}
\begin{theorem}
$W_2$ defines a distance in $\mathcal{P}_2(\mathbb{R}^d)$.
\proof
Let $\mu,\nu,\sigma \in \mathcal{P}_2(\mathbb{R}^d)$ and let $\boldsymbol{\gamma}\in \Gamma_0(\mu,\nu)$ and $\boldsymbol{\eta} \in \Gamma_0(\nu,\sigma)$. 
The theorem \ref{JL} ensure the existence of $\boldsymbol{\lambda} \in \mathcal{P}(\mathbb{R}^d \times \mathbb{R}^d \times \mathbb{R}^d)$ such that \begin{equation*}
(\pi^1,\pi^2)_\#\boldsymbol{\lambda}=\boldsymbol{\gamma}, \qquad (\pi^2,\pi^3)_\#\boldsymbol{\lambda}=\boldsymbol{\eta}.
\end{equation*}
Then, as \begin{equation*}
\pi^1_\#(\pi^1,\pi^3)_\#\boldsymbol{\lambda}=\pi^1_\#\boldsymbol{\lambda}=\pi^1_\#\boldsymbol{\gamma}=\mu, \qquad 
\pi^2_\#(\pi^1,\pi^3)_\#\boldsymbol{\lambda}=\pi^3_\#\boldsymbol{\lambda}=\pi^2_\#\boldsymbol{\gamma}=\sigma,
\end{equation*}
we obtain that $(\pi^1,\pi^3)_\#\boldsymbol{\lambda}\in \Gamma(\mu,\sigma)$, hence \begin{equation*}
W_p(\mu,\sigma)\leq \left( \int_{X \times X}\text{d}^p(x_1,x_3)d(\pi^1,\pi^3)_\#\boldsymbol{\lambda}\right)^{1 / p} = \bigl\|\text{d}(x_1,x_3)\bigr\|_{L^p(\boldsymbol{\lambda})}.
\end{equation*}
As $\text{d}(x_1,x_3)\leq \text{d}(x_1,x_2) + \text{d}(x_2,x_3)$ and  \begin{eqnarray*}
  \bigl\|\text{d}(x_1,x_2)\bigr\|_{L^p(\boldsymbol{\lambda})} & = & \bigl\|\text{d}(x_1,x_2)\bigr\|_{L^p(\boldsymbol{\gamma})} = \; W_p(\mu,\nu), \\
  \bigl\|\text{d}(x_2,x_3)\bigr\|_{L^p(\boldsymbol{\lambda})} & = & \bigl\|\text{d}(x_2,x_3)\bigr\|_{L^p(\boldsymbol{\eta})} =  \; W_p(\nu,\sigma),
\end{eqnarray*} 
the triangle inequality $W_p(\mu,\sigma)\leq W_p(\mu,\nu) + W_p(\nu,\sigma)$ follows by the standard triangle inequality in $L^p(\boldsymbol{\lambda})$.
\endproof
\end{theorem}

The next result shows a characterization of the convergence induced by the distance $W_2$.
\begin{theorem}
$\mathcal{P}_2(\mathbb{R}^d)$, endowed with the Wasserstein distance, is a complete and separable metric space. A set 
$\mathcal{K} \subset \mathcal{P}_2(\mathbb{R}^d)$ 
is relatively compact iff it is 2-uniformly integrable and tight. Furthermore, for a given sequence 
$(\mu_n) \subset \mathcal{P}_2(\mathbb{R}^d)$ the following implication holds: 
\begin{equation*}
\lim_{n \to \infty}W_2(\mu_n,\mu)=0 
\end{equation*} 
\begin{equation}
\Updownarrow 
\end{equation}
\begin{equation*}
\begin{cases}
\mu_n \text{ narrowly converge to } \mu, \\ (\mu_n) \text{ has uniformly integrable quadratic moments}
\end{cases}
\end{equation*} 
\end{theorem}
The next short example shows why the narrow convergence is different from the convergence in the Wasserstein metric.
\begin{example}
Consider a sequence $(x_n) \in \mathbb{R}$ such that $\lim_{n \to \infty} x_n=\infty$ and the sequence of probability measures $\mu_n=\alpha_n \delta_{x_n} +
(1-\alpha_n) \delta_{x_0}$. If $(\alpha_n)$ is infinitesimal it's easy to see that $\mu_n \to \delta_{x_0}$ narrowly. Calculating 
the distance between $\delta_{x_0}$ and $\mu_n$ we get $W_2(\delta_{x_0},\mu_n)=\sqrt{\alpha_n}|x_n-x_0|$, so the convergence in 
$\mathcal{P}_2(\mathbb{R})$ requires the stronger condition that $\alpha_n|x_n-x_0|^2 \to 0$.
\end{example}

\subsection*{The real case}
\label{altrasubse}
In this subsection is introduced an useful change of variable that holds only in the one dimensional case. The result obtained here 
permits to simplify the calculus for the Wasserstein distance and for more general integral, as will be done in the third chapter. 
Starting from a probability measure $\mu$ on $\mathbb{R}$ we can define its cumulative distribution:
\begin{definition}
The cumulative distribution function $M_{\mu}$ associated to a probability measure $\mu \in \mathcal{P}(\mathbb{R})$ is \begin{equation}
M_{\mu}(x):=\mu((-\infty,x])=\int_{-\infty}^xd\mu(x).
\end{equation}
\end{definition}
Every probability measure $\mu \in \mathcal{P}(\mathbb{R})$ can be represented by its monotone rearrangement $X_{\mu}$ which is the 
pseudo-inverse of the the distribution function $M_{\mu}$.
\begin{definition} \label{defimonoreaaq}
The monotone rearrangement $X_{\mu}$ is defined by \begin{equation}
\label{pseudoinve}
X_{\mu}(s):=\inf \{ x: M_{\mu}(x)>s \} \qquad s \in (0,1).
\end{equation}
\end{definition}  
The map $\mu \mapsto X_{\mu}$ is an isometry between $\mathcal{P}_2(\mathbb{R})$ (endowed with the Wasserstein distance) and the convex cone 
$\mathcal{K}$ of non decreasing functions in the Hilbert space $L^2(0,1)$ and the following theorem holds: 
\begin{theorem}
Given a probability measure $\mu \in \mathcal{P}(\mathbb{R})$ and its pseudo-inverse $X_{\mu}$ defined as in (\ref{defimonoreaaq}) we have 
that \begin{equation}
\int_{\mathbb{R}} \xi(x) d\mu(x)= \int_0^1\xi(X_{\mu}(s))ds
\end{equation}
for every  nonnegative Borel map $\xi : R \to [0,+\infty]$.
\end{theorem}
Thanks to previous result  (see \cite{SGANGB,LNS}) the Wasserstein distance can be rewritten in an easier form thanks to the 
change of variable formula (\ref{changevar}) and to the pseudo-inverse (\ref{pseudoinve}):
\begin{equation} \label{changevarutil}
\begin{split}
W_2^2(\mu^1,\mu^2)=\int_0^1 \bigl|X_{\mu^1}(s)-X_{\mu^2}(s)\bigr|^2ds.
\end{split}
\end{equation}

\begin{figure}[p]
\centering
\subfloat[][\emph{Probability measure}.]
{\includegraphics[width=.60\columnwidth]{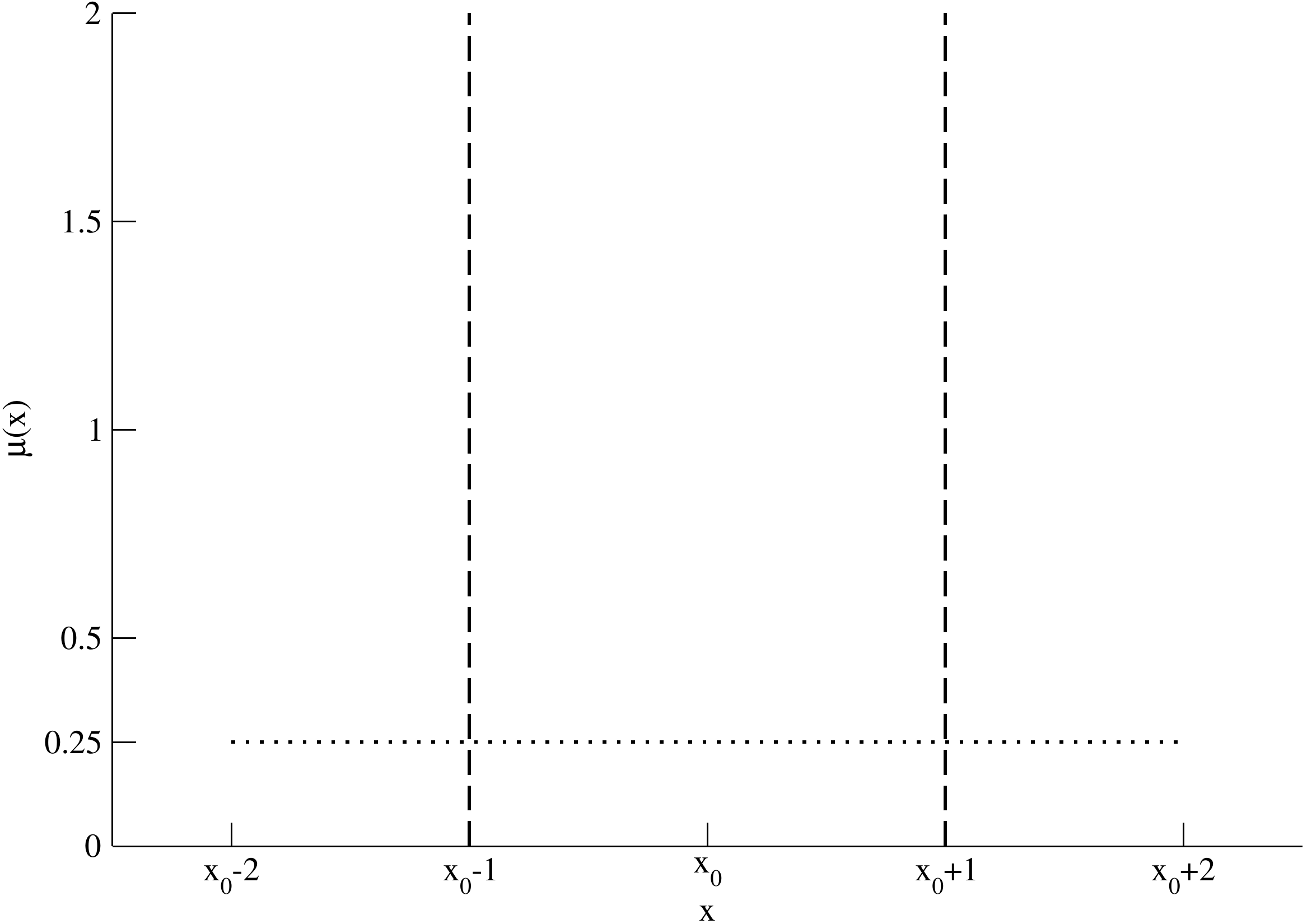}} \quad
\subfloat[][\emph{Cumulative distribution function}.]
{\includegraphics[width=.60\columnwidth]{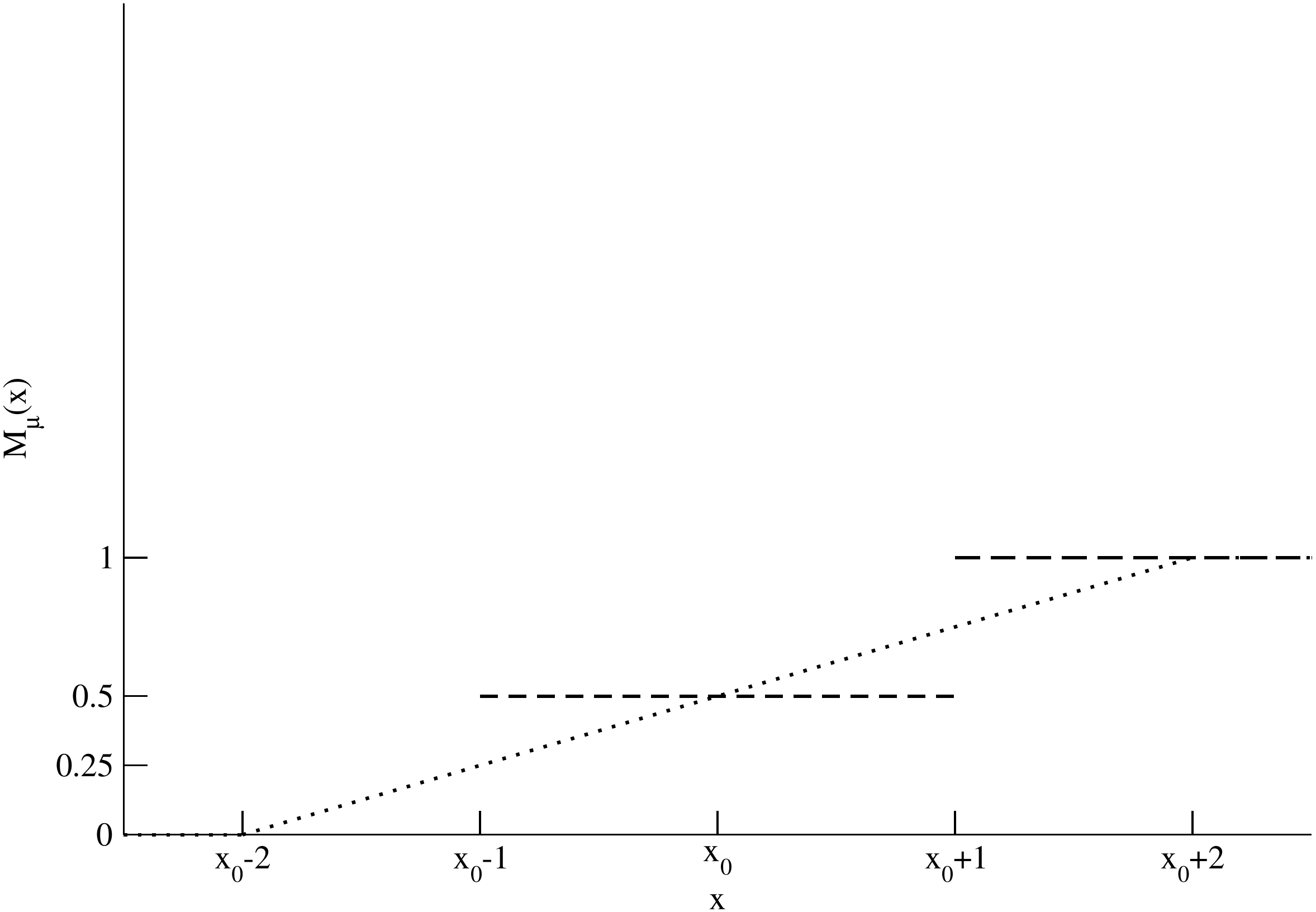}} \\
\subfloat[][\emph{Monotone rearrangement}.]
{\includegraphics[width=.60\columnwidth]{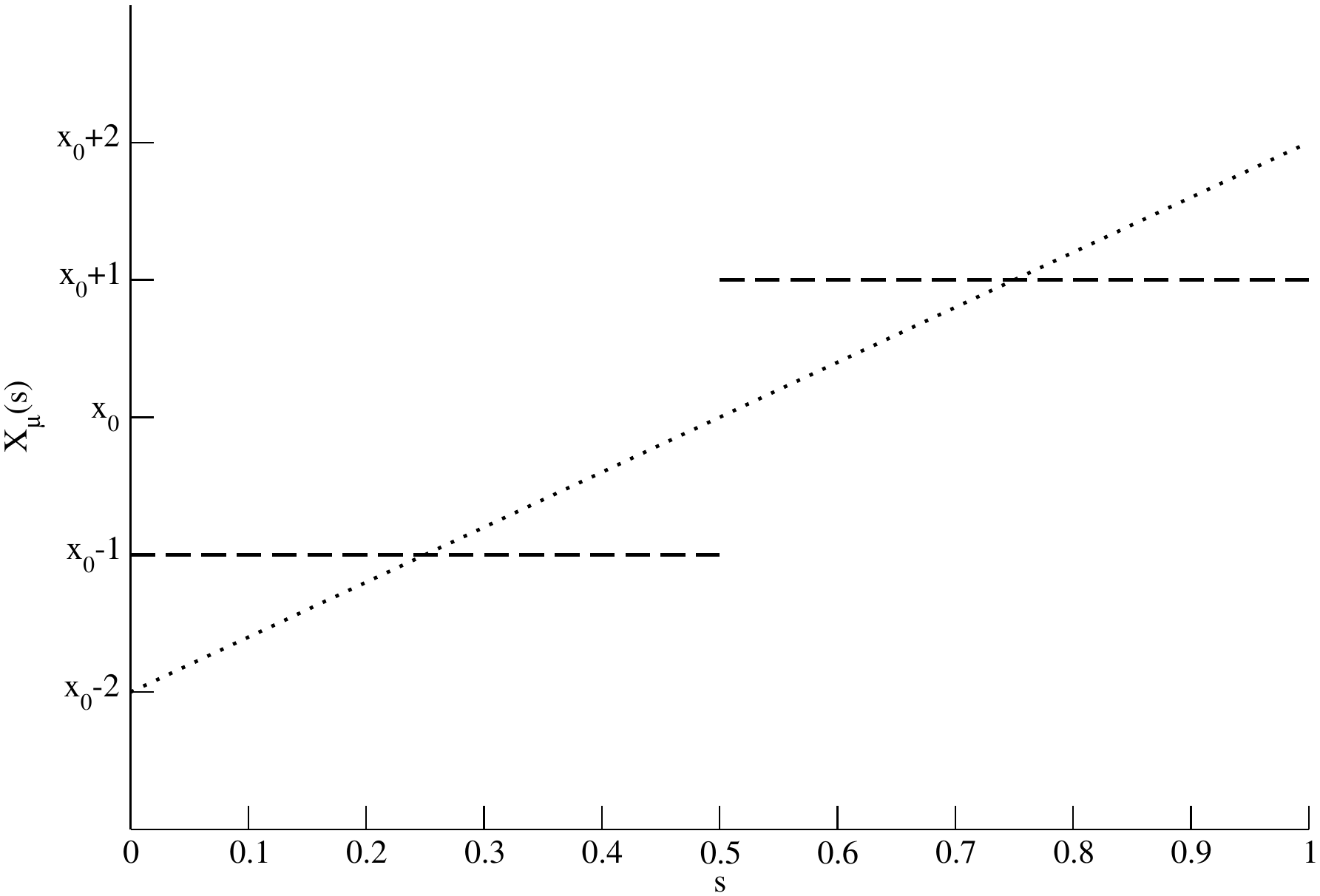}}
\caption{In the images above are represented two probability measures (the representative of the interval $[x_0-2,x_0+2]$ (dotted line), two 
dirac delta (dashed line)) and their respective cumulative distribution functions and monotone 
rearrangement.}
\end{figure}

\clearpage

\lhead[\thepage]{\small{\textsc{Continuity equation and gradient flows in $\mathcal{P}_2(\mathbb{R}^d)$}}}
\rhead[\small{\textsc{Continuity equation and gradient flows in $\mathcal{P}_2(\mathbb{R}^d)$} }] {\thepage}	

\clearpage{\pagestyle{empty}\cleardoublepage}

\chapter{Continuity equation and gradient flows in $\mathcal{P}_2(\mathbb{R}^d)$}
\markboth{Continuity equation and gradient flows in $\mathcal{P}_2(\mathbb{R}^d)$}{Continuity equation and gradient flows in $\mathcal{P}_2(\mathbb{R}^d)$}

\section{Continuity equation}

The continuity equation \begin{equation}
\label{ceq}
\partial_t\mu_t+\nabla \cdot (\boldsymbol{v}_t\mu_t)=0 \text{ in }\mathbb{R}^d\times(0,T),
\end{equation} 
is a differential equation that describes the transport of some kind of conserved quantity along the characteristic curves tangent to 
the vector field $\boldsymbol{v}_t$. The variety of phenomena described depends on the definition of the velocity field. 
Applications in physics as transport of mass, electric charge or momentum,  in population dynamics or in biology (as in chemotaxis 
models) also fit in this framework as we will briefly show in the following. \\ 

Let us briefly recall a few results concerning the continuity equation.

\begin{definition}
A family of  Borel probability measures $\mu_t$ on $\mathbb{R}^d$ defined for $t$ 
in the open interval $I:=(0,T)$ with $\boldsymbol{v}:(x,t)\mapsto \boldsymbol{v}_t(x) \in \mathbb{R}^d$ a Borel velocity field such that 
\begin{equation} \label{condv}
\int_0^T\int_{\mathbb{R}^d}|\boldsymbol{v}_t(x)|d\mu(x)dt < +\infty,
\end{equation}
solves the equation \ref{ceq} in the sense of distribution if \begin{equation}\label{distsense}
\begin{split}
\int_0^T\int_{\mathbb{R}^d}  \Bigl( \partial_t  \varphi(x,& t)+ \langle \boldsymbol{v}_t(x),\nabla_x \varphi(x,t) \rangle \Bigr)  d\mu(x)dt = 0, \\
& \forall \varphi \in C_c^\infty (\mathbb{R}^d \times (0,T)).
\end{split}
\end{equation}
\end{definition}

By a simple regularization argument via convolution, it is easy to show that the last equation holds if $\varphi \in C_c^1(\mathbb{R}^d \times (0,T))$ 
as well. An equivalent formulation of (\ref{distsense}) is \begin{equation} \label{distsense2}
\frac{d}{dt}\int_{\mathbb{R}^d} \zeta (x)d \mu_t(x)=\int_{\mathbb{R}^d} \langle \nabla \zeta(x),\boldsymbol{v}_t(x) \rangle d\mu_t(x) 
\qquad \forall \zeta \in C_c^{\infty} (\mathbb{R}^d) 
\end{equation}
in the sense of distributions in $(0,T)$; it corresponds to the choice $\varphi(t,x):=\eta(t)\zeta(x)$ with $\eta \in C_c^{\infty}(0,T)$. 
Using this formulation the following result is obtained:
\begin{lemma}[Continuous representative] Let $\mu_t$ be a Borel family of probability measures satisfying (\ref{distsense2}) for a 
Borel velocity field $\boldsymbol{v}_t$ satisfying (\ref{condv}). Then there exists a narrowly continuous curve $t\in[0,T]\mapsto 
\tilde{\mu}_t \in \mathcal{P}(\mathbb{R}^d)$ such that $\mu_t=\tilde{\mu}_t$ for $\mathscr{L}^1$-a.e. $t \in (0,T)$. \\ Moreover, if 
$\varphi \in C_c^1(\mathbb{R}^d \times [0,T])$ and $t_1 \geq t_2 \in [0,T]$ it's true that \begin{equation}
\begin{split}
\int_{\mathbb{R}^d} \varphi (x,t_2) d\tilde{\mu}_{t_2}(x)-\int_{\mathbb{R}^d}\varphi (x,t_1)d\tilde{\mu}_{t_1}(x) \\
 = \int_{t_1}^{t_2}\int_{\mathbb{R}^d}\left( \partial_t \varphi + \langle \nabla \varphi , \boldsymbol{v}_t \rangle \right) d\mu_t(x)dt
\end{split}
\end{equation}
\label{contire}
\end{lemma}
A more explicit formula representing solutions of (\ref{ceq}) can be found if the velocity field exhibits some regularity properties:
it can be obtained by the classical method of characteristics, which also provides existence and uniqueness for the 
solution of the continuity equation. 
\begin{lemma}[The characteristics system of ODE] Let $\boldsymbol{v}_t$ be a Borel vector field such that for every compact set 
$B \subset \mathbb{R}^d$ \begin{equation}
\label{condvtlip}
\int_0^T\left( \sup_B |\boldsymbol{v}_t| + \text{Lip}(\boldsymbol{v}_t,B) \right) dt < +\infty
\end{equation}
Then, for every $x \in \mathbb{R}^d$ and $s \in [0,T]$, the ODE \begin{equation}\label{charode}
X_s(x,s)=x, \qquad \frac{d}{dt}X_t(x,s)=\boldsymbol{v}_t\left( X_t(x,s) \right),
\end{equation}
admits a unique maximal solution defined in an interval $I(x,s)$ relatively open in $[0,T]$ and containing $s$ as (relatively) internal 
point. \\ Furthermore, if $t\mapsto |X_t(x,s,)|$ is bounded in the interior of $I(x,s)$ then $I(x,s)=[0,T]$; finally, if $\boldsymbol{v}$ 
satisfies the global bounds analogous \begin{equation}
S:=\int_0^T\left( \sup_{\mathbb{R}^d} |\boldsymbol{v}_t| + \text{Lip}(\boldsymbol{v}_t,\mathbb{R}^d) \right) dt < +\infty
\end{equation}
then the flow map $X$ satisfies \begin{equation}
\int_0^T \sup_{x \in \mathbb{R}^d}|\partial_tX_t(x,s)|dt\leq S, \qquad \sup_{t,s \in [0,T]}\text{Lip}\left( X_t(\cdot ,s),\mathbb{R}^d \right) \leq e^S.
\end{equation}
\end{lemma}

\begin{proposition}[Representation formula for the continuity equation] \label{repforfortheconequ}
Let $\mu_t$, $t\in [0,T]$, be a narrowly continuous family of Borel probability measures solving the continuity equation (\ref{ceq}) w.r.t. 
a Borel vector field $\boldsymbol{v}_t$ satisfying (\ref{condvtlip}) and (\ref{condv}). Then for $\mu_0$-a.e. $x \in \mathbb{R}^d$ the 
characteristic system (\ref{charode}) admits a globally defined solution $X(t)$ in $[0,T]$ and \begin{equation}
\mu_t=(X_t)_{\#}\mu_0 \; \forall t \in [0,T],
\end{equation}
moreover, if \begin{equation} \label{pintegr}
\int_0^T\int_{\mathbb{R}^d}|\boldsymbol{v}_t(x)|^pd\mu_t(x)dt < +\infty \;\text{ for some }p > 1,
\end{equation}
then the velocity field $\boldsymbol{v}_t$ is the time derivative of $X_t$ in the $L^p$-sense \begin{equation}
\lim_{h\downarrow 0}\int_0^{T-h}\int_{\mathbb{R}^d}\left| \frac{X_{t+h}(x)-X_t(x)}{h}-\boldsymbol{v}_t(X_t(x)) \right|^p d\mu_o(x)dt=0
\end{equation}
\begin{equation}
\lim_{h \to 0}\frac{X_{t+h}(x,t)-x}{h}=\boldsymbol{v}_t(x) \text{ in } L^p(\mu_t;\mathbb{R}^d) \text{ for } \mathscr{L}^1\text{-a.e. }t\in(0,T). 
\end{equation}
\end{proposition}

The following result shows a useful approximation technique which will be used in the proof of theorem \ref{exicurmaxslo}.
\begin{lemma}[Approximation by regular curves] \label{appbyregcur}
Let $p\geq 1$ and let $\mu_t$ be a time-continuous solution of (\ref{ceq}) w.r.t. a velocity field satisfying the $p$-integrability 
condition (\ref{pintegr}). Let $(\rho_{\varepsilon})\subset C^{\infty}(\mathbb{R}^d)$ be a family of strictly positive mollifiers in the $x$ 
variable, (e.g. $\rho_{\varepsilon}(x)=(2\pi\varepsilon)^{-d / 2}exp(-|x|^2/2\varepsilon$), and set \begin{equation}
\mu_t^{\varepsilon}:=\mu_t \ast \rho_{\varepsilon}, \; E_t^{\varepsilon}:=(\boldsymbol{v}_t\mu_t)\ast \rho_{\varepsilon}, \; 
\boldsymbol{v}_t^{\varepsilon} := \frac{E_t^{\varepsilon}}{\mu_t^{\varepsilon}}
\end{equation}
Then $\mu_t^{\varepsilon}$ is a continuous solution of (\ref{ceq}) w.r.t. $\boldsymbol{v}_t^{\varepsilon}$, which satisfy the local regularity 
assumptions (\ref{pintegr}) and the uniform integrability bounds \begin{equation}
\int_{\mathbb{R}^d}|\boldsymbol{v}_t^{\varepsilon}(x)|^pd\mu_t^{\varepsilon}(x)
 \leq \int_{\mathbb{R}^d}|\boldsymbol{v}_t(x)|^pd\mu_t(x) \qquad \forall t \in (0,T).
\end{equation}
Moreover, $E_t^{\varepsilon}\to \boldsymbol{v}_t\mu_t$ narrowly and \begin{equation}
\lim_{\varepsilon \downarrow 0}\|\boldsymbol{v}_t^{\varepsilon}\|_{L^p(\mu_t^{\varepsilon};\mathbb{R}^d)}=
\|\boldsymbol{v}_t\|_{L^p(\mu_t;\mathbb{R}^d)} \qquad \forall t \in (0,T).
\end{equation}
\end{lemma}

\subsection*{Continuity equation and absolutely continuous curve in $\mathcal{P}_2(\mathbb{R}^d)$}

Here we define absolutely continuous curve and their metric derivative. A following theorem will show that this class of curves coincides 
with (distributional) solutions of the continuity equation.

\begin{definition}[Absolutely continuous curve]\label{ascurvedef}
Thanks to the fact that $(\mathcal{P}_2(\mathbb{R}^d),W_2)$ is a complete metric space and letting $v:(a,b)\to \mathcal{P}_2(\mathbb{R}^d)$ 
be a curve; $v$ belongs to $AC^2(a,b;\mathcal{P}_2(\mathbb{R}^d))$ if there exists $m \in L^2(a,b)$, such that \begin{equation}
\text{d}(v(s),v(t))\leq \int_s^t m(r)d(r) \qquad \forall a < s \leq t < b.
\label{ascurve}
\end{equation} 
\end{definition}

Any curve in $AC^2 (a, b; \mathcal{P}_2(\mathbb{R}^d) )$ is uniformly continuous; if $a > - \infty $ (resp. $b < +\infty$) denoting by 
$v(a+)$ (resp. $v(b-)$) the right (resp. left) limit of $v$, which exists since $\mathcal{P}_2(\mathbb{R}^d)$ is complete. The above limit 
exist even in the case $a = -\infty$ (resp. $b = +\infty$) if $v \in AC(a, b; \mathcal{P}_2(\mathbb{R}^d) )$. Among all the possible 
choices of $m$ there exists a minimal one, which is provided by the following theorem (see \cite{7,8,20}): 
\begin{theorem}[Metric derivative]\label{mdert}
For any curve $v$ in $AC^2 (a, b; \mathcal{P}_2(\mathbb{R}^d) )$ the limit \begin{equation}\label{mder}
|v'|(t):=\lim_{s \to t}\frac{\text{d}(v(s),v(t))}{|s-t|} 
\end{equation}
exists for $\mathscr{L}^1$-a.e. $t\in(a,b)$ and it is called metric derivative. Moreover the function $t \mapsto |v'|(t)$ belongs to 
$L^2(a,b)$, it is an admissible integrand for the right hand side of (\ref{ascurve}) and it is minimal in the following sense: 
\begin{equation} \begin{split}
|v'|(t)\leq m(t) \qquad \text{for $\mathscr{L}^1$-a.e. } t\in (a,b), \\
\text{for each function $m$ satisfying (\ref{ascurve}).}
\end{split}
\end{equation}
\end{theorem}

With these two definition the following theorem can be enounced:
\begin{theorem}
Let $I$ be an open interval in $\mathbb{R}$, let $\mu_t:I\to \mathcal{P}_2(\mathbb{R}^d)$ be an absolutely continuous curve and let 
$|\mu'|\in L^1(I)$ be its metric derivative, given by theorem \ref{mdert}. Then there exists a Borel vector field 
$\boldsymbol{v}:(x,t)\mapsto \boldsymbol{v}_t(x)$ such that \begin{equation}
\boldsymbol{v}_t \in L^2(\mu_t;\mathbb{R}^d), \qquad \|\boldsymbol{v}_t\|_{L^2(\mu_t;\mathbb{R}^d)}\leq |\mu'|(t) \qquad for \mathscr{L}^1\text{-a.e. }t\in I,
\end{equation} 
and the continuity equation \ref{ceq} holds in the sense of distributions (\ref{distsense}).  Moreover, for $\mathscr{L}^1$-a.e. $t \in I$, 
$\boldsymbol{v}_t$ belongs to the closure in $L^2(\mu_t,\mathbb{R}^d)$ of the subspace generated by the gradients $\nabla \varphi$ with 
$\varphi \in C_c^{\infty}(\mathbb{R}^d)$.\\ Conversely, if a narrowly continuous curve $\mu_t:I\to \mathcal{P}_2(\mathbb{R}^d)$ satisfies the continuity 
equation for some Borel velocity field $\boldsymbol{v}_t$ with $\|\boldsymbol{v}_t\|_{L^2(\mu_t; \mathbb{R}^d)} \in L^1(I)$ then 
$\mu_t:I \mapsto \mathcal{P}_2(\mathbb{R}^d)$ is absolutely continuous and $|\mu'|(t) \leq \|\boldsymbol{v}_t\|_{L^2(\mu_t; \mathbb{R}^d)}$ 
for $\mathscr{L}^1$-a.e. $t\in I$.
\end{theorem}
The theorem shows in particular that among all velocity fields $\boldsymbol{v}_t$, which produce the same flow $\mu_t$, there is an unique 
optimal one with the smallest $L^2(\mu_t;\mathbb{R}^d)$-norm, equal to the metric derivative of $\mu_t$ up to a negligible set of times 
in $(a,b)$; this optimal field can be viewed 
as the tangent vector field to the curve $\mu_t$. Moreover the minimality of the $L^2$ norm of 
$\boldsymbol{v}_t$ is equivalent to the structure property 
\begin{equation}
\boldsymbol{v}_t \in \overline{\nabla \varphi : \varphi \in C_c^{\infty}(\mathbb{R}^d)}^{L^2(\mu;\mathbb{R}^d)} 
\qquad \text{ for } \mathscr{L}^1 \text{-a.e.} \; t \in (a,b).
\label{ts}
\end{equation}
and notice that: \begin{equation}
|\mu'|(t)= \text{min}\biggl\{ \|\boldsymbol{v}_t\|_{L^2(\mu_t)}: \boldsymbol{v}_t \text{ solves }(\ref{ceq}) \text{ in the sense of distributions}\biggr\}.
\end{equation}

\section{Gradient flows in $\mathcal{P}_2(\mathbb{R}^d)$}

\subsection{A variational definition of gradient flows}

Working in a finite-dimensional smooth setting, the gradient flow of a function $f:\mathbb{M}^d \to \mathbb{R}$ defined on a 
Riemannian manifold $\mathbb{M}^d$ simply means the family of solutions $u:\mathbb{R}\to \mathbb{M}^d$ of the Cauchy problem 
associated to the differential equation \begin{equation}\label{gf1}
\frac{d}{dt}u(t)=-\nabla f\bigl(u(t\bigr) \; \text{ in }T_{u(t)}\mathbb{M}^d,t\in\mathbb{R}; \; u(0)=u_0 \in \mathbb{M}^d
\end{equation}
In fact the equation, living in the tangent space, impose that the velocity vector of the curve ($\boldsymbol{v}_:=\frac{d}{dt}u(t)$) 
is equal to the opposite of the gradient of $f$ at $u(t)$. \\ 

The theory of gradient flows has been extended to more general framework (see \cite{61b,33b,21b,22b}). First moving from the Riemannian 
setting $\mathbb{M}^d$ towards a Hilbert space $H$, where the function $f$ is a proper, convex and  lower semicontinuous functional 
$\phi:H\to(-\infty,+\infty]$. In this case the tangent space is $H$ itself. The result is that $\phi$ admits only a subdifferential 
(possibly defined not everywhere). The equation 
\ref{gf1} become a subdifferential inclusion on the real line: \begin{equation}
u'(t) \in -\partial \phi(u(t)), \; t > 0; \; u(0)=u_0 \in \overline{D(\phi)}.
\end{equation}

A further step of the theory was the extension to more general metric spaces and/or to non-smooth perturbations of a convex functional 
(a series of paper started with \cite{64} and culminated in \cite{64b}). \\

In order to motivate the general definition \ref{defcms}, let us briefly recall the well known case of a dynamical system in a 
Euclidean space. We refer to \cite{AGS} section 1.3-1.4 . 

The gradient $\nabla \phi$ of a smooth real functional $\phi: \mathcal{S}\to \mathbb{R}$ can be defined taking the derivative along regular 
curves, i.e. \begin{equation}
\boldsymbol{g}=\nabla\phi \iff \begin{aligned}
& (\phi \circ v)' = \langle \boldsymbol{g}(v),v' \rangle \\ \text{for every reg}&\text{ular curve }v:(0,+\infty) \to \mathcal{S}
\end{aligned}
\end{equation}
and the modulus $|\nabla \phi|$ is characterized: \begin{equation} \label{themodisc}
g \geq |\nabla\phi| \iff \begin{aligned}
& |\phi \circ v'| \leq g(v)|v'| \\ \text{for every reg}&\text{ular curve }v:(0,+\infty) \to \mathcal{S}
\end{aligned}
\end{equation}
With this notation a steepest descent curve, i.e. a solution of the equation \begin{equation}\label{steepest}
u'(t)=-\nabla \phi(u(t)) \qquad t >0
\end{equation}
satisfy the two scalar condition \begin{gather}
(\phi \circ u)'=-|\nabla \phi||u'|, \\ |u|'=|\nabla \phi|;
\end{gather}
that are equivalent to the single equation\footnote{thanks to Young inequality} \begin{equation} \label{thesingeq}
(\phi \circ u)' = -\frac{1}{2}|u'|^2 -\frac{1}{2}|\nabla \phi|^2 \qquad \text{in }(0,+\infty).
\end{equation}
The formulation is purely metric and can be extended to more general metric spaces when these equations are changed into inequalities. 
Notice that it is sufficient to check just one inequality in (\ref{thesingeq}), since the opposite one is verified along any curve by 
(\ref{themodisc}). The definition of gradient flows needs the following preliminar result, and notice that the framework is the metric 
space $\mathcal{P}_2(\mathbb{R}^d)$.

\begin{definition}[Upper gradients, \cite{92,51}] Given a proper functional $\phi: \mathcal{P}_2(\mathbb{R}^d) \to (-\infty,+\infty]$, 
a function $g: \mathcal{P}_2(\mathbb{R}^d) \to [0,+\infty]$ is a upper gradient for $\phi$ if for every absolutely continuous curve 
$v \in AC^2(a,b;\mathcal{P}_2(\mathbb{R}^d))$ the function $g \circ v$ is Borel and \begin{equation} \label{123}
|\phi(v(t))-\phi(v(s))| \leq \int_s^tg(v(r))|v'|(r)dr \qquad \forall a < s \leq t < b.
\end{equation}
In particular, if $g \circ v' \in L^1(a,b)$ then $\phi \circ v$ is absolutely continuous and \begin{equation}
|(\phi \circ v)'(t)| \leq g(v(t)|v'|(t) \qquad \text{for }\mathscr{L}^1\text{-a.e. }t \in (a,b). 
\end{equation}
\end{definition}

\begin{definition}[Metric gradient flows] \label{defcms} A locally absolutely continuous map $\mu:(a,b) \to \mathcal{P}_2(\mathbb{R}^d)$ is 
a curve of maximal slope for the functional $\phi$ with respect to its upper gradient $g$, if 
$\phi \circ \mu$ is $\mathscr{L}^1$-a.e. equal to a non-increasing map $\varphi$ and \begin{equation}
\varphi ' \leq -\frac{1}{2}|\mu'|^2(t) - \frac{1}{2}g^2(\mu_t) \qquad \text{for $\mathscr{L}^1$-a.e. in $(a,b)$}.  
\end{equation}   
\end{definition}

\subsection{Gradient flows}

We firstly need  to define:
\begin{definition}[Extended Fr\'{e}chet subdifferential] \label{extfresub}
Given a functional $\phi$ satisfying (\ref{1sts},\ref{2nds}) and a measure $\mu^1\in D(\phi)$. A plan 
$\boldsymbol{\gamma}\in \mathcal{P}_{2}(\mathbb{R}^d\times \mathbb{R}^d)$ 
belongs to the Fr\'{e}chet subdifferential $\boldsymbol{\partial}\phi (\mu^1)$ if $\pi_{\#}^1\boldsymbol{\gamma}=\mu^1$ and \begin{equation}
\phi(\mu^3)-\phi(\mu^1) \geq \inf_{\boldsymbol{\mu}\in \Gamma(\boldsymbol{\gamma},\mu^3)}\int_{X^3} \langle x_2,x_3-x_1 \rangle 
d\boldsymbol{\mu} + o\bigl( W_2(\mu^1,\mu^3) \bigl).
\end{equation}
\end{definition}

An important property ensures that the element with minimal norm is concentrated on the graph of a vector field. It can be seen 
 in \cite{AGS} following the approach of sections 10.3-10.4 and in the demonstration of theorem 11.1.3). So 
 $\boldsymbol{\partial}^0\phi(\mu_t)$ is concentrated on the graph of the transport map $-\boldsymbol{v}_t$ for $\mathscr{L}^1$-a.e. 
 $t > 0$, even if the measure $\mu_t$ do not satisfy any regularity assumption. In this case the subdifferential can be defined 
\begin{definition} \label{defsubdiff}
A vector field $\boldsymbol{k} \in L^2(\mu)$ is said to be an element of the subdifferential of $\phi$ at $\mu$ ($\boldsymbol{k} \in 
\partial \phi(\mu)$) if \begin{equation}
\phi(\nu) - \phi(\mu) \geq \inf_{\boldsymbol{\gamma}_0 \in \Gamma_0(\mu,\nu)}\int_{\mathbb{R}^d\times \mathbb{R}^d} 
\boldsymbol{k}(x)\cdot (y-x)d\boldsymbol{\gamma}_0(x,y) +o(W_2(\nu,\mu)).
\end{equation}  
\end{definition} 

After this definition we can finally define a gradient flow:
\begin{definition}[Gradient flows] \label{definizgraflo}
Given a proper and l.s.c functional $\phi:\mathcal{P}_2(\mathbb{R}^d)\to (-\infty,+\infty]$, a map $\mu_t \in 
AC_{loc}^2((0,+\infty);\mathcal{P}_2(\mathbb{R}^d))$ is a solution of the gradient flow equation \begin{equation}
\boldsymbol{v}_t\in-\partial \phi(\mu_t) \qquad t > 0,
\end{equation} 
if denoting by $\boldsymbol{v}_t\in \text{Tan}_{\mu_t}\mathcal{P}_2(\mathbb{R}^d)$ its velocity field, it 
belongs to the subdifferential \ref{defsubdiff} of $\phi$ at $\mu_t$ for $\mathscr{L}^1$-a.e. $t>0$.
\end{definition}

\subsubsection{Equivalence}

Gradient flows and curves of maximal slope are equivalent under certain conditions. Assuming that \begin{equation} \label{1sts}
\phi:\mathcal{P}_2(\mathbb{R}^d)\to (-\infty,+\infty],\text{ proper and lower semi-continuous},
\end{equation}
is such that\footnote{See next section for a detailed description of this assumption}
 \begin{equation} \label{2nds} \begin{split}
\nu \mapsto \Phi(\tau,\mu;\nu)=\frac{1}{2\tau}W_2^2(\mu,\nu)+\phi(\nu) \; \text{admits at least} \\ 
\text{a minimum point }\mu_t, \text{ for all }\tau \in (0,\tau^{\ast}) \text{ and }\mu \in \mathcal{P}_2(\mathbb{R}^d),
\end{split}
\end{equation}
and that the functional satisfy the following regularity property:
\begin{definition}[Regular functionals] \label{regfun}
A given functional satisfying (\ref{1sts}) is regular whenever the strong subdifferential 
$\boldsymbol{\gamma}_n\in \boldsymbol{\partial}\phi(\mu_n),\varphi_n=\phi(\mu_n)$ satisfy \begin{equation}
\begin{split}
\varphi_n \to \varphi \in \mathbb{R},\qquad &  \mu_n \to \mu \qquad \text{in }\mathcal{P}_2(X), \\
\sup_n|\boldsymbol{\gamma}_n|< + \infty,\qquad & \boldsymbol{\gamma}_n \to \boldsymbol{\gamma} \qquad \text{in }\mathcal{P}_2(X \times X).
\end{split}
\end{equation}
then $\boldsymbol{\gamma} \in \boldsymbol{\partial}\phi(\mu)$ and $\varphi= \phi(\mu)$,
\end{definition}

the following general theorem ensure the equivalence: 
\begin{theorem}[Curves of maximal slope coincide with gradient flows] \label{eqcmsgf}
Let $\phi$ be a regular functional according to definition \ref{regfun}, satisfying (\ref{1sts}) and (\ref{2nds}). Then $\mu_t:(0,+\infty) 
\to \mathcal{P}_2(\mathbb{R}^d)$ is a curve of maximal slope w.r.t $|\partial \phi|$ (according to definition \ref{defcms}) iff $\mu_t$ 
is a gradient flow and $t\mapsto \phi(\mu_t)$ is $\mathscr{L}^1$-a.e. equal to a function of bounded variation. In this case the tangent 
vector field $\boldsymbol{v}_t$ to $\mu_t$ satisfies the minimal selection principle \begin{equation}
\boldsymbol{v}_t=-\partial^0\phi(\mu_t) \qquad \text{for }\mathscr{L}^1\text{-a.e. }t>0. 
\end{equation}
\end{theorem}

\subsection{Examples}

The definition \ref{definizgraflo} is equivalent to asking that exists a Borel vector field $\boldsymbol{v}_t \in \text{Tan}_{\mu}\mathcal{P}_2(\mathbb{R}^d)$ 
for $\mathscr{L}^1$-a.e. $t>0,\|\boldsymbol{v}_t\|_{L^2(\mu_t)}\in L_{loc}^2(0,+\infty)$, the continuity equation holds in the sense of 
distribution and \begin{equation}
\boldsymbol{v}_t \in -\partial \phi(\mu_t) \qquad \text{for }\mathscr{L}^1\text{-a.e. }t>0.
\end{equation} 
This equivalence is evident looking at the following evolutionary PDE of diffusion type. \\

Given the parabolic equation in the space-time open cylinder
\begin{equation} \label{stoceq}
\frac{\partial \rho}{\partial t}-\nabla \cdot \left[ \left( \nabla \frac{\delta \mathscr{F}}{\delta \rho} \right) \rho \right] = 0 \; 
\text{in } \mathbb{R}^d \times (0,+\infty),
\end{equation}
where \begin{equation}
\frac{\delta \mathscr{F}(\rho)}{\delta \rho}=-\nabla \cdot F_p(x,\rho,\nabla \rho)+F_z(x,\rho, \nabla \rho).
\end{equation}
is the first variation of a integral functional \begin{equation}
\mathscr{F}(\rho)=\int_{\mathbb{R}^d}F(x,\rho(x),\nabla \rho(x))dx
\end{equation}
associated with a Lagrangian $F=F(x,z,p):\mathbb{R}^d\times[0,+\infty)\times \mathbb{R}^d\to \mathbb{R}$. \\ 
Identifying $\rho$ with the measure $\mu_t:=\rho(\cdot,t)\cdot \mathscr{L}^d$ and considering the functional $\mathscr{F}$ defined on 
$\mathcal{P}_2(\mathbb{R}^d)$; $\mu_t(x)$ is the gradient flow solution for the functional $\mathscr{F}$ if $\rho$ is the solution of the 
system that arise from the equation \ref{stoceq}: \begin{align} 
\label{la234} \partial_t \rho + \nabla \cdot (\rho \boldsymbol{v})=0 & & \text{(continuity equation)}, \\
\rho \boldsymbol{v}=\rho \nabla \psi & & \text{(tangent condition)},\\
\label{la236} \psi =-\frac{\delta \mathscr{F}(\rho)}{\delta \rho} & & \text{(differential inclusion)}.
\end{align}
Here the equivalence is obviously true only for smooth data, but theory of gradient flows presents numerous improvements for the research 
of the solution. 
Many well-known equations for probability densities can be recast in the formalism of gradient flows (see \cite{VILTO}). One has the following 
correspondence between energy functionals on the one hand, and gradient 
flows with respect to the differential structure induced by optimal transportation on the other hand \begin{align*}
& \mathcal{E}[\mu]=\int \mu log(\mu)& & \partial_t \mu = \Delta \mu \\
& \mathcal{E}[\mu]=\int \mu log(\mu) + \int \mu V & & \partial_t \mu = \Delta \mu + \nabla \cdot (\mu \nabla V)\\
& \mathcal{E}[\mu]= \frac{1}{m-1} \int \mu^m & &\partial_t \mu = \Delta \mu^m\\
& \mathcal{E}[\mu]=\frac{1}{2} \int W(x-y) d\mu(x)d\mu(y) & &\partial_t \mu = \nabla \cdot (\mu \nabla (W \ast \mu)) 
\end{align*}

The energy $\int \mu log(\mu)$ is Boltzmann's famous $H$ functional, which has the physical meaning of the negative of an \emph{entropy}. 
The partial differential equations above are known under the respective names of \emph{heat} equation (see \cite{HEAT}), 
\emph{linear Fokker-Planck} equation (see \cite{AS,JKO}), \emph{porous media} equation (see \cite{OT1,OT2}) 
and a particular case of the \emph{McKean-Vlasov} equation. These equations come from the more general Vlasov 
equation in kinetic theory, with no velocity variable in the phase space (see \cite{Vlasov}): \begin{equation*}
\frac{\partial \rho}{\partial t}= \Delta F(\rho) + \nabla \cdot (\rho \nabla V) + \nabla \cdot (\rho \nabla(\rho \ast W)), 
\end{equation*}
related to the functional, sum of a sort of internal, potential and interaction energy: \begin{equation*}
\mathcal{F}[\mu]=\int F(\rho)dx + \int Vd\mu + \frac{1}{2}\int Wd\mu \times \mu \; \text{ if } \mu=\rho \mathscr{L}^d.
\end{equation*}

\subsection{Minimizing movement scheme}

In this section will be analyzed, among the most useful tools of the theory of gradient flows, an approximation procedure for obtaining 
solutions: the Minimizing Movement Scheme. The main idea comes from similarity to the Implicit Euler Method. 
It is a first-order numerical procedure for solving equation like (\ref{steepest}) in an $\mathbb{R}^k$. Given a sequence of time 
steps $\boldsymbol{\tau}:=\{\tau_n\}_{n=1}^{+\infty}$ with $|\boldsymbol{\tau}|=\sup_n(\tau_n)$ and partitioning the interval $(0,+\infty)$ 
\begin{equation}
\begin{split}
P_{\boldsymbol{\tau}}:=\{0= & t_0<t_1<t_2<\cdots t_n < \cdots \}, \qquad I_n:=(t_{n-1},t_n) \\
& t_n:=t_{n-1}+\tau_n, \qquad \lim_{n \to \infty}t_n=\sum_{k=1}^{+\infty}\tau_k=+\infty
\end{split}
\end{equation}
can be found an approximate solution $U_{\boldsymbol{\tau}}^n=u(t_n), \; n=1,\cdots,$ by solving iteratively the following equation, 
with the starting point $U_{\boldsymbol{\tau}}^0 \approx u_0$ \begin{equation}
\frac{U_{\boldsymbol{\tau}}^n-U_{\boldsymbol{\tau}}^{n-1}}{\tau_n}=\nabla \phi(U_{\boldsymbol{\tau}}^n).
\end{equation} 
This is the Euler equation associated to the functional in the variable $V$ for the whom are searched 
minimum points\begin{equation} \label{functional}
\Phi(\tau_n,U_{\boldsymbol{\tau}}^{n-1};V):=\frac{1}{2\tau_n}|V-U_{\boldsymbol{\tau}}^{n-1}|^2+\phi(V) \qquad V \in \mathbb{R}^k.
\end{equation}
The following recursive scheme is so defined:
\begin{equation}	\label{scheme}
\begin{cases}
U_{\boldsymbol{\tau}}^0\text{ is given; whenever }U_{\boldsymbol{\tau}}^1,\cdots,U_{\boldsymbol{\tau}}^{n-1}\text{ are known} \\
\text{find }U_{\boldsymbol{\tau}}^n\in \mathbb{R}^k: \Phi(\tau_n,U_{\boldsymbol{\tau}}^{n-1};U_{\boldsymbol{\tau}}^n) \leq 
\Phi(\tau_n,U_{\boldsymbol{\tau}}^{n-1};V) \; \forall V \in \mathbb{R}^k
\end{cases}
\end{equation}
Now can be defined the \begin{definition}[Discrete solution] If, for a choice of $\boldsymbol{\tau}$ and $U_{\boldsymbol{\tau}}^0\in 
\mathbb{R}^k$ a sequence $\{ U_{\boldsymbol{\tau}}^n \}_{n=1}^{+\infty}$ solving the recursive scheme exists, so the discrete values can 
be interpolated by the piecewise constant function $\overline{U}_{\boldsymbol{\tau}}$, defined by \begin{equation} \label{discsol}
\overline{U}_{\boldsymbol{\tau}}(0)=U_{\boldsymbol{\tau}}^0, \qquad \overline{U}_{\boldsymbol{\tau}}(t) \equiv U_{\boldsymbol{\tau}}^n \;
 \text{if }t\in(t^{n-1},t^{n}], \; n \leq 1.
\end{equation}
$\overline{U}_{\boldsymbol{\tau}}$ is called a discrete solution corresponding to the partition $P_{\boldsymbol{\tau}}$.
\end{definition}

\begin{figure} \centering 
\includegraphics[scale=1]{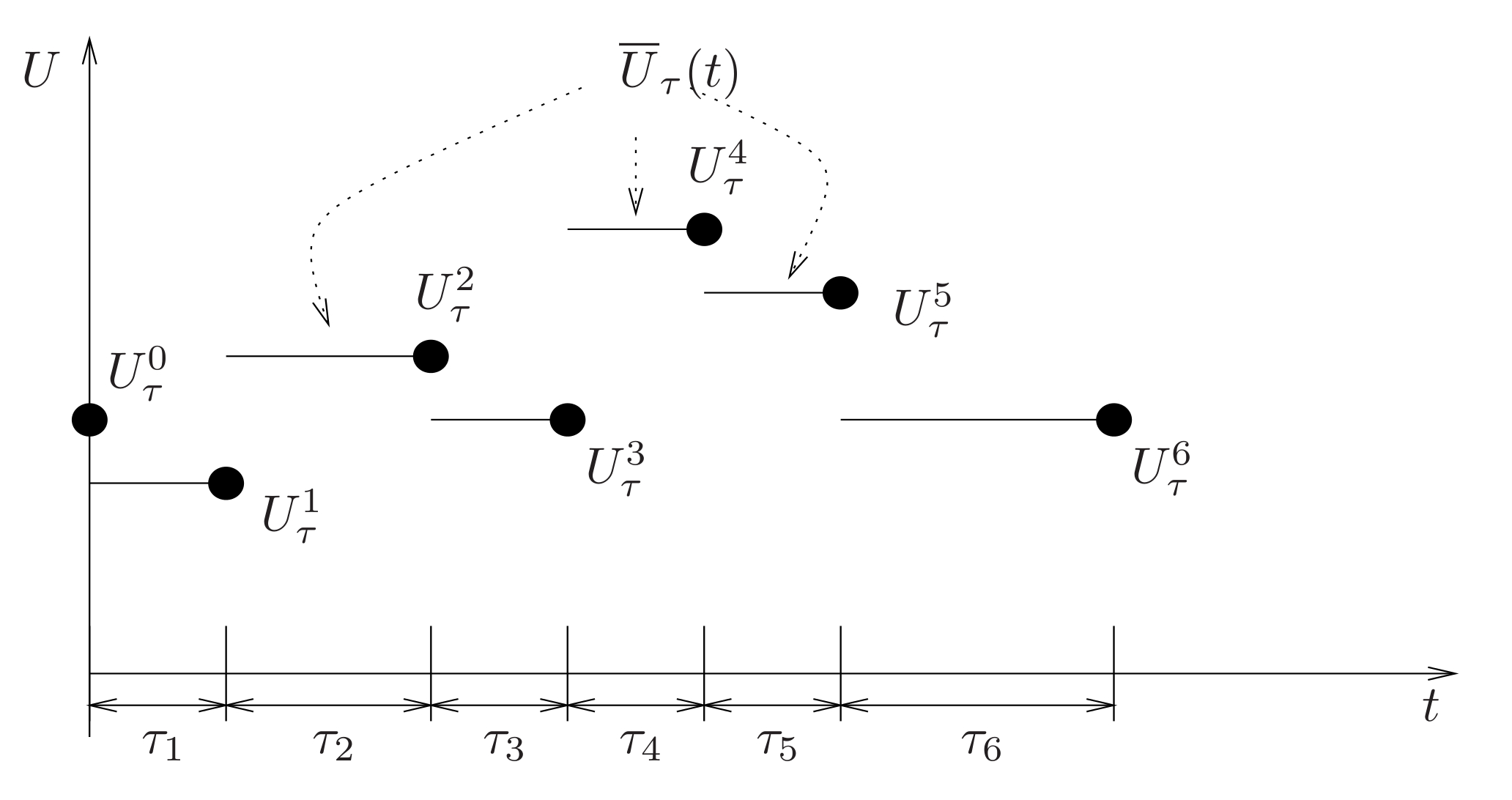}
\caption{Pointwise constant interpolation (image taken from \cite{AGS})}
\end{figure}

The scheme can be enlarged to a more general metric context. We replace $\mathbb{R}^k$ with $\mathcal{P}_2(\mathbb{R}^d)$ and the modulus 
in the equation \ref{functional} with the Wasserstein distance between measures.

After this remark can be defined: 
\begin{definition}[Minimizing movements]
Given the functional $\Phi$ defined as in (\ref{functional}) and an initial datum $u_0 \in \mathcal{P}_2(\mathbb{R}^d)$, a curve 
$u:[0,+\infty) \to \mathcal{P}_2(\mathbb{R}^d)$ is a minimizing movement for $\Phi$ starting from $u_0$ if for every partition 
$\boldsymbol{\tau}$ (with sufficiently small $|\boldsymbol{\tau}|$) there exists a discrete solution $\overline{U}_{\boldsymbol{\tau}}$ 
defined as in (\ref{scheme}),(\ref{discsol}) such that \begin{equation}
\begin{split}
\lim_{|\boldsymbol{\tau}| \downarrow 0}\phi (U_{\boldsymbol{\tau}}^0)&=\phi(u_0), \qquad \limsup_{|\boldsymbol{\tau}|\downarrow 0}
\text{d}((U_{\boldsymbol{\tau}}^0),u_0)< +\infty , \\
&\overline{U}_{\boldsymbol{\tau}}(t)\to u(t) \text{ narrowly } \; \forall t \in [0,+\infty).
\end{split}
\end{equation}
The collection of all the minimizing movements for $\Phi$ starting from $u_0$ is denoted by $MM(\Phi;u_0)$. \\ Analogously, a curve 
$u:[0,+\infty)\to \mathcal{P}_2(\mathbb{R}^d)$ is a generalized minimizing movement for $\Phi$ starting from $u_0$ if there exists a 
sequence of partitions $\boldsymbol{\tau}_k$ with $|\boldsymbol{\tau}_k| \downarrow 0$ and a corresponding sequence of discrete solutions 
$\overline{U}_{\boldsymbol{\tau}_k}$ defined as in (\ref{scheme}),(\ref{discsol}) such that \begin{equation}
\begin{split}
\lim_{k \to \infty }\phi (U_{\boldsymbol{\tau}_k}^0)&=\phi(u_0), \qquad \limsup_{k \to \infty} \text{d}((U_{\boldsymbol{\tau}_k}^0),u_0)< +\infty , \\
& \overline{U}_{\boldsymbol{\tau}_k}(t)\to u(t) \text{ narrowly }\; \forall t \in [0,+\infty).
\end{split}
\end{equation}
\end{definition}


\subsection{A general existence, uniqueness and convergence result}
\label{subsection}

Three main results are recalled here. Two of them are technical: one is about the convergence of the MMS and the other about the existence 
of the subdifferential. The third is a theorem about existence and uniqueness of gradients flows.

\subsubsection{Convergence of MMS}

The existence of discrete solutions and their convergence to a absolutely continuous curve can be demonstrated  (see corollary 2.2.2. and 
proposition 2.2.3 in \cite{AGS}) with the following assumptions on the functional:
\begin{itemize}
\item \textbf{Lower semicontinuity} $\phi$ is sequentially $\sigma$-lower semicontinuous on $d$-bounded sets \begin{equation}
\sup_{n,m}d(u_n,u_m) < +\infty, \qquad u_n \rightharpoonup u \Rightarrow \liminf_{n \to \infty}\phi(u_n) \geq \phi(u).
\end{equation}
\item \textbf{Coercivity} There exists $\tau_{\ast} > 0$ and $u_{\ast} \in \mathcal{P}_2(\mathbb{R}^d)$ such that \begin{equation}
\phi_{\tau_{\ast}}(u_{\ast}):= \inf_{v \in \mathcal{P}_2(\mathbb{R}^d)}\Phi(\tau_{\ast},u_{\ast};v) > -\infty.
\end{equation}
\item \textbf{Compactness} Every bounded set contained in a sublevel of $\phi$ is relatively $\sigma$-sequentially compact: i.e. 
\begin{equation}
\begin{split}
\text{if }(u_n)& \subset \mathcal{P}_2(\mathbb{R}^d) \text{ with }\sup_n\phi(u_n) < +\infty, \qquad \sup_{n,m}d(u_n,u_m) < +\infty \\
& \text{then }(u_n) \text{ admits a }\sigma\text{-convergent subsequence.}
\end{split}
\end{equation}
\end{itemize} 

\subsubsection{Existence of the subdifferential}

Following the approach of \cite{AGS} chapter 10, is considered a functional \begin{equation} \label{cond1}
\phi:\mathcal{P}_2(\mathbb{R}^d)\to (-\infty,+\infty] \text{ proper and lower semicontinuous}
\end{equation}
such that \begin{equation} \label{cond2} \begin{split}
\nu \mapsto \Phi(\tau,\mu;\nu)=\frac{1}{2\tau}W_2^2(\mu,\nu) + \phi(\nu) \text{ admits at least} \\
\text{a minimum point }\mu_{\tau}, \text{ for all }\tau \in (0,\tau_{\ast}) \text{ and }\mu \in \mathcal{P}_2(\mathbb{R}^d).
\end{split}
\end{equation}

These assumptions lead to a result of existence for the subdifferential in the minimizing point:
\begin{lemma} \label{illemmone}
Let $\phi:\mathcal{P}_2(\mathbb{R}^d)\to (-\infty,+\infty]$ be satisfying (\ref{cond1}) and let $\mu_{\tau}$ be a minimizer of 
(\ref{cond2}); if $\hat{\boldsymbol{\gamma}}_{\tau} \in \Gamma_0(\mu_{\tau},\mu)$, then the rescaled plans \begin{equation}
\boldsymbol{\gamma}_{\tau}:=(\boldsymbol{\rho}_{\tau})_{\#}\hat{\boldsymbol{\gamma}}_{\tau} \text{ with } \boldsymbol{\rho}_{\tau}(x_1,x_2):=
\biggl( x_1,\bigl(\frac{x_2-x_1}{\tau}\bigr)\biggr)
\end{equation}
and the associated plans $\boldsymbol{\mu}_{\tau} \in \Gamma_0(\boldsymbol{\gamma}_{\tau},\mu^3)$ satisfy \begin{equation} \label{precedente}
\phi(\mu^3)-\phi(\mu_{\tau})-\int_{\mathbb{R}^{3d}}\langle x_2,x_3-x_1 \rangle d\boldsymbol{\mu}_{\tau} \geq 
o\bigl(W_2(\mu^3,\mu_{\tau})\bigr).
\end{equation}
and $\boldsymbol{\gamma}_{\tau}$ is a subdifferential, according to \ref{extfresub}.
\end{lemma}

This result is very general but the subdifferential can be characterized explicitly under some assumptions.
We consider the functional $\phi:\mathcal{P}_2(\mathbb{R}^d)\to (-\infty,+\infty]$
given by the sum of internal, potential and interaction energy:
\begin{equation}
\phi(\mu)=\int_{\mathbb{R}^d} F(\rho)dx + \int_{\mathbb{R}^d} Vd\mu + \frac{1}{2}\int_{\mathbb{R}^d \times \mathbb{R}^d}
 Wd\mu \times \mu \; \text{ if } \mu=\rho \mathscr{L}^d.
\end{equation}
setting $\phi(\mu)=+\infty$ if $\mu \in \mathcal{P}_2(\mathbb{R}^d) \setminus \mathcal{P}_2^a(\mathbb{R}^d)$. 
Recalling the "doubling condition" ($ \exists C > 0 : K(x+y) \leq C(1+K(x)+K(y))$) we make the following assumptions on $F$, $V$ and $W$:
\begin{itemize}
\item $F : [0,+\infty) \to \mathbb{R}$ is a doubling, convex differentiable function with superlinear growth satisfying 
\begin{gather*} 
\text{the map }s \mapsto s^d F (s^{-d}) \text{ is convex and non increasing in }(0,+\infty), \\
F(0)=0, \liminf_{s \downarrow 0}\frac{F(s)}{s^{\alpha}}>-\infty \text{ for some }\alpha > \frac{d}{d+2};
\end{gather*}\item $V: \mathbb{R}^d \to (-\infty,+\infty]$ is a l.s.c. $\lambda$-convex function with proper domain $D(V)$
with nonempty interior $\Omega \subset \mathbb{R}^d$;
\item $W: \mathbb{R}^d \to [0,+\infty)$ is a convex, differentiable, even function satisfying the doubling condition.
\end{itemize}
We have the following characterization of the minimal selection in the subdifferential $\partial^0 \phi(\mu)$:
\begin{theorem}[Minimal subdifferential of $\phi$] A measure $\mu=\rho \mathscr{L}^d \in D(\phi) \subset \mathcal{P}_2(\mathbb{R}^d)$ 
belongs to $D(|\partial \phi|)$ if and only if  $\bigl( zF'(z)-F(z) \bigr)=L_F(z) \in W^{1,1}_{loc}(\Omega)$ and \begin{equation} \label{ecchedefinit}
\rho \boldsymbol{w}= \nabla L_F(\rho)+\rho \nabla V + \rho (\nabla W \ast \rho) \text{ for some }\boldsymbol{w}\in L^2(\mu,\mathbb{R}^d).
\end{equation}
In this case the vector $\boldsymbol{w}$ defined $\mu$-a.e. by (\ref{ecchedefinit}) is the minimal selection in $\partial \phi(\mu)$, 
i.e. $\boldsymbol{w}=\partial^0\phi(\mu)$.
\end{theorem}

\subsubsection{Existence and uniqueness of gradient flows}

In this section we are considering the case of a \begin{equation} \label{eccheassu}
\begin{split}
\text{proper, l.s.c. and } & \text{coercive functional} \phi: \mathcal{P}_2(\mathbb{R}^d)\to (-\infty,+\infty], \\
\text{which is }\lambda & \text{-convex along generalised geodesics,}
\end{split}
\end{equation} 
according to the following definition:
\begin{definition}[Generalized geodesics] \label{gengeo}
A generalized geodesics joining $\mu^2$ to $\mu^3$ (with base $\mu^1$)is a curve of the type 
\begin{equation*}
\mu_t^{2 \to 3}=(\pi_t^{2 \to 3})_{\#}\boldsymbol{\mu} \qquad t \in [0,1],
\end{equation*}
where \begin{equation} \label{congengeo}
\boldsymbol{\mu}\in \Gamma(\mu^1,\mu^2,\mu^3) \text{ and } \pi_{\#}^{1,2}\boldsymbol{\mu}\in \Gamma_0(\mu^1,\mu^2), \; 
\pi_{\#}^{1,3}\boldsymbol{\mu}\in \Gamma_0(\mu^1,\mu^3).
\end{equation}
\end{definition}

\begin{definition}[Convexity along generalized geodesics] \label{laconvexxgener} Given $\lambda \in \mathbb{R}$, $\phi$ is $\lambda$-convex along generalized 
geodesics if for any $\mu^1,\mu^2,\mu^3 \in D(\phi)$ there exists a generalized geodesic 
$\mu_t^{2 \to 3}$ induced by a plan $\boldsymbol{\mu}\in\Gamma(\mu^1,\mu^2,\mu^3)$ satisfying (\ref{congengeo}) such that \begin{equation}
\phi(\mu_t^{2 \to 3}) \leq (1-t)\phi(\mu^2)+y\phi(\mu^3)-\frac{\lambda}{2}t(1-t)\int_{\mathbb{R}^3d}|x^2-x^3|d\boldsymbol{\mu}.
\end{equation}
\end{definition}

Under this assumptions can be enounced the theorem: 
\begin{theorem}[Existence and uniqueness of gradient flows] \label{exiunigraflo}
 Let us suppose that $\phi:\mathcal{P}_2(\mathbb{R}^d)\to (-\infty,+\infty]$ satisfies (\ref{eccheassu}) and let $\mu_0 \in D(\phi)$. 
The discrete solution $\overline{U}_{\boldsymbol{\tau}}$ given by (\ref{discsol}) converges locally uniformly
to a locally Lipschitz curve $\mu \in \mathcal{P}_2(\mathbb{R}^d) $ which is the unique gradient
flow of $\phi$ with $\mu(0+)=\mu_0$.
\end{theorem}

\clearpage

\lhead[\thepage]{\small{\textsc{Granular media and GF of smooth interaction potential}}}
\rhead[\small{\textsc{Granular media and GF of smooth interaction potential} }] {\thepage}	

\clearpage{\pagestyle{empty}\cleardoublepage}

\chapter{Granular media and GF of smooth interaction potential}
\markboth{Granular media and GF of smooth interaction potential}{Granular media and GF of smooth interaction potential}

\section{A kinetic equation for granular media}

Over the last years, due to industrial application and to the evolution of the trends in theoretical physics, a lot of attention was 
given to the modelling of granular material (sand, powders, heaps of cereals, grains, molecules, snow, or even asteroids...).
The main features of these systems are the possibility of the occurrence of inelastic collapses (namely infinitely many collisions in 
a finite time) and the tendency of  the system to clusterize, that is to create states of concentration of the density, as
sand grains over a shaken sheet of paper. The model is described here because it presents a gradient flow representation, as will be 
shown in the end of the section. \\

Here we will focus on the one-dimensional particle systems, following the derivation of the model made 
in \cite{i4}. The literature of this model is huge, see for example \cite{VILNEW,MCNAMA,46,6,11,LN,LNS,RAOUL,FELRAO,FELRAO2,LT,28,i18}.\\

The equation is derived considering a one-dimensional system constituted by N particles on the line,
colliding inelastically. Then we rescale suitably the degree of inelasticity, as well
as the total number of particles (which is assumed to diverge), to obtain a kinetic
equation for the one-particle probability density. \\

The dynamics of the system is defined in the following way: denote by $x_1,\cdots,x_N$ and $v_1,\cdots,v_N$ positions and velocities 
of the particles. The
particles goes freely up to the first instant in which two of them are in the same
point. Then they collide according to the rule: \begin{equation} \label{equavvpri}
v'=v_1+\varepsilon(v-v_1), \qquad v_1'=v-\varepsilon(v-v_1), 
\end{equation}
where $v',v_1'$ and $v,v_1$ are the outgoing and ingoing velocities respectively and $\varepsilon$ is a real parameter measuring the degree 
of inelasticity of the collision. Notice that the total momentum is conserved in the collision, while the modulus of the relative
velocity decreases by a fixed rate for any collision. Then the particles go on up to
the instant of the next collision which is performed by the same rule and so on.
Since the particles are assumed to be identical, the physics does not change if we replace the law \ref{equavvpri} by the following one 
\begin{equation} \label{equavvprimo}
v'=v-\varepsilon(v-v_1), \qquad v_1'=v_1+\varepsilon(v-v_1),
\end{equation}
which is the same as equation \ref{equavvpri} with the names of the particles exchanged after the collision, but it is often easier to do 
computations using (\ref{equavvprimo}).

The ordinary differential equation governing the time evolution of the system is: \begin{equation}
\dot{x}_i=v_i, \qquad \dot{v}_i=\varepsilon\sum_{j=1}^N\delta(x_i-x_j)(v_j-v_i)|v_j-v_i|. 
\end{equation}

Notice that $\varepsilon (v_j-v_i)$ is the jump performed by the particle $i$ after a collision with the particle $j$, while 
$\delta(x_i-x_j)|v_j-v_i|=\delta(t-t_{i,j})$, being $t_{i,j}$ the instant of the impact between the particle $i$ with the particle $j$.

Let $\mu^N=\mu^N(x_1,v_1,\cdots,x_N,v_N)$ be a probability density for the system. The Liouville equation describing its time evolution 
read as: \begin{equation}
\begin{split}
&(\partial_t+\sum_{i=1}^Nv_i\partial_{x_i})\mu^N(x_1,v_1,\cdots,x_N,v_N)= \\
-\varepsilon \sum_{i \neq j}&\delta(x_i-x_j)\partial_{v_i}[\phi(v_j-v_i) \mu^N(x_1,v_1,\cdots,x_N,v_N)]
\end{split}
\end{equation}
where $\phi(\overline{v}-v)=(\overline{v}-v)|\overline{v}-v|$.

Proceeding as in the derivation of the BBKGY hierarchy for Hamiltonian systems, we introduce the $j$-particle distribution functions: 
\begin{equation} \label{equazpardisfun}
f_j^N(x_1,v_1,\cdots,x_j,v_j)=\int dx_{j+1}dv_{j+1}\cdots dx_Ndv_N\mu^N(x_1,v_1,\cdots,x_N,v_N)
\end{equation}
and integrating (\ref{equazpardisfun}) over the last variables, we obtain the following hierarchy of equations
\begin{gather} \label{launopuntosei}
\nonumber (\partial_t+\sum_{i=1}^jv_i\partial_{x_i})f_j^N(x_1,v_1,\cdots,x_j,v_j) = \\
\nonumber -\varepsilon \sum_{i \neq k}^j\delta(x_i-x_k)\phi(v_k-v_i) \partial_{v_i}f_j^N(x_1,v_1,\cdots,x_j,v_j)] + \\
-\varepsilon(N-j) \sum_{i =1}^j\partial_{v_i}\int dv_{j+1} \phi(v_{j+1}-v_i) f_{j+1}^N(x_1,v_1,\cdots,x_{j+1},v_{j+1}).
\end{gather}

An inspection of (\ref{launopuntosei}) suggest the scaling limit, namely $\varepsilon \to 0, N \to \infty$ in such a way that 
$N\varepsilon\to \kappa$, where $\kappa$ is a positive parameter. If $f_j^N$ have a limit (say $f_j$) they are expected to 
satisfy the following (infinite) hierarchy of equations: \begin{gather} \label{launopuntosette}
\nonumber (\partial_t+\sum_{i=1}^jv_i\partial_{x_i})f_j(x_1,v_1,\cdots,x_j,v_j) = \\
-\kappa \sum_{i =1}^j\partial_{v_i}\int dv_{j+1} \phi(v_{j+1}-v_i) f_{j+1}(x_1,v_1,\cdots,x_{j+1},v_{j+1}).
\end{gather}

Finally, it the initial state is chaotic, namely if initially: \begin{equation*}
f_j(x_1,v_1,\cdots,x_j,v_j)=\prod_{i=1}^jf_0(x_i,v_i), 
\end{equation*}
then we expect that the dynamics does not create correlations (propagation of chaos) so that:
\begin{equation*}
f_j(x_1,v_1,\cdots,x_j,v_j,t)=\prod_{i=1}^jf_0(x_i,v_i,t), 
\end{equation*}
by which we obtain, for the one particle distribution function, the kinetic equation: \begin{equation} \label{kinectdist}
(\partial_t+ v \partial_x )f(x,v)=-\kappa\partial_v (Ff),
\end{equation}
where: \begin{equation}
F(x,v,t)=\int d\overline{v}\phi(\overline{v}-v)f(x,\overline{v},t).
\end{equation}
In facts products of solution of (\ref{kinectdist}) are solutions of the hierarchy (\ref{launopuntosette}) as follows by a simple 
algebraic computation. The analysis of (\ref{kinectdist}) is considerably simplified whenever the medium is considered to be 
spatially homogeneous. In this case, setting $\kappa=1$, we have: \begin{equation} \label{laequazionedelref}
\partial_tf(v,t)+\partial_v(Ff)(v,t)=0.
\end{equation} 
We stress that the analysis of a spatially homogeneous medium is not academic. Indeed it is interesting to investigate carefully 
what happens locally, when the homogeneous regime is dominant.
The equation \ref{laequazionedelref} can be modified with the introduction of other terms. They often represent collision or other 
phenomenon of interests. Here we will focus on terms that can be interpreted as a gradient flow.

Many experiments about granular material include shaking, as a way to input energy into the system, counterbalancing the freezing due 
to energy loss. A rather trivial but seemingly not so absurd model consists in a heat bath, or white noise forcing: this amounts to 
adding to the right-hand side of the kinetic equation the term \begin{equation*}
T_e\partial_v^2f, \qquad T_e >0,
\end{equation*}
where $T_e$ is an "external temperature".

Particles may experience extra friction forces if they are going through a viscous fluid or something to that effect. Here again, 
there is a trivial model consisting in adding a drift term to the equation; the most simple case being that of a linear drift,
\begin{equation*}
\alpha \partial_v(vf), \qquad \alpha > 0.
\end{equation*}

The complete equation now is
\begin{equation}
 \partial_t f(x,v,t)=-\kappa \partial_v(Ff)+ \alpha \partial_v(vf) + T_e \partial_v^2f
\end{equation}
with $F=\partial_v W\ast f$.

It can be written in the sense of gradient flow, i.e.:\begin{equation} \label{stoceqeeli}
\frac{\partial f}{\partial t}=\nabla \cdot \bigl( f(f\ast \nabla W) \bigr) 
+ \beta \nabla \cdot (vf) + \sigma \Delta f=\nabla \cdot \biggl( f \nabla\frac{\delta \mathcal{E}}{\delta f} \biggr),
\end{equation}
with the functional $E[f]$ defined in the following way: \begin{equation}
\mathcal{E}[f]=\frac{1}{2}\int f(v)f(w)W(v-w)dvdw+\frac{\beta}{2}\int f|v|^2dv+\sigma \int f log(f).
\end{equation}

One advantage of identifying a gradient flow structure is that it yields interesting recipes for computing, say derivatives of 
functionals along the flow, in terms of gradients and Hessians; this is what is described in \cite{AGS,VILTO} as Otto's calculus. 
Another advantage of such a formalism is that the convexity properties of the energy functional might help the study of convergence
to equilibrium. For instance, if $\mathcal{E}$ is $\lambda$-uniformly convex, then trajectories get closer to each other like 
$O(e^{-\lambda t})$, and there is exponential convergence to the unique infimum of $\mathcal{E}$. This also comes with automatic 
useful energy inequalities.

\section{The case of the interaction functional}

First it must be noticed that with respect to the previous section, here the variable $v$ is substituted with $x$. 
Here we study an equation like (\ref{stoceqeeli}) with only the interaction functional. Following the work of \cite{IL} the standard 
assumptions of the gradient flow theory
can be slightly weakened. Here is presented the framework, the steps of the proofs for the relaxed assumptions and an aggregation result for a 
particular class of potential.  

\subsubsection{The framework}
Consider a mass distribution of particles $\mu\geq 0$, interacting under a continuous interaction potential $W$. 
The associated interaction energy is defined as \begin{equation}
\mathcal{W}[\mu]:=\frac{1}{2}\int_{\mathbb{R}^d\times \mathbb{R}^d}W(x-y)d\mu(x)d\mu(y).
\end{equation}
The equation takes the form \begin{equation}
\label{theeq}
\frac{\partial \mu}{\partial t}=\nabla \cdot \left[ \left( \nabla \frac{\delta \mathcal{W}}{\delta \mu} \right) \mu \right]=\nabla \cdot 
\left[ \left( \nabla W \ast \mu \right) \mu \right] 
\end{equation}
\begin{equation*}
\mu(0)=\mu_0 \qquad x\in \mathbb{R}^d,t\geq 0
\end{equation*}

The velocity field $\boldsymbol{v}_t$ in the equations is defined to be $-(\nabla W \ast \mu)(t,x)$ and it represents the interaction, at 
the point $x$, between particles through the interaction potential $W$ (these are pure non local interactions). Without loss of generality 
masses are normalized to 1 because of an invariance of the equation: if $\mu(t)$ is a solution, so is $M\mu(Mt)$ for every $M > 0$.\\

The choice of $W$ strongly depends on the phenomenon studied. In the framework of population
dynamics is described the evolution of a density of individuals. In a first approximation, it's reasonable that the interaction between 
two individuals depends only on the distance between them. Because of this $W$ can be choosen to be a radial function, 
i.e. $W (x) = w(|x|)$. Moreover, when $w (r) > 0$ the force is attractive among the particles
(or individuals), when $w (r) < 0$ a repulsive force is acting.
The equation \ref{theeq} can be applied in physics and biology. Applications with a regular potential are e.g.:
a simplified inelastic interaction models for granular media studied in \cite{i4,i18} with $W = \frac{|x|^3}{3}$ and in 
\cite{LT,28} with $W = |x|^\alpha$ $\alpha > 1$, flocking and swarming that are described with the quadratic attractive-repulsive Morse potential 
$W (x) = -C_a e^{-|x|^2 / l_a} + C_r e^{\-|x|^2 / l_r}$ (see for example \cite{41,42}). When the potential is singular at 
$x=0$ and attractive, a relevant application is chemotaxis in $2D$ with $W(x) = \frac{1}{2\pi} log|x|$ (see \cite{BDP,BCC}). Other 
application can be found for a repulsive and singular at $x=0$ potential: swarming described with a Morse potential 
$W (x) = -C_a e^{-|x|/l_a} + C_r e^{-|x|/l_r}$ (see \cite{21,41}) or phenomenon in Physics related to a Lennard-Jones type potentials 
(see \cite{The}). \\

The potential must satisfy the following assumptions:
\begin{description}
\item[A] $W$ is continuous, $W(x)=W(-x)$, $W(0)=0$ and $W \in C^1(\mathbb{R}^d\setminus \{0\})$. 
\item[B] There exists a constant $C > 0$ such that \begin{equation*}
W(z) \leq C(1+|z|^2), \qquad \forall z \in \mathbb{R}^d
\end{equation*}
\item[C] $W$ is $\lambda$-convex for some $\lambda \leq 0$, i.e. $W(x)-\frac{\lambda}{2}|x|^2$ is convex. 
\end{description}

\begin{remark} \label{qweretyty}
With respect to the hypothesis of \cite{AGS,AS} here is not assumed that $W$ is differentiable at the origin and it can presents a  
negative quadratic behaviour at infinity. Because of the weaker assumptions is not trivial to find the explicit form of the subdifferential. 
The potential $\mathcal{W}$ is well-defined in $\mathcal{P}_2(\mathbb{R}^d)$ and the continuity of $W$ eliminates possible problems with 
singular measures. The main idea of interest is the strategy used to show the convergence of the scheme. In fact the functional is 
$\lambda$-convex with respect to generalized geodesics, see definition \ref{laconvexxgener}. It follows directly from \cite{AGS} in  
proposition 9.3.5, with the decomposition $\mathcal{W} = \overline{\mathcal{W}} + \mathcal{Q}$ (with $\overline{\mathcal{W}}$ generated 
by a convex potential and $\mathcal{Q}$ generated by $\lambda|x|^2$). Exploiting this fact, the existence of solutions for the discrete 
scheme and the convergence of the scheme easily follows. The method provided by the article is of interest in itself because can be used   
to study a more general situation where $\lambda$-convexity fails.
\end{remark}

Thanks to Lemma \ref{contire} the solution can be defined in this way: \begin{definition} \label{definiz1punto2}
A locally absolutely continuous curve $\mu:[0,+\infty) \ni t \mapsto \mathcal{P}_2(\mathbb{R}^d)$ is said to be a weak measure 
solution to (\ref{theeq}) with initial datum $\mu_0 \in \mathcal{P}_2(\mathbb{R}^d)$ if $\partial^0W \ast \mu$ belongs to 
$L^1_{loc}([0,+\infty);L^2(\mu_t))$ and \begin{equation}
\begin{split}
\int_0^{+ \infty} \int_{\mathbb{R}^d} \frac{\partial \varphi}{\partial t} d\mu_t(x)dt + \int_{\mathbb{R}^d}\varphi(x,0)d\mu_0(x) = \\
\int_0^{+ \infty} \int_{\mathbb{R}^d \times \mathbb{R}^d}\nabla \varphi (x,t) \cdot \partial^0 W(x-y) d\mu_t(x)d\mu_t(y)dt,
\end{split}
\end{equation}
for all test function $\varphi \in C_c^{\infty}([0,+\infty)\times \mathbb{R}^d)$.
\end{definition}
The idea is to search solution in the sense of \emph{curve of maximal slope} and in the end is demonstrated the equivalence to a 
\emph{gradient flow} solution. \begin{definition} \label{lalorodef} A locally absolutely 
continuous curve $[0,T] \ni t \mapsto \mu_t \in \mathcal{P}_2(\mathbb{R}^d)$ is a curve of maximal slope for the functional 
$\mathcal{W}$ if $t \mapsto \mathcal{W}[\mu_t]$ is an absolutely continuous function, and the following inequality holds for every 
$0 \leq s \leq t \leq T$: \begin{equation} \label{energidisug}
\frac{1}{2}\int_s^t|\mu'|^2(r)dr+\frac{1}{2}\int_s^t|\partial \mathcal{W}|^2[\mu_r]dr \leq \mathcal{W}[\mu_s] - \mathcal{W}[\mu_t].
\end{equation}
\end{definition}

\subsubsection{The proof}

The first step of the demonstration is to characterize the minimal element of the subdifferential of $\mathcal{W}$: 
$\partial^0\mathcal{W}[\mu]$, see definition \ref{defsubdiff}.
The characterization of $\partial^0 \mathcal{W}[\mu]$ is achieved through \begin{proposition} \label{classssubdif}
Given a potential satisfying \textbf{A - C}, the vector field \begin{equation}
k(x):=(\partial^0 W \ast \mu)(x)=\int_{y \neq x}\nabla W(x-y)d\mu(y)
\end{equation}
is the unique element of minimal $L^2(\mu)$-norm in the subdifferential of $\mathcal{W}$, i.e. $\partial^0 W \ast \mu = \partial^0 
\mathcal{W}[\mu]$\footnote{Obviously $\partial^0 W$ is the element of minimal norm of the subdiff. of the \emph{function} $W(x)$}.
\end{proposition}
The next step is based on the minimizing movement scheme. Given an initial measure $\mu_0 \in \mathcal{P}_2(\mathbb{R}^d)$ and a time step 
$\tau > 0 $, is considered a sequence $\mu_k^{\tau}$ recursively defined by $\mu_0^{\tau}=\mu_0$ and \begin{equation}
\mu_{k+1}^{\tau}\in \text{argmin}_{\mu \in \mathcal{P}_2(\mathbb{R}^d)}  \biggl\{ \mathcal{W}[\mu]+\frac{1}{2\tau}W_2^2(\mu_k^{\tau},\mu) \biggr\}.
\end{equation}

Because of the quadratic behaviour and of the pointy singularity the existence 
of the minimizer in the scheme is not trivial. A preliminary result is needed: \begin{lemma}[Weak lower semi-continuity of the penalized 
interaction energy] Suppose $W$ satisfies \textbf{A - C}. Then, for a fixed $\overline{\mu}\in \mathcal{P}_2(\mathbb{R}^d)$, the penalized 
interaction energy functional \begin{equation}
\mathcal{P}_2(\mathbb{R}^d) \ni \mu \mapsto \mathcal{W}[\mu]+\frac{1}{2\tau}W_2^2(\mu,\overline{\mu})
\end{equation}
is lower semi-continuous with respect to the narrow topology of $\mathcal{P}(\mathbb{R}^d)$ for all $\tau > 0 $ such that $8 \tau \lambda^- 
\leq 1$, where $\lambda^-:=\text{max} \left\{ 0,-\lambda \right\}$. 
\end{lemma}
After this lemma can be demonstrated that: \begin{proposition}[Existence of minimizers] Suppose $W$ satisfies \textbf{A - C}. Then, there 
exists $\tau_0 > 0$ depending only on $W$ such that, for all $0 < \tau < \tau_0$ and for a given $\overline{\mu} \in \mathcal{P}_2(\mathbb{R}^d)$, 
there is $\mu_{\infty} \in \mathcal{P}_2(\mathbb{R}^d)$ such that \begin{equation}
\mathcal{W}[\mu_{\infty}]+\frac{1}{2\tau}W_2^2(\mu_{\infty},\overline{\mu}) = \min_{\mu \in \mathcal{P}_2(\mathbb{R}^d)} \biggl\{ 
\mathcal{W}[\mu]+\frac{1}{2\tau}W_2^2(\mu,\overline{\mu}) \biggr\}.
\end{equation} 
\end{proposition}

Introducing the piecewise constant interpolation \begin{equation*}
\begin{cases}
\mu^{\tau}(0):=\mu_0   \\
\mu^{\tau}(t):=\mu_k^{\tau} \text{ if } t \in \bigl((k-1)\tau,k\tau\bigr], \; k \leq 1
\end{cases} 
\end{equation*}
the convergence can be proved. The proof can be found in \cite{AGS} proposition 2.2.3. The hypothesis are compactness, coercivity and 
lower semi-continuity that have been demonstrated in the previous proposition. There is another request about the the starting point of 
each discrete solution that can be slightly different from the original $\mu_0$, but in this case is taken to be the same and so no problem 
arises. 
\begin{proposition} \label{curvalimite} Suppose $W$ satisfies \textbf{A - C}. There exist a sequence $\tau_n \searrow 0$, and a limit curve $\mu \in AC^2_{loc}
\bigr([0,+\infty);\mathcal{P}_2(\mathbb{R}^d)\bigr)$, such that \begin{equation*}
\mu_n(t):=\mu^{\tau_n}(t) \to \mu(t), \text{   narrowly as } n\to +\infty
\end{equation*}
for all $t \in [0,+\infty)$.
\end{proposition}

The last step of the procedure is to check that the limit curve obtained is a curve of maximal slope for $W$ according to definition 
\ref{lalorodef}. \\

Denoting by $\tilde{\mu}^n(t)$ the De Giorgi variational interpolation thanks to a technical lemma (see \cite{AGS} Lemma 3.2.2, "a priori 
estimates") the following energy inequality holds: \begin{equation}\label{apriori}
\mathcal{W}[\mu_0] \geq \frac{1}{2}\int_0^T\|\boldsymbol{v}^n(t)\|^2_{L^2(\mu^n(t))}dt+\frac{1}{2}\int_0^T|\partial \mathcal{W}| (\tilde{\mu}^n(t))^2dt + \mathcal{W}[\mu^n(T)]
\end{equation}
for all $T>0$, where on any interval $[(k-1)\tau_n,k\tau_n]$ the curve $\mu^n(t)$ is a Wasserstein geodesic connecting $\mu_{k-1}^{\tau_n}$ 
to $\mu_k^{\tau_n}$, and $\boldsymbol{v}^n$ is its velocity field. The continuity equation for $\mu^n$ holds and up to a subsequence both 
$\mu^n$ and $\tilde{\mu}^n$ narrowly converge to the same limit curve. A technical lemma is needed to pass to the limit the slope term in (ref{apriori}).

\begin{lemma}[Lower semi-continuity of the slope] \label{lowersemislope} \begin{equation*}
\liminf_{n \to \infty}\int_0^T|\partial\mathcal{W}|^2(\tilde{\mu}^n(t))dt \geq \int_0^T|\partial \mathcal{W}|^2(\mu(t))dt.
\end{equation*}
\proof 
By using the representation formula give in proposition \ref{classssubdif}, is needed to prove that  \begin{equation*}
\liminf_{n \to +\infty}\int_0^T \int_{\mathbb{R}^d} |k^n(x,t)|d\mu^n(t)(x)dt \geq \int_0^T \int_{\mathbb{R}^d} |k(x,t)|d\mu(t)(x)dt,
\end{equation*}
where \begin{equation*}
k^n(x,t):= \partial_0W \ast \mu^n(x,t), \qquad k(x,t):= \partial_0 W \ast \mu(x,t).
\end{equation*}
Without loss of generality, up to passing to a subsequence can be assumed that \begin{equation*}
\sup_n \int_0^T \int_{\mathbb{R}^d} |k^n(x,t)|d\mu^n(t)(x)dt < +\infty.
\end{equation*}
Hence, using theorem 5.4.4 of \cite{AGS} on the measure space $X:=\mathbb{R}^d \times [0,T]$ with the family of measures $\mu^n \otimes dt$, 
the lemma is proved if $k^n$ converges weakly to $k$, i.e. that for any vector field $\phi \in C_c^{\infty}(\mathbb{R}^d \times [0,T]; 
\mathbb{R}^d)$ \begin{equation} \label{easirecov}
\int_0^T \int_{\mathbb{R}^d} \phi(x,t) \cdot k^n(x,t) d\mu^n(t)(x)dt \to \int_0^T \int_{\mathbb{R}^d} \phi(x,t) k^n(x,t) d\mu^n(t)(x)dt
\end{equation}
as $n \to +\infty$. To show this is observed that the term on the left-hand side is given by \begin{multline*}
\int_0^T \int_{\mathbb{R}^d} \phi(x,t) k^n(x,t)d\mu^n(t)(x)dt = \\ 
= \int_0^T \int_{x \neq y} \phi(x,t) \cdot \nabla W(x-y)d\mu^n(t)(y)d\mu^n(t)(x)dt = \\
= \frac{1}{2} \int_0^T \int_{x \neq y} (\phi(x,t) - \phi(y,t)) \cdot \nabla W(x-y)d\mu^n(t)(y)d\mu^n(t)(x)dt,
\end{multline*}
where for the second equality is used the fact that $\nabla W$ is odd, so the expression in the integral can be symmetrized. \\

Thanks to the "a priori estimates", the sequence $\mu^n$ has uniformly bounded second moments. Therefore, with the linear growth control on 
the gradient of $W$ 
\begin{equation} \label{linegrocon}
|\nabla W(x)| \leq 2K + (2K + |\lambda|)|x| \qquad \text{ with }K>0,
\end{equation}
the function $(\phi(x,t)-\phi(y,t))\cdot \nabla W(x-y)$ is uniformly integrable with respect to $\mu^n\otimes \mu^n \otimes dt$, and 
(\ref{easirecov}) can be proved with weak convergence arguments.
\endproof
\end{lemma}

And now everything is prepared for the proof of existence of a solution. 
\begin{theorem}[Existence of curves of maximal slope] \label{exicurmaxslo}
Let $W$ satisfy the assumptions \textbf{A - C}. Then, there exists at least one curve of maximal slope for the functional $\mathcal{W}$, i.e. 
there exists at least one curve $\mu \in AC^2_{loc}([0,+\infty);\mathcal{P}_2(\mathbb{R}^d))$ such that the energy inequality 
\begin{multline} \label{eq223}
\mathcal{W}[\mu_0] \geq \frac{1}{2}\int_0^T\|\boldsymbol{v}(t)\|^2_{L^2(\mu^n(t))}dt + \\
+ \frac{1}{2}\int_0^T \int_{\mathbb{R}^d}\biggl|\int_{x \not= y}\nabla W(x-y)d\mu(t)(y)\biggr|^2d\mu(t)(x)dt + \mathcal{W}[\mu(T)]
\end{multline}
is satisfied, where $\boldsymbol{v}(t) \in L^2(\mu(t))$ is the minimal velocity field of $\mu$.
\proof
It is needed to prove that the curve $\mu(t)$ provided by proposition \ref{curvalimite} satisfies the desired condition. As a consequence 
of (\ref{apriori}) and of lemma \ref{lowersemislope}, showing that \begin{equation} \label{inequaliti224} 
\liminf_{n\to \infty} \frac{1}{2}\int_0^T\|v^n(t)\|^2_{L^2(\mu^n(t))}dt+\mathcal{W}[\mu^n(T)] \geq 
\frac{1}{2}\int_0^T\|v(t)\|^2_{L^2(\mu(t))}dt+\mathcal{W}[\mu(T)]
\end{equation}
all the remaining part of the proof of the convergence of the scheme to a solution goes through like in the case when $\mathcal{W}$ is 
lower semicontinuous with respect to the narrow topology, see \cite{AGS} chapter 3.
To prove the inequality, the solutions of $\partial_t \mu^n(t)+div(v^n(t)\mu^n(t))=0 \; \forall n$ are regularized as follows: 
\begin{eqnarray*}
& v^{n,\varepsilon}:=\frac{(v^n(t)\mu^n(t))\ast \eta_{\varepsilon}}{\mu^n(t) \ast \eta_{\varepsilon}} , \qquad & \mu^{n,\varepsilon}=\mu^n(t)\ast \eta_{\varepsilon}, \\
& v^{\varepsilon}:=\frac{(v(t)\mu(t))\ast \eta_{\varepsilon}}{\mu(t) \ast \eta_{\varepsilon}} , \qquad & \mu^{\varepsilon}=\mu(t)\ast \eta_{\varepsilon},
\end{eqnarray*} 
where $\eta_{\varepsilon}=\frac{1}{\varepsilon^d}\eta(\frac{\cdot}{\varepsilon}\in C^{\infty}(\mathbb{R}^d)$ is a smooth convolution kernel with support 
the whole $\mathbb{R}^d$, say a gaussian. This is done like in lemma \ref{appbyregcur}. Applying proposition \ref{repforfortheconequ} the 
measures $\mu^{n,\varepsilon}(t),\mu^{\varepsilon}(t)$ are given by the formula $\mu^{n,\varepsilon}(t)=(X^{n,\varepsilon}(t))_{\#}\mu_0$ and 
$\mu^{\varepsilon}(t)=(X^{\varepsilon}(t))_{\#}\mu_0$, where $X^{n,\varepsilon}(t)$ and $X^{n,\varepsilon}(t)$ denote the flows of $v^{n,\varepsilon}(t)$ 
and $v^{\varepsilon}(t)$ respectively, more precisely 
\begin{eqnarray*}
& \frac{d}{dt}X^{n,\varepsilon}(t,x)=v^{n,\varepsilon}(t,X^{n,\varepsilon}(t,x)), \qquad & X^{n,\varepsilon}(0,x)=x, \\
& \frac{d}{dt}X^{\varepsilon}(t,x)=v^{\varepsilon}(t,X^{\varepsilon}(t,x)), \qquad & X^{\varepsilon}(0,x)=x. 
\end{eqnarray*}
Now can be defined the map from $\mu^{\varepsilon}(T)$ to $\mu^{n,\varepsilon}(T)$ as $T^{\varepsilon}_n:=X^{n,\varepsilon}(T) \circ 
(X^{\varepsilon}(T))^{-1}$ and get \begin{align*}
& W_2^2(\mu^{\varepsilon}(T),\mu^{n,\varepsilon}(T)) \leq \int_{\mathbb{R}^d}|T_n^{\varepsilon}(x)-x|^2d\mu^{\varepsilon}(T)(x) = \\
& = \int_{\mathbb{R}^d}|X^{n,\varepsilon}(T) \circ (X^{\varepsilon}(T))^{-1} (x)- (X^{\varepsilon}(T))^{-1}(x) + \\
& \qquad \qquad \qquad \qquad \qquad \qquad \qquad \qquad +(X^{\varepsilon}(T))^{-1}(x) - x|^2d\mu^{\varepsilon}(T)(x) = \\ 
& = \int_{\mathbb{R}^d}\biggl| \int_0^T [v^{n,\varepsilon}(t,X^{n,\varepsilon}(T) \circ (X^{\varepsilon}(T))^{-1}(x)) -  \\
& \qquad \qquad \qquad \qquad \qquad  - v^{\varepsilon}(t,(X^{\varepsilon}(T)) \circ (X^{\varepsilon}(T))^{-1}(x)]dt \biggr|^2 d\mu^{\varepsilon}(T)(x) = \\ 
& = \int_{\mathbb{R}^d}\biggl| \int_0^T [v^{n,\varepsilon}(t,X^{n,\varepsilon}(t,x)) - v^{\varepsilon}(t,X^{\varepsilon}(t,x)]dt \biggr|^2 d\mu_0^{\varepsilon}(x).
\end{align*}
By H\"older's inequality and expanding the squares is obtained \begin{align} \label{eq225}
\nonumber & W_2^2(\mu^{\varepsilon}(T),\mu^{n,\varepsilon}(T)) \leq T \int_{\mathbb{R}^d} \int_0^T 
|v^{n,\varepsilon}(t,X^{n,\varepsilon}(t,x)) - v^{\varepsilon}(t,X^{\varepsilon}(t,x)|^2dt d\mu_0^{\varepsilon}(x)\\
\nonumber & \leq T\int_0^T \int_{\mathbb{R}^d}|v^{n,\varepsilon}(t,x)|^2d\mu^{n,\varepsilon}(t)(x)dt + T\int_0^T \int_{\mathbb{R}^d} |v^{\varepsilon}
(t,x)|^2d\mu^{\varepsilon}(t)(x)dt \\
& -2T\int_0^T \int_{\mathbb{R}^d} v^{n,\varepsilon}(t,X^{n,\varepsilon}(t,x))\cdot v^{\varepsilon}(t,X^{\varepsilon}(t,x)]d\mu_0(x)dt.
\end{align}
Thanks to lemma \ref{appbyregcur} the following inequality holds
\begin{equation} \label{eq226}
\int_0^T \int_{\mathbb{R}^d} |v^{n,\varepsilon}(t,x)|^2d\mu^{n,\varepsilon}(t)(x)dt \leq \int_0^T \int_{\mathbb{R}^d} |v^{n}(t,x)|^2d\mu^n(t)(x)dt 
\; \forall \varepsilon > 0.
\end{equation}
Moreover, thanks to the weak convergence of $(\mu^n (t), v^n (t)\mu^n (t))$ to $(\mu (t),$\\$ v(t)\mu(t))$, which
is a consequence of the linear growth control of the gradient of W in (\ref{linegrocon}) and the fact
that $\mu^{n,\varepsilon}(t)$ and $\mu^{\varepsilon}(t)$ are uniformly (in $n \in \mathbb{N}$) bounded away from zero on compact sets
of $\mathbb{R}^d$ , is proved that \begin{equation}\label{eq227}
v^{n,\varepsilon}(t)\to v^{\varepsilon} \text{ in }L^1([0,T],C^{\infty}_c(\mathbb{R}^d)).
\end{equation}
Indeed \begin{gather*}
D^{\alpha}[v^{n,\varepsilon}-v^{\varepsilon}]=\frac{D^{\alpha}\eta^{\varepsilon}\ast (v^n\mu^n)}{\mu^{n,\varepsilon}}-\frac{D^{\alpha}\eta^{\varepsilon}\ast (v\mu)}{\mu^{\varepsilon}} \\
=D^{\alpha}\eta^{\varepsilon}\ast (v^n\mu^n)\biggl( \frac{\mu^{\varepsilon}-\mu^{n,\varepsilon}}{\mu^{\varepsilon}\mu^{n,\varepsilon}} \biggr) + \frac{1}{\mu^{\varepsilon}} D^{\alpha}\eta^{\varepsilon}\ast (v\mu-v^n\mu^n)
\end{gather*}
and $v^n$ is uniformly bounded in $L^2 (\mu^n )$ with respect to $n$. Since the flows $X^{n,\varepsilon} (t)$ and
$X^{\varepsilon} (t)$ are globally defined (see proposition 2.1.1), (\ref{eq227}) easily implies
that for any $t\in [0,T]$ \begin{equation} \label{eq228}
X^{n,\varepsilon} (t)\to X^{\varepsilon} (t) \text{ locally uniformly on compact subset of } \mathbb{R}^d.
\end{equation}
This fact, together with the fact that $v^{n,\varepsilon}(t,X^{n,\varepsilon} (t))$ are uniformly bounded in $L^2(\mu_0 \otimes dt)$ thanks to 
(\ref{eq226}), implies that \begin{equation} \label{eq229}
\begin{split}
& \lim_{n \to \infty}\int_{\mathbb{R}^d}\int_0^T v^{n,\varepsilon}(t,X^{n,\varepsilon}(t,x))\cdot v^{\varepsilon}(t,X^{\varepsilon}(t,x))dtd\mu_0(x) \\
& = \int_{\mathbb{R}^d}\int_0^T |v^{\varepsilon}(t,X^{\varepsilon}(t,x))|^2dtd\mu_0(x)=\int_{\mathbb{R}^d}\int_0^T |v^{\varepsilon}(t,x)|^2d\mu_0(x)dt.
\end{split}
\end{equation}
To prove (\ref{eq229}), split the integral on the left-hand side as follows \begin{equation*} \begin{split}
& \int_{\mathbb{R}^d}\int_0^T v^{n,\varepsilon}(t,X^{n,\varepsilon}(t,x))\cdot v^{\varepsilon}(t,X^{\varepsilon}(t,x))dtd\mu_0(x)= \\
& = \int_{|x|>R}\int_0^T v^{n,\varepsilon}(t,X^{n,\varepsilon}(t,x))\cdot v^{\varepsilon}(t,X^{\varepsilon}(t,x))dtd\mu_0(x)= \\
& + \int_{|x| \leq R}\int_0^T v^{n,\varepsilon}(t,X^{n,\varepsilon}(t,x))\cdot v^{\varepsilon}(t,X^{\varepsilon}(t,x))dtd\mu_0(x)=: I_1+I_2. 
\end{split}
\end{equation*}
Now, thanks to (\ref{eq226})and the fact that $v^n$ is uniformly bounded in $L^2(\mu^n)$ with respect to $n$, can be done the following estimates 
\begin{equation*}
\begin{split}
I_1^2 & \leq \int_{\mathbb{R}^d}\int_0^T |v^{n,\varepsilon}(t,X^{n,\varepsilon}(t,x))|^2dtd\mu_0(x) \int_{|x|>R}\int_0^T |v^{\varepsilon}(t,X^{\varepsilon}(t,x))|^2dtd\mu_0(x) \\
& \leq \int_{|x|>R}\int_0^T |v^{\varepsilon}(t,X^{\varepsilon}(t,x))|^2dtd\mu_0(x) \\
\end{split}
\end{equation*}
for some constant $C$ independent on $n$. Hence, one can choose $R$ large enough such that
$|I_1|< \eta$ for an arbitrarily small $\eta>0$. On the other hand, (\ref{eq227}) and (\ref{eq228}) imply \begin{equation*}
I_2 \to \int_{|x|\leq R}\int_0^T|v^{\varepsilon}(t,x)|^2d\mu_0(x)dt
\end{equation*}
as $n\to +\infty$ and (\ref{eq229}) follows by letting $R \to +\infty$. Therefore, by combining (\ref{eq229}) with (\ref{eq225}) and (\ref{eq226}) is obtained that \begin{multline} \label{eq230}
\liminf_{n \to \infty}W_2^2(\mu^{\varepsilon}(T),\mu^{n,\varepsilon}(T))+2T\mathcal{W}[\mu^n(T)] \leq \liminf_{n \to \infty}T \biggl[ 2\mathcal{W}[\mu^n(T)] \biggr.   \\
 \biggl. + \int_0^T \int_{\mathbb{R}^d}|v^n(t,x)|^2d\mu^n(t)(x)dt -  \int_0^T \int_{\mathbb{R}^d}|v^{\varepsilon}(t,x)|^2d\mu^{\varepsilon}(t)(x)dt \biggr].
\end{multline}
We now claim that there exists a constant $C_0 > 0$, depending only on the convolution
kernel $\eta$, such that for any $\mu \in \mathcal{P}(\mathbb{R}^d)$ \begin{equation} \label{eq231}
W_2^2(\mu,\mu \ast \eta_{\varepsilon}) \leq C_0\varepsilon^2.
\end{equation}
Indeed is suffices to consider the transport plan $\gamma^{\varepsilon} \in \Gamma(\mu,\mu \ast \eta_{\varepsilon})$ defined as \begin{equation*}
\int_{\mathbb{R}^d \times \mathbb{R}^d}f(x,y) d\gamma^{\varepsilon}(x,y):=
\int_{\mathbb{R}^d \times \mathbb{R}^d}f(x,y) \eta_{\varepsilon}(y-x)dyd\mu(x) \; \forall f \in C^0_b(\mathbb{R}^d \times \mathbb{R}^d),
\end{equation*}
to get that \begin{equation*}
\int_{\mathbb{R}^d \times \mathbb{R}^d}|y-x|^2 d\gamma^{\varepsilon}(x,y)=
\int_{\mathbb{R}^d}|z|^2 \eta_{\varepsilon}dz=\varepsilon^2\int_{\mathbb{R}^d}|z|^2\eta(z)dz,
\end{equation*}
which proves (\ref{eq231}). Finally can be observed that \begin{equation}
\label{eq232}
\liminf_{\varepsilon \to 0}\int_0^T \int_{\mathbb{R}^d}|v^{\varepsilon}(t,x)|^2d\mu^{\varepsilon}(t)(x)dt \geq 
\int_0^T \int_{\mathbb{R}^d}|v(t,x)|^2d\mu(t)(x)dt.
\end{equation} 
Combining (\ref{eq230}) and (\ref{eq231}) \begin{equation*} 
\begin{split}
& \liminf_{n \to \infty}W_2^2(\mu(T),\mu^n(T))+2T\mathcal{W}[\mu^n(T)] \leq \liminf_{n \to \infty}T \biggl[ 2\mathcal{W}[\mu^n(T)] \biggr.   \\
& \biggl.+ \int_0^T \int_{\mathbb{R}^d}|v^n(t,x)|^2d\mu^n(t)(x)dt -  \int_0^T \int_{\mathbb{R}^d}|v^{\varepsilon}(t,x)|^2d\mu^{\varepsilon}(t)(x)dt  \biggr] + O(\varepsilon),
\end{split}
\end{equation*}
so, that letting $\varepsilon \to 0$, thanks to (\ref{eq232}), 
\begin{equation} \label{eq233}
\begin{split}
& \liminf_{n \to \infty}W_2^2(\mu(T),\mu^n(T))+2T\mathcal{W}[\mu^n(T)] \leq \liminf_{n \to \infty}T \biggl[ 2\mathcal{W}[\mu^n(T)] \biggr. \\
& \biggl. + \int_0^T \int_{\mathbb{R}^d}|v^n(t,x)|^2d\mu^n(t)(x)dt  -  \int_0^T \int_{\mathbb{R}^d}|v(t,x)|^2d\mu(t)(x)dt \biggr].
\end{split} 
\end{equation}
Moreover, in view of the lower semicontinuity of the slope, \begin{equation} \label{eq234}
\liminf_{n \to \infty}W_2^2(\mu(T),\mu^n(T))+2T\mathcal{W}[\mu^n(T)] \geq 2T\mathcal{W}[\mu(T)]
\end{equation}
for $T$ small enough. Combining (\ref{eq234}) with (\ref{eq233}), is obtained that (\ref{inequaliti224}) holds provided
$T$ is sufficiently small (but independent on the initial datum $\mu_0$), and this allows to prove
the existence of a curve of maximal slope on a small time interval $[0,T]$. Iterating now
the construction via minimizing movements on $[T,2T], \; [2T,3T]$ and so on, and adding
the energy inequalities (\ref{eq223}) on each time interval, the desired result is finally proved.
\endproof
\end{theorem}

To summarize, the three notions of solutions are shown to be equivalent because, thanks to previous results, the hypothesis of theorem 
\ref{eqcmsgf} are satisfied. With that result and the existence of curves of maximal slope the following theorem holds: 
\begin{theorem}[Existence of the gradient flow]
Let $W$ satisfy the assumptions \textbf{A - C}. Given any $\mu_0 \in \mathcal{P}_2(\mathbb{R}^d)$, then there exists a gradient flow solution, 
i.e. a curve $\mu \in AC^2_{loc}([0,+\infty);\mathcal{P}_2(\mathbb{R}^d))$ satisfying \begin{equation*}
\begin{split}
& \frac{\partial \mu(t)}{\partial t} + div(\boldsymbol{v}(t)\mu(t))=0 \text{ in } \mathcal{D}'([0,+\infty) \times \mathbb{R}^d), \\
& \boldsymbol{v}(t)=-\partial^0\mathcal{W}[\mu(t)]=-\partial^0W \ast \mu(t), \\
& \| \boldsymbol{v}(t) \|_{L^2(\mu(t))}=|\mu|'(t) \qquad a.e. \; t>0,
\end{split}
\end{equation*}
with $\mu(0)=\mu_0$. Moreover, the energy identity \begin{equation*}
\int_a^b\int_{\mathbb{R}^d}|\boldsymbol{v}(t,x)|^2d\mu(t)(x) + \mathcal{W}[\mu(b)] = \mathcal{W}[\mu(a)]
\end{equation*}
holds for all $0 \leq a \leq b < +\infty$.
\end{theorem}
Moreover, thanks to remark \ref{qweretyty}, he following result about unicity of the solution follows readily from \cite{AGS}, theorem 11.1.4.
\begin{theorem}
Let W satisfy the assumptions \textbf{A - C}. Given two gradient flow solutions $\mu^1(t)$ and $\mu^2(t)$ in the sense of the theorem 
above, we have
\begin{equation}
W_2(\mu^1(t),\mu^2(t)) \leq e^{-\lambda t}W_2(\mu_0^1,\mu_0^2)
\end{equation}
for all $t \geq 0$. In particular, the gradient flow solution starting from any given $\mu_0 \in \mathcal{P}_2(\mathbb{R}^d)$
is unique. Moreover, this solution is characterized by a system of evolution variational inequalities:
\begin{equation} \label{systemevi}
\frac{1}{2}\frac{d}{dt}W_2^2(\mu(t),\sigma)+\frac{\lambda}{2}W_2^2(\mu(t),\sigma) \leq \mathcal{W}[\sigma] - \mathcal{W}[\mu(t)] \qquad \text{a.e.} t>0,
\end{equation}
for all $\sigma \in \mathcal{P}_2(\mathbb{R}^d)$.
\end{theorem}

\subsubsection{Finite time aggregation}

Here we take in account for an attractive non-Osgood potential $W$ that, in addition to \textbf{A - C}, satisfies the finite time 
blow-up condition: 
\begin{description}
\item[D] $W$ is radial, i.e. $W (x) = w(|x|)$, $W \in C^2(\mathbb{R}^d\setminus \{ 0\})$ with $w' (r) > 0$
for $r > 0$ and satisfying the following monotonicity condition: either \textbf{(a)} $w'(0+ ) >
0$, or \textbf{(b)} $w'(0+ ) = 0$ with $w''(r)$ monotone decreasing on an interval $(0, \varepsilon_0 )$.
Moreover, the potential satisfies the integrability condition \begin{equation} \label{lequa18}
\int_0^{\varepsilon_1}\frac{1}{w'(r)}dr< +\infty, \qquad \text{ for some } \varepsilon_1>0.
\end{equation}
\end{description} 
The condition of monotonicity of $w''(r)$ is not too restrictive. It is 
actually automatically satisfied by any potential which satisfies (\ref{lequa18}) and whose second
derivative does not oscillate badly at the origin. Examples of this type of potentials are the ones
having a local behavior at the origin like $w'(r)\simeq r^{\alpha}$ with $0 \leq \alpha < 1$ or $w'(r)\simeq r\text{log}^2 r$. \\

Here are reported two important result about aggregation of solutions (see article \cite{IL} for a more detailed dissertation). 
The first is that, in presence of a non-Osgood potential, solutions tend to aggregate in finite time.
\begin{proposition}[Finite time total collapse] \label{finitimetot}
Assume W satisfies \textbf{A - C} and \textbf{D}. Let $\mu(t)$ denote the unique gradient flow solution starting from the probability 
measure $\mu_0$ with center of mass \begin{equation*}
x_c:=\int_{\mathbb{R}^d}xd\mu_0,
\end{equation*}
supported in $\overline{B}(x_c,R_0)$. Then there exists $T^{\ast}$, depending only on $R_0$, such that $\mu(t)=\delta_{x_c}$ 
for all $t \geq T^{\ast}$.
\end{proposition}

The second proposition shows that, if we start with a measure which has some atomic part, then the atoms can only increase their mass.
\begin{proposition}[Dirac delta can only increase mass] \label{finitimetot2}
Let $\mu(t)$ denote the unique gradient flow solution starting from the probability 
measure $\sum_{i=1}^N m_i\delta_{x_i^0}+\nu_0$, and define the curves $t \mapsto x_i(t), \; i=1,\cdots,N$, as the solution of the ODE 
\begin{equation}
\dot{x}_i(t)=-(\partial^0W\ast \mu(t))(x_i(t)). 
\end{equation} 
Then $\mu(t)\geq \sum_{i=1}^N m_i\delta_{x_i(t)}$ for all $t \geq 0$, with possibly $x_i(t)=x_j(t)$ for some $t>0,i\neq j$.
\end{proposition}

\clearpage

\lhead[\thepage]{\small{\textsc{A non smooth repulsive interaction potential}}}
\rhead[\small{\textsc{A non smooth repulsive interaction potential}}] {\thepage}	

\clearpage{\pagestyle{empty}\cleardoublepage}

\chapter{A non smooth repulsive interaction potential}
\markboth{A non smooth repulsive interaction potential}{A non smooth repulsive interaction potential}

In this chapter is consider again a mass distribution of particles $\mu_t \in \mathcal{P}_2(\mathbb{R}^d), \; \mu\geq 0$, evolving 
under the action of a continuous interaction potential $W$. The associated interaction energy is defined as \begin{equation}
\mathcal{W}[\mu]:=\frac{1}{2}\int_{\mathbb{R}^d\times \mathbb{R}^d}W(x-y)d\mu(x)d\mu(y).
\end{equation}
The equation in $\mathbb{R}^d$ takes the form \begin{equation} 
\label{theeeeeq}
\frac{\partial \mu}{\partial t}=\nabla \cdot \left[ \left( \nabla \frac{\delta \mathcal{W}}{\delta \mu} \right) \mu \right]=\nabla \cdot 
\left[ \left( \nabla W \ast \mu \right) \mu \right] \; t\geq 0,
\end{equation}
with the initial condition $\mu(0)=\mu_0$.

We assume that the potential satisfies the following properties:
\begin{description}
\item[A] $W$ is continuous, $W(x)=W(-x)$, $W(0)=0$ and $W \in C^1(\mathbb{R}^d\setminus \{0\})$. 
\item[B] There exists a constant $C > 0$ such that \begin{equation*}
W(z) \leq C(1+|z|^2) \qquad \forall z \in \mathbb{R}^d.
\end{equation*}
\item[C'] $\exists \lambda',\lambda'' \geq 0$ such that $W(x)+\frac{\lambda''}{2}|x|^2+\lambda'|x|$ is convex.  
\end{description}

The assumption \textbf{C'}, without loss of generality, can be stated with the same constant $\lambda=max(\lambda',\lambda'')$ for each 
term. The existence of a solution has been demonstrated in \cite{RAOUL} in the one dimensional case for a starting probability measure $\rho_0 \in W^{2,\infty}$ 
with compact support, thanks to an "a priori estimate for $\rho$" (is shown that $\|\rho(t,\cdot)\|_{W^{2,\infty}} \leq \|\rho^0\|_{W^{2,\infty}}$) 
and with an inductive scheme to obtain the solution. Here every probability measure are admitted and is made use of the 
minimizing movement scheme. \\

The idea of working under this assumption comes naturally looking at the equation (in the 1 dimensional case) with the simple choice of 
the potential $W(x)=|x|$; in this case the solution exists and it is unique, and it has been shown by \cite{IL} (see also 
 propositions \ref{finitimetot} and \ref{finitimetot2}) that particles aggregate in finite time. Moreover all the distributions with 
 the same center of mass aggregate in the same final distribution: the Dirac delta of the center of mass. What happens if 
the time reversal is applied to the equation? The answer is not trivial since infinitely many distributional solutions are allowed. It is 
easy to check that the equation generated by the time reversal corresponds to the gradient flow driven by the potential $W(x)=-|x|$. 
It has to be noticed that all the potential satisfying the assumption \textbf{C'} also satisfy the weaker assumption \textbf{C}. 
Now the subdifferential of $W$ is empty at the origin ($\partial W(0) = \varnothing $) and the crucial $\lambda$-convexity property no 
more holds. Our objective is to show that a nice solution can still be found by adopting the gradient flow representation.  \\

Three main steps characterize the structure of this chapter. \begin{itemize}
\item We first compare equation \ref{theeeeeq} with an ODE system.
\item We study the MMS.
\item We characterize the Wasserstein subdifferential (achieved only in the one dimensional case) and prove existence and uniqueness of 
the gradient flow.
\end{itemize}

\section{The ODE system}

The starting point is the comparison of the continuity equation with respect to an ODE system. \\

Let $x_i(t),i=1,\cdots,N,$ be $C^1$-solutions (at least for a short initial time interval) of the ODE system  \begin{equation}
\dot{x}_1=-\sum_{j \not= i}m_j\nabla W(x_i-x_j), \qquad i=1,\cdots ,N,
\end{equation}
with $m_i>0$ and $\sum_im_i=1$. Then it is straightforward to check that $\mu(t):=\sum_{i=1}^Nm_i\delta_{x_i(t)}$ is a solution of the continuity 
equation in the distributional sense. Conversely, if $\mu(t)$ of the above form solves the PDE and $x_i(t)$ are $C^1$ curves for $i=1,\cdots,N$, 
then $x_i(t)$ solve the ODE system. 

If the particles collide the solutions of the PDE can still be represented by solving an ODE after the collision. A sketch of the proof 
(following \cite{IL}, remark 2.10) is reported here. We consider absolutely continuous solutions of \begin{gather}
\dot{x}_i=-\sum_{j \in C(i)}m_j \nabla W(x_i-x_j), \qquad i=1,\cdots,N, \\
\text{with }C(i):=\{ j \in \{ 1,\cdots, N \} : j \neq i, x_j(t) \neq x_i(t) \}
\end{gather}
More precisely we consider the solutions of the associated integral equation. If $C(i)$ is
empty, then all particles have collapsed to a single particle. We then define the right hand
side to be zero (i.e. we define the sum over empty set of indexes to be zero). The right
hand side of this ODE system is bounded and Lipschitz-continuous in space on short time
intervals. Thus the ODE system has a unique Lipschitz-continuous solutions on short time
intervals. The estimate (\ref{linegrocon}) then implies that the solutions are uniformly bounded. Note that
the solutions are Lipschitz (in time) on bounded time intervals. Also note that collisions
of particles can occur, but that we do not relabel the particles when they collide. Since
the number of particles is $N$ there exist $0\leq k \leq N-1$ times $0=:T_0 <T_1<T_2<\cdots < T_k$ 
at which collisions occur. Note that $\mu(t)=\sum_{i=1}^Nm_i\delta_{x_i(t)}$ is a solution
of the PDE on the time intervals $[T_l,T_{l+1})$. Furthermore, the Lipschitz continuity of $x_i$
 implies that $\mu$ is an absolutely continuous curve in $\mathcal{P}_2(\mathbb{R}^d)$. It is then straightforward
to verify that $\mu$ is a weak solution according to definition \ref{definiz1punto2}. Since the solution to the
PDE is unique (at least when $W$ is $\lambda$-convex) the converse claim also holds. \\

When $W$ is not $\lambda$-convex, e.g. in the case $W(x)=-|x|$, it is easy to construct a simple example, showing existence of
many solutions to (\ref{theeeeeq}) with the same initial datum.
\begin{example} \label{theesempio}
Consider a starting distribution with all particles concentrated in a point, $x_i(0)=x_0\;\forall i$ with $W(x)=-|x|$. The ODE system 
become \begin{equation}
\dot{x}_i=\sum_{j \not= i}m_j sign(x_i-x_j).
\end{equation} 
It's a straightforward calculation to check that a solution can be: 
\begin{itemize}
\item $ x_i=x_0 \ \forall i$
\item $ x_i(t)=x_0+(-1)^i\frac{t}{2} \ i=1,2 \ , \ m_i=\frac{1}{2}$
\item $ x_i(t)=x_0+sign(i-2)\frac{2t}{3} \ i=1,2,3 \ , \ m_i=\frac{1}{3}$
\item $ x_i=x_0+(-1)^i\frac{t-\overline{t}}{2}H(t-\overline{t}), \ i=1,2 \ , \ m_i=\frac{1}{2} \ , \ \overline{t}>0$
\item $ \dots $
\end{itemize}
So there are infinite distributional solutions for the PDE related to the solutions of the ODE system. In fact this example represents the 
time reversal of the problem with the attractive potential $W(x)=|x|$ . Difficulties arise after the  
complete aggregation of the initial distribution in the center of mass. So it is impossible to reconstruct in an unique way the starting distribution.
\end{example}

Starting from this example, we look if we can select a "preferred" solution to (\ref{theeeeeq}) as a 
gradient flow equation. A similar study has been done in \cite{FELRAO,FELRAO2} 
while studying the stability of the steady state solution in the one dimensional case. In these articles the solution is considered 
in a distributional sense, here we adopt the gradient flow point of view.

\section{Minimizing movement scheme}

The proof follows the steps of \cite{IL}, apart from the characterization of the subdifferential, which has to be proved in a different way. 
The existence of the minimizing sequence is the first result that can be obtained . \\

Here is recalled the definition of the minimizing movement scheme and then proved the weak lower semi-continuity property: given an 
initial measure $\mu_0 \in \mathcal{P}_2(\mathbb{R}^d)$ and a time step 
$\tau > 0 $, a sequence $\mu_k^{\tau}$ is recursively defined by $\mu_0^{\tau}=\mu_0$ and \begin{equation}
\mu_{k+1}^{\tau}\in \text{argmin}_{\mu \in \mathcal{P}_2(\mathbb{R}^d)}  \biggl\{ \mathcal{W}[\mu]+\frac{1}{2\tau}W_2^2(\mu_k^{\tau},\mu) \biggr\}.
\end{equation}

\begin{lemma}[Weak lower semi-continuity of the penalized interaction energy] 
Suppose $W$ satisfies \textbf{A - C'}. Then, for a fixed $\overline{\mu}\in \mathcal{P}_2(\mathbb{R}^d)$, the penalized 
interaction energy functional \begin{equation}
\mathcal{P}_2(\mathbb{R}^d) \ni \mu \mapsto \mathcal{W}[\mu]+\frac{1}{2\tau}W_2^2(\mu,\overline{\mu})
\end{equation}
is lower semi-continuous with respect to the narrow topology of $\mathcal{P}(\mathbb{R}^d)$ for all $\tau > 0 $ such that $12 \tau \lambda^- 
\leq 1$, where $\lambda^-:=\text{max} \left\{0, \lambda',\lambda'' \right\}$. 
\proof
Let $\{\mu_n\} \subset \mathcal{P}_2(\mathbb{R}^d)$ such that $\lim_{n\to \infty}\mu_n=\mu_{\infty}$ narrowly. Is needed to prove that 
\begin{equation}
\liminf_{n\to \infty}\biggl[\mathcal{W}[\mu_n]+\frac{1}{2\tau}W_2^2(\mu_n,\overline{\mu})\biggr] \geq \mathcal{W}[\mu_{\infty}]+
\frac{1}{2\tau}W_2^2(\mu_{\infty},\overline{\mu}).
\end{equation}
The conditions $W(0)=0$, $\lambda = -\text{max}\left\{\lambda',\lambda'' \right\} \leq 0$ and the \textbf{C'} property of $W$ implies that \begin{equation} \label{lowerest}
\begin{split}
W(x-y) & \geq \frac{\lambda}{2}|x-y|^2 + \lambda |x-y| \geq \lambda(|x|^2+|y|^2) + \lambda(|x|+|y|) \\
& \geq \lambda (|x|^2+|y|^2+\frac{|x|^2}{2} + \frac{|y|^2}{2} + 1) \geq \frac{3\lambda}{2}(|x|^2+|y|^2) + \lambda
\end{split}
\end{equation}
and so \begin{equation*}
h(x,y):=W(x-y)-\frac{3\lambda}{2}(|x|^2+|y|^2) - \lambda
\end{equation*}
is a nonnegative continuous function. Therefore, \begin{equation} \begin{split}
\mathcal{W}[\mu_n]+\frac{1}{2\tau}W_2^2(\mu_n,\overline{\mu})& = \frac{3\lambda}{2}\int_{\mathbb{R}^d\times \mathbb{R}^d}
(|x|^2+|y|^2)d\mu_n(x)d\mu_n(y) +  \lambda\\
& + \int_{\mathbb{R}^d\times \mathbb{R}^d}h(x,y)d\mu_n(x)d\mu_n(y) + \frac{1}{2\tau}W_2^2(\mu_n,\overline{\mu}) \\
& = \int_{\mathbb{R}^d\times \mathbb{R}^d}h(x,y)d\mu_n(x)d\mu_n(y) \\
& + \frac{1}{2\tau}W_2^2(\mu_n,\overline{\mu}) + 3\lambda\int_{\mathbb{R}^d}|x|^2d\mu_n(x) +  \lambda 
\end{split}
\end{equation}
Since $h \geq 0$ \begin{equation}
\liminf_{n \to \infty}\int_{\mathbb{R}^d\times \mathbb{R}^d}h(x,y)d\mu_n(x)d\mu_n(y) \geq \int_{\mathbb{R}^d\times \mathbb{R}^d}h(x,y)d\mu_{\infty}(x)d\mu_{\infty}(y).
\end{equation}
Therefore, to get the desired assertion it suffices to prove that \begin{equation} \label{suffices} \begin{split}
& \liminf_{n \to \infty} \frac{1}{2\tau}W_2^2(\mu_n,\overline{\mu}) + 3\lambda\int_{\mathbb{R}^d}|x|^2d\mu_n(x) \\
& \geq \frac{1}{2\tau}W_2^2(\mu_{\infty},\overline{\mu}) + 3\lambda\int_{\mathbb{R}^d}|x|^2d\mu_{\infty}(y).
\end{split}
\end{equation}
Now, let $\boldsymbol{\gamma}_n \in \Gamma_0(\overline{\mu},\mu_n)$. Then, \begin{equation}
\frac{1}{2\tau}W_2^2(\mu_n,\overline{\mu}) + 3\lambda\int_{\mathbb{R}^d}|x|^2d\mu_n(x) =
\int_{\mathbb{R}^d\times \mathbb{R}^d}\biggl(\frac{1}{2\tau}|x-y|^2+3\lambda|y|^2\biggr)d\boldsymbol{\gamma}_n(x,y).
\end{equation}
Stability of optimal transportation plans (see \cite{CV} theorem 5.20) implies that there exists a subsequence, that is here assumed to be 
the whole sequence, such that $\boldsymbol{\gamma}_n$ converges narrowly to an optimal plan $\boldsymbol{\gamma}_{\infty} \in 
\Gamma_0(\overline{\mu},\mu_{\infty})$. As a consequence of \begin{equation*}
\int_{\mathbb{R}^d \times \mathbb{R}^d}|x|^2d\boldsymbol{\gamma}_n(x,y)=\int_{\mathbb{R}^d}|x|^2d\overline{\mu}(x)
= \int_{\mathbb{R}^d \times \mathbb{R}^d}|x|^2d\boldsymbol{\gamma}_{\infty}(x,y)
\end{equation*}
and the elementary inequality $|y|^2 \leq 2|x-y|^2 + 2|x|^2$ which implies \begin{equation} \label{elementary}
\frac{1}{2\tau}|x-y|^2+3\lambda|y|^2+\frac{1}{2\tau}|x|^2 \geq \biggl(\frac{1}{4\tau}+3\lambda\biggr)|y|^2 \geq 0 \; \text{ if } 
\tau \leq \frac{1}{12\lambda^-},
\end{equation}
can be obtained the inequality \ref{suffices}: \begin{multline*}
\liminf_{n \to \infty} \int_{\mathbb{R}^d\times \mathbb{R}^d}\biggl(\frac{1}{2\tau}|x-y|^2+3\lambda|y|^2\biggr)d\boldsymbol{\gamma}_n(x,y)  \\
= -\frac{1}{2\tau} \int |x|^2d\boldsymbol{\gamma}_{\infty}(x,y) + 
\liminf_{n \to \infty} \int \biggl(\frac{1}{2\tau}|x-y|^2+3\lambda|y|^2 + \frac{1}{2\tau}|x|^2\biggr)d\boldsymbol{\gamma}_n(x,y)  \\
\geq -\frac{1}{2\tau} \int |x|^2d\boldsymbol{\gamma}_{\infty}(x,y) + 
\int \biggl(\frac{1}{2\tau}|x-y|^2+3\lambda|y|^2 + \frac{1}{2\tau}|x|^2\biggr)d\boldsymbol{\gamma}_{\infty}(x,y)  \\
= \int_{\mathbb{R}^d \times \mathbb{R}^d} \biggl(\frac{1}{2\tau}|x-y|^2+3\lambda|y|^2 \biggr)d\boldsymbol{\gamma}_{\infty}(x,y). 
\end{multline*} 
\endproof
\end{lemma}

The next proposition can be proved thanks to the results presented in section \ref{subsection}: 
\begin{proposition}[Existence of minimizers] Suppose $W$ satisfies \textbf{A - C'}. Then, there 
exists $\tau_0 > 0$ depending only on $W$ such that, for all $0 < \tau < \tau_0$ and for all given $\overline{\mu} \in \mathcal{P}_2(\mathbb{R}^d)$, 
there exists $\mu_{\infty} \in \mathcal{P}_2(\mathbb{R}^d)$ such that \begin{equation}
\mathcal{W}[\mu_{\infty}]+\frac{1}{2\tau}W_2^2(\mu_{\infty},\overline{\mu}) = \min_{\mu \in \mathcal{P}_2(\mathbb{R}^d)} \biggl\{ 
\mathcal{W}[\mu]+\frac{1}{2\tau}W_2^2(\mu,\overline{\mu}) \biggr\}.
\end{equation}
\proof 
\underline{Compactness}: Given a measure $\overline{\mu} \in \mathcal{P}_2(\mathbb{R}^d)$ and a time step $\tau >0$, it's considered 
a minimizing sequence $(\mu_n) \subset \mathcal{P}_2(\mathbb{R}^d)$, i.e. \begin{equation*}
\inf_{\mu \in \mathcal{P}_2(\mathbb{R}^d)} \biggl\{ \mathcal{W}[\mu]+\frac{1}{2\tau}W_2^2(\mu,\overline{\mu}) \biggr\} = 
\lim_{n \to \infty} \biggl\{ \mathcal{W}[\mu_n]+\frac{1}{2\tau}W_2^2(\mu_n,\overline{\mu}) \biggr\}
\end{equation*}
Since $\mu_n$ is a minimizing sequence, the following inequality holds: \begin{equation}
\mathcal{W}[\mu_n] + \frac{1}{2\tau}W_2^2(\mu_n,\overline{\mu}) \leq C_1
\end{equation} 
for some constant $C_1$. Then, thanks to the inequalities (\ref{lowerest}),(\ref{elementary}), the following relations holds: 
\begin{gather*}
\mathcal{W}[\mu_n] + \frac{1}{2\tau}W_2^2(\mu_n,\overline{\mu}) \geq \int_{\mathbb{R}^d}3\lambda |y|^2d\mu_n(y)  
+ \frac{1}{2\tau}\int_{\mathbb{R}^d \times \mathbb{R}^d}|x-y|^2d\boldsymbol{\gamma}_n(x,y) + \lambda \\
\geq \int_{\mathbb{R}^d \times \mathbb{R}^d} \biggl(6\lambda + \frac{1}{2\tau}\biggr)|x-y|^2 d\boldsymbol{\gamma}_n(x,y) + 
\int_{\mathbb{R}^d}6\lambda|x|^2d\overline{\mu}(x) + \lambda \\
\text{and so} \\
\biggl(6\lambda + \frac{1}{2\tau}\biggr)W_2^2(\mu_n,\overline{\mu}) \leq 
C_1 - \lambda - \int_{\mathbb{R}^d}6\lambda|x|^2d\overline{\mu}(x)
\end{gather*}
The right side of the equation is constant and independent on $n$. For $\tau$ small enough $(6\lambda + \frac{1}{2\tau})$ is positive, so 
the Wasserstein distance is uniformly bounded with respect to $n$. Prokhorov's compactness theorem (see \ref{prok}) then implies that 
the sequence $\{\mu_n\}_n$ is tight. \\

\underline{Coercivity}: Is needed to prove that \begin{equation*}
\liminf_{n \to \infty}\biggl[\mathcal{W}[\mu_n]+ \frac{1}{2\tau}W_2^2(\mu_n,\overline{\mu})\biggr] \geq C_0 W_2^2(\mu_{\infty},\overline{\mu})-C_1
\end{equation*}
for some positive constant $C_0,C_1$ independent on $n$. The proof is similar to the previous step, the only difference is that has to 
be applied the infimum limit to both side of the inequality. \\

\underline{Passing to the limit by lower semi-continuity}: this is a consequence of the previous lemma.
\endproof 
\end{proposition}

The convergence to a limit curve $\mu(t)\in AC^2_{loc}\bigr([0,+\infty);\mathcal{P}_2(\mathbb{R}^d)\bigr)$ follows exactly as done in 
previous chapter up to the definition of the constant piecewise interpolation.

\begin{proposition} Suppose $W$ satisfies \textbf{A - C'}. There exist a sequence $\tau_n \searrow 0$, and a limit curve $\mu \in AC^2_{loc}
\bigr([0,+\infty);\mathcal{P}_2(\mathbb{R}^d)\bigr)$, such that \begin{equation*}
\mu_n(t):=\mu^{\tau_n}(t) \to \mu(t), \text{   narrowly as } n\to +\infty
\end{equation*}
for all $t \in [0,+\infty)$.
\end{proposition}

\section{Characterization of the subdifferential}

The characterization of the subdifferential is an important tool for proving the existence of a solution. We first consider the general 
case  of arbitrary space dimension $d \geq 1$, giving only a partial result. The next subsection will be focused 
on the one dimensional case, that exhibits a remarkable regularity property, essential to solve the problem.

\subsection{Dimension $d \geq 1$}

Here we prove just a first property satisfied by the subdifferential. Notice that we are considering only one implication 
and we do not claim that any element satisfying (\ref{eluieluiesi}) belongs to $\partial \mathcal{W} [\mu]$
\begin{proposition} \label{proposiminimal}
If $\mu \in\mathcal{P}_2(\mathbb{R}^d)$ and $\partial \mathcal{W}(\mu) \neq \varnothing$, then \begin{equation} \label{eluieluiesi}
k(x)=\int_{x \not= y}\nabla W(x-y)d\mu(y)
\end{equation}
 is the element of minimal norm of $\partial \mathcal{W}[\mu]$.
\proof
Fixed a vector field $\xi \in C_c^{\infty}(\mathbb{R}^d,\mathbb{R}^d)$ and observing that \\ $W(x-z+t(\xi(x)-\xi(z)))=W(x-z)=0$ when $x=z$,
 can be shown that
\begin{equation} \label{proposi}
\begin{split}
\nonumber & \lim_{t \to 0} \frac{\mathcal{W}[(id+t\xi)_{\#}\mu]-\mathcal{W}[\mu]}{t}\\
\nonumber & = \lim_{t\to 0}\frac{1}{2}\int_{\mathbb{R}^d \times \mathbb{R}^d} \frac{W((x-z)+t(\xi(x)-\xi(z)))-W(x-z)}{t} d\mu(x)d\mu(z) \\
\nonumber & = \lim_{t\to 0}\frac{1}{2}\int_{x \not= z} \frac{W((x-z)+t(\xi(x)-\xi(z)))-W(x-z)}{t} d\mu(x)d\mu(z) \\
\nonumber & = \frac{1}{2}\int_{x \not= z}\nabla W(x-z) \cdot (\xi(x)-\xi(z)) d\mu(x)d\mu(z) \\
& = \int_{\mathbb{R}^d}k(x)\cdot \xi(x)d\mu(x).
\end{split}
\end{equation}
Hence, since the definition of slope easily implies \begin{equation}
\liminf_{t \searrow 0}\frac{\mathcal{W}[(id+t\xi)_{\#}\mu]-\mathcal{W}[\mu]}{W_2((id+t\xi)_{\#}\mu,\mu)} \geq -|\partial \mathcal{W}|(\mu),
\end{equation}
using (\ref{proposi}) and the upper bound \begin{equation*}
\limsup_{\epsilon \searrow 0}\frac{W_2((id+\epsilon \xi)_{\#}\mu,\mu)}{\epsilon} \leq \|\xi\|_{L^2(\mu)},
\end{equation*}
is obtained
\begin{equation*} \begin{split}
\int_{\mathbb{R}^d}k(x)\cdot \xi(x)d\mu(x) & \geq -|\partial \mathcal{W}|(\mu)\liminf_{t \to 0} 
\frac{W_2((id+t \xi)_{\#}\mu,\mu)}{t} \\ & \geq -|\partial \mathcal{W}|(\mu)\|\xi\|_{L^2(\mu)}.
\end{split}
\end{equation*}
Changing $\xi$ with $-\xi$ gives \begin{equation*}
\biggl|\int_{\mathbb{R}^d}k(x)\cdot \xi(x)d\mu(x)\biggr| \leq |\partial \mathcal{W}| \| \xi \|_{L^2(\mu)},
\end{equation*}
so the arbitrariness of $\xi$ implies that $\| k \|_{L^2(\mu)} \leq |\partial \mathcal{W}|(\mu)$, and therefore $k$ is the element 
of minimal norm.
\endproof
\end{proposition} 

The problem is that the subdifferential may be empty. Consider, e.g., $d=1$ and $W(x)=-|x|$. Choosing $\mu=\delta_0$ the previous formula 
yields $k\equiv 0$. But we will show that $\mathcal{W}$ is convex, so $0\notin \partial \mathcal{W}[\mu]$ since otherwise $\mathcal{W}$ 
would have  minimum (not a maximum!) in $\delta_0$. \\

This difficulty could be circumvented by introducing a sort of relaxed version of $\partial \mathcal{W}$ (see e.g. \cite{AGS}), with the 
help of results like lemma \ref{illemmone}, which proves the existence of a subdifferential in the minimum point of the MMS.

\subsection{Case d=1}

In this section is analyzed the one dimensional case with the tools provided by section \ref{altrasubse}. The change of variable is applied  
 to the functional with the monotone rearrangement, obtaining: 
\begin{equation} \label{changevarutilli}
\mathcal{W}[\mu]=\int_{\mathbb{R} \times \mathbb{R}}W(x-y)d\mu(x)d\mu(y)=\int_0^1\int_0^1 W(X_{\mu}(s)-X_{\mu}(r))dsdr.
\end{equation}
Notice that the potential $W$ can be divided into the sum of two contributions 
$W(x)=\overline{W}(x)+\eta |x|$ with $\overline{W}$ satisfying the assumptions \textbf{A - C} and of class $C^1$, so that $\eta$ is 
uniquely determined. With this notation the 
functional itself turns out to be the sum of two terms: \begin{equation} \label{defdielle}
\mathcal{W}[\mu]=\overline{\mathcal{W}}+\mathcal{L}=\int_{\mathbb{R}\times \mathbb{R}}\overline{W}d\mu(x)d\mu(y)+\eta 
\int_{\mathbb{R}\times \mathbb{R}} |x-y|d\mu(x)d\mu(y).
\end{equation}
The subdifferential for $\overline{\mathcal{W}}$ is well defined in the previous chapter (see proposition \ref{classssubdif}) and here 
is denoted by $\overline{k}(x)$. 
Here we focus on the second contribution. \begin{proposition}
The functional $\mathcal{L}$ defined in (\ref{defdielle}) is geodesically convex for every $\eta \in \mathbb{R}$.
\proof
Given $\mu^1$ and $\mu^2$, is defined $\mu_{\theta}^{1 \to 2}=((1-\theta)\pi^1+\theta \pi^2)_{\#}\boldsymbol{\mu}$ with 
$\boldsymbol{\mu} \in \Gamma_0(\mu^1,\mu^2)$. Then, thanks to (\ref{changevarutilli}) \begin{equation}
\mathcal{L}[\mu_{\theta}^{1 \to 2}]=\eta\int_{\mathbb{R}\times \mathbb{R}} |x-y|d\mu_{\theta}(x)d\mu_{\theta}(y) = 
\eta \int_0^1 \int_0^1 |X_{\theta}(s) - X_{\theta}(r)|dsdr
\end{equation}
with \begin{equation}
X_{\theta}=(1-\theta)X_1+\theta X_2.
\end{equation}
Then, thanks to the monotonicity of the maps $X_i$, \begin{equation*}
X_i(s) > X_i(r) \text{ if } s > r \text{ and viceversa,}
\end{equation*}  
is straightforward to see that \begin{multline} \label{laconvessita}
\mathcal{L}[\mu_{\theta}^{1 \to 2}]=\eta \int_{s > r}X_{\theta}(s)-X_{\theta}(r)dsds - \eta \int_{s < r}X_{\theta}(s)-X_{\theta}(r)dsdr \\
=\eta \int_{s > r}\theta X_2(s)-\theta X_2(r) + (1-\theta)X_1(s)-(1-\theta)X_1(r)dsdr \\
- \eta \int_{s < r}\theta X_2(s)-\theta X_2(r) + (1-\theta)X_1(s)-(1-\theta)X_1(r)dsdr \\
= \eta \theta \int \bigl|  X_2(s)- X_2(r) \bigr| dsdr + \eta (1-\theta) \int \bigl| X_1(s)- X_1(r) \bigr| dsdr \\
= \theta \mathcal{L}[\mu^2] + (1-\theta)\mathcal{L}[\mu^1].
\end{multline}
This prove the convexity.
\endproof
\end{proposition}

It is proved that the functional $\mathcal{W}$ is proper, l.s.c., coercive and $\eta$-geodesically convex; this leads to the following theorem:
\begin{theorem}[Existence and uniqueness of gradient flows]
Let $W$ satisfy the assumptions \textbf{A - C'} with $\eta<0$ and let $\mu_0 \in D(\mathcal{W})$. 
The discrete solution $\overline{U}_{\boldsymbol{\tau}}$ given by (\ref{discsol}) converges locally uniformly
to a locally Lipschitz curve $\mu \in \mathcal{P}_2(\mathbb{R}) $, such that $\mu(t)$ is diffuse for a.e. $t>0$, which is the unique 
gradient flow of $\mathcal{W}$ with $\mu(0+)=\mu_0$, i.e. a curve $\mu \in AC^2_{loc}([0,+\infty);\mathcal{P}_2(\mathbb{R}))$ satisfying 
\begin{equation*}
\begin{split}
& \frac{\partial \mu(t)}{\partial t} + \frac{\partial}{\partial x}(\boldsymbol{v}(t)\mu(t))=0 \text{ in } \mathcal{D}'([0,+\infty) \times \mathbb{R}), \\
& \boldsymbol{v}(t)=-\partial^0\mathcal{W}[\mu(t)]=\int_{x \neq y}\partial_x W(x-y)d\mu(y) \text{ for a.e. }t>0 \\
& \| \boldsymbol{v}(t) \|_{L^2(\mu(t))}=|\mu|'(t) \qquad a.e. \; t>0,
\end{split}
\end{equation*}
with $\mu(0)=\mu_0$. Moreover, the energy identity \begin{equation*}
\int_a^b\int_{\mathbb{R}}|\boldsymbol{v}(t,x)|^2d\mu(t)(x) + \mathcal{W}[\mu(b)] = \mathcal{W}[\mu(a)]
\end{equation*}
holds for all $0 \leq a \leq b < +\infty$.
Moreover, this solution is characterized by a system of evolution variational inequalities:
\begin{equation}
\frac{1}{2}\frac{d}{dt}W_2^2(\mu(t),\sigma)+\frac{\eta}{2}W_2^2(\mu(t),\sigma) \leq \mathcal{W}[\sigma] - \mathcal{W}[\mu(t)] \qquad \text{a.e.} t>0,
\end{equation}
for all $\sigma \in \mathcal{P}_2(\mathbb{R})$.
\end{theorem} 
 
We prove now a characterization of the subdifferential of $\mathcal{W}$ in the case when $\mu$ is diffuse, i.e. $\mu(\{ x \})=0$ for 
every $x \in \mathbb{R}$. This is the case of any measure which is absolutely continuous w.r.t. the Lebesgue measure on $\mathbb{R}$. Note 
that no concentrated masses are allowed. \begin{proposition}
Given a potential satisfying \textbf{A - C'}, the vector field \begin{equation} \label{definizionedikappa}
k(x):=\int_{y \neq x}\partial_x W(x-y)d\mu(y)
\end{equation}
is the unique element of minimal $L^2(\mu)$-norm in the subdifferential of $\mathcal{W}$.
\proof
The factorization of $\mathcal{W}=\overline{\mathcal{W}}+\mathcal{L}$ is needed here. The subdifferential of $\overline{\mathcal{W}}$ has 
been already obtained:
\begin{equation} \label{laprima}
\overline{\mathcal{W}}[\mu^2] - \overline{\mathcal{W}}[\mu^1] \geq \inf_{\boldsymbol{\gamma}_0 \in \Gamma_0}
\int_{\mathbb{R}\times \mathbb{R}} \overline{k}(x)\cdot (y-x)d\boldsymbol{\gamma}_0(x,y) +o(W_2(\mu^2,\mu^1)).
\end{equation} 
Using the equation \ref{laconvessita} and some calculation is obtained the following: \begin{multline}
\mathcal{L}[\mu^2]-\mathcal{L}[\mu^1]=\frac{\mathcal{L}[\mu_{\theta}]-\mathcal{L}[\mu^1]}{\theta}=\frac{\eta}{2}\int 
\frac{\bigl|X_{\theta}(s)-X_{\theta}(r)\bigr| - \bigl|X_1(s)-X_1(r)\bigr|}{\theta} \\
= \frac{\eta}{2 \theta} \int_{s > r} \bigl( (1-\theta)X_1(s)+\theta X_2(s)-(1-\theta)X_1(r)-\theta X_2(r)-X_1(s)+X_1(r) \bigr) dsdr \\
- \frac{\eta}{2 \theta} \int_{s < r} \bigl( (1-\theta)X_1(s)+\theta X_2(s)-(1-\theta)X_1(r)-\theta X_2(r)-X_1(s)+X_1(r) \bigr) dsdr \\
= \frac{\eta}{2} \int_{s > r} \bigl( -X_1(s)+X_2(s)+X_1(r)-X_2(r) \bigr) dsdr - \frac{\eta}{2} \int_{s < r} \bigl( \cdots \bigr) dsdr \\
= \eta \int_{s > r} \bigl( -X_1(s)+X_2(s) \bigr) dsdr - \eta \int_{s < r} \bigl( -X_1(s)+X_2(s) \bigr) dsdr \\
= \eta \int_0^1 \bigl( X_2(s) - X_1(s) \bigr) \cdot \biggl( \int_0^1 \text{sign}(s-r) dr \biggr) ds \text{ with sign}(0)=0.
\end{multline}
When $\mu$ is diffuse then $X_{\mu}$ is strictly increasing, so that $sign(s-r)=sign(X_{\mu}(s)-X_{\mu}(r))$ and 
applying the change of variable formula in the opposite sense is obtained \begin{equation} \label{laseconda}
\mathcal{L}[\mu^2]-\mathcal{L}[\mu^1]=\int_{\mathbb{R} \times \mathbb{R}} \biggl( \int_{x \neq y} \eta \text{sign}(x-z)d\mu^1(z) \biggr)\cdot (y-x)
d\boldsymbol{\gamma}(x,y).
\end{equation}
In the end, adding (\ref{laprima}) to (\ref{laseconda}) is obtained
\begin{equation} 
\mathcal{W}[\mu^2] - \mathcal{W}[\mu^1] \geq \inf_{\boldsymbol{\gamma}_0 \in \Gamma_0}
\int_{\mathbb{R}\times \mathbb{R}} k(x)\cdot (y-x)d\boldsymbol{\gamma}_0(x,y) +o(W_2(\mu^2,\mu^1)),
\end{equation}
and so $k \in \partial \mathcal{W}[\mu]$.\\
The minimality of $k$ has already been proved in proposition \ref{proposiminimal}.
\endproof
\end{proposition}

\begin{theorem}
If there exists $\overline{x} \in \mathbb{R}$ such that $\mu(\{ \overline{x} \}) > 0$ then the function $k$ defined by (\ref{definizionedikappa}) 
does not belong to $\partial \mathcal{W} [\mu]$.
\proof
Let us assume that $\mu_1(\{ \overline{x} \})= \delta > 0$ and, for the sake of simplicity, $\mu_1(\{ x \})=0$ if $x\neq \overline{x}$. \\
Then there exists $0 \leq r_1 < r_2 \leq 1$ such that \begin{equation*}
r_2=r_1+\delta, \; X_1(r)\equiv \overline{x} \; \text{ for every }r\in [r_1,r_2].
\end{equation*}
Arguing as before, for every measure $\mu_2 \in \mathcal{P}_2(\mathbb{R})$ we get \begin{gather*}
\mathcal{L}[\mu_2]-\mathcal{L}[\mu]-\eta \int_{\mathbb{R} \times \mathbb{R}}\biggl( \int_{x \neq y}sign(x-z)d\mu(z) \biggr)
(y-x)d\boldsymbol{\gamma}(x,y) \\
=\eta \int_0^1\bigl(X_2(s)-X_1(s)\bigr) \biggl( \int_0^1 sign(s-r) - sign(X_1(s)-X_1(r)) dr \biggr) ds .
\end{gather*}
Notice that $sign(s-r)=sign(X_1(s)-X_1(r))$ if $s\notin [r_1,r_2]$ or $r \notin [r_1,r_2]$. \\
Therefore the previous integral becomes \begin{multline*}
\eta \int_{r_1}^{r_2}\bigl( X_2(s)-X_1(s) \bigr)\biggl( \int_{r_1}^{r_2}sign(s-r)dr \biggr) ds\\
=\eta \int_{r_1}^{r_2} \bigl( X_2(s)-X_1(s) \bigr)(2s-r_1-r_2)ds \\
=2\eta \int_{r_1}^{r_2} \bigl( X_2(s)-\overline{x} \bigr)(s-\frac{r_1+r_2}{2})ds \\ 
=-\eta \int_{r_1}^{r_2}  X_2'(s)\cdot(s-r_1)(s-r_2)ds.
\end{multline*} 
Since $X_2' \geq 0$ and $(s-r_1)(s-r_2) \leq 0$ in $[r_1,r_2]$ the previous quantity is nonnegative for if $\eta \geq 0$ (convex case) 
and non positive if $\eta<0$ (concave case). in this case, $k$ could not belong to $\partial \mathcal{W}[\mu]$.
\endproof
\end{theorem}


\clearpage

\begin{example}
Going back to the example \ref{theesempio}, now can be calculated the behaviour of the solution. In the one dimensional case the equation 
\ref{theeeeeq} can be rewritten (see \cite{LT,CT,BDIF}) \begin{equation}
\partial_t X(t,z)=\int_0^1W'\bigl( X(\xi)-X(z) \bigr)d\xi, \qquad \forall z \in [0,1].
\end{equation}
In this case, with $W(x)=-|x|$ it become \begin{equation}
\partial_t X(t,z)=\int_0^1 -\text{sign}( \xi - z )d\xi,
\end{equation}
solving the integral and integrating with respect to $t$ (assuming that $X$ is strictly increasing) gives the following pseudo inverse: 
\begin{equation}
X(t,z)=t(2z-1)+X(0,z)
\end{equation}
Three examples of the problem with $W(x)=-|x|$ and different starting point are reported here: \begin{itemize}
\item if the starting point is $\mu(0,x)=\delta_{x_0}$, the solution is \begin{equation}
\mu(t,x)=\frac{1}{2t}\mathds{1}(x_0-t,x_0+t);
\end{equation}
\item if the starting point is $\mu(0,x)=\frac{1}{2}\delta_{x_{1}}+\frac{1}{2}\delta_{x_{2}}$, the solution is \begin{equation}
\mu(t,x)=\frac{1}{2t}\mathds{1}(x_1-t,x_1)+\frac{1}{2t}\mathds{1}(x_2,x_2+t);
\end{equation}
\item if the starting point is $\mu(0,x)=\frac{1}{4}\delta_{x_{1}}+\frac{1}{4}\delta_{x_{2}}+\frac{1}{2}\delta_{x_{3}}$, 
the solution is \begin{equation}
\mu(t,x)=\frac{1}{2t}\biggl[ \mathds{1}(x_1-t,x_1-1/2t)+\mathds{1}(x_2-1/2t,x_2)+\mathds{1}(x_3,x_3+t)\biggr].
\end{equation}
\end{itemize}
Taking in account the other solution found in example \ref{theesempio}, it's easy to see that they don't satisfy the energy inequality 
\ref{energidisug} and so they are not curve of maximal slope.
\end{example}

\begin{figure}[p]
\centering
\subfloat[][\emph{Probability measure}.]
{\includegraphics[width=.80\columnwidth]{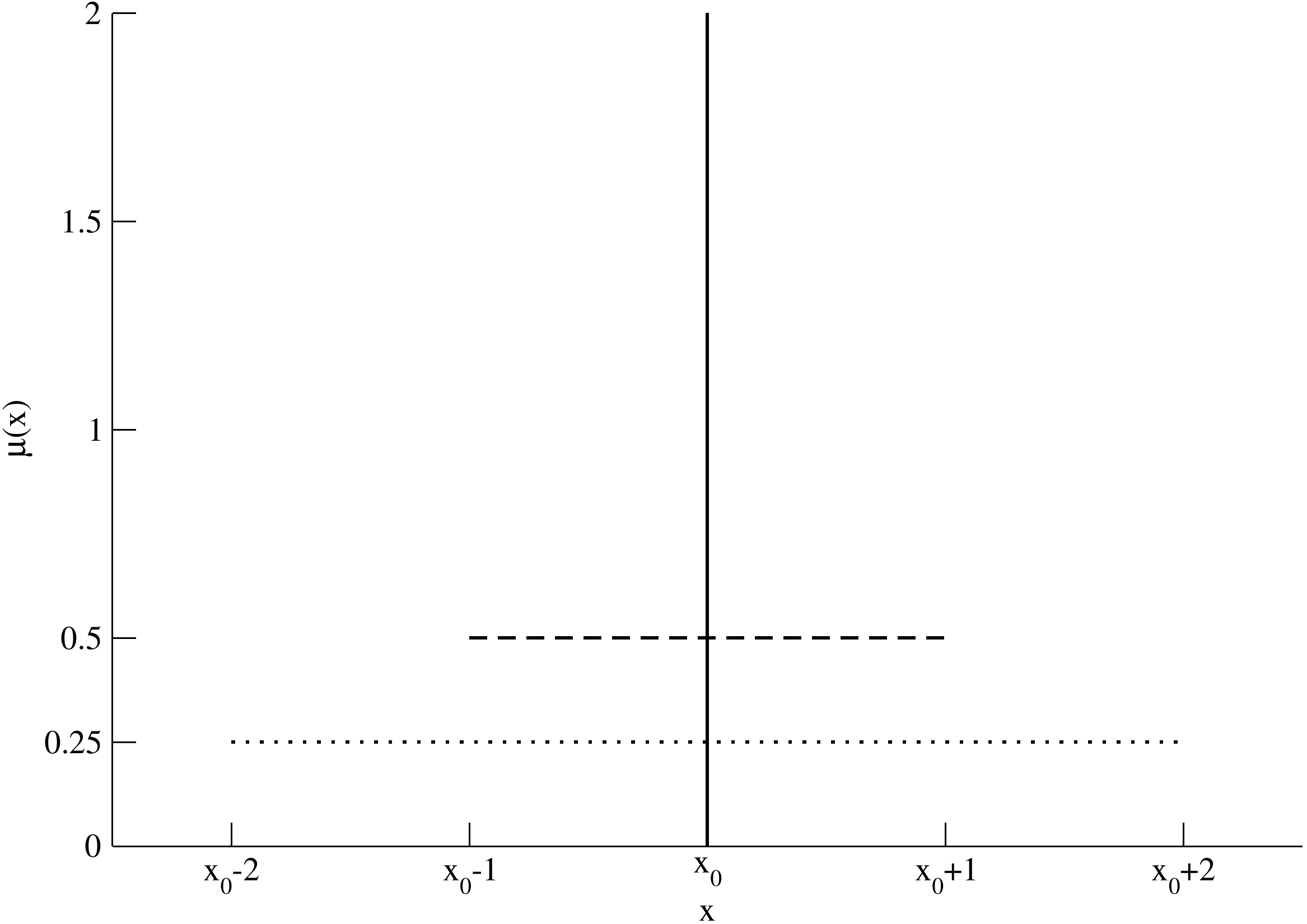}} \\
\subfloat[][\emph{Monotone rearrangement}.]
{\includegraphics[width=.80\columnwidth]{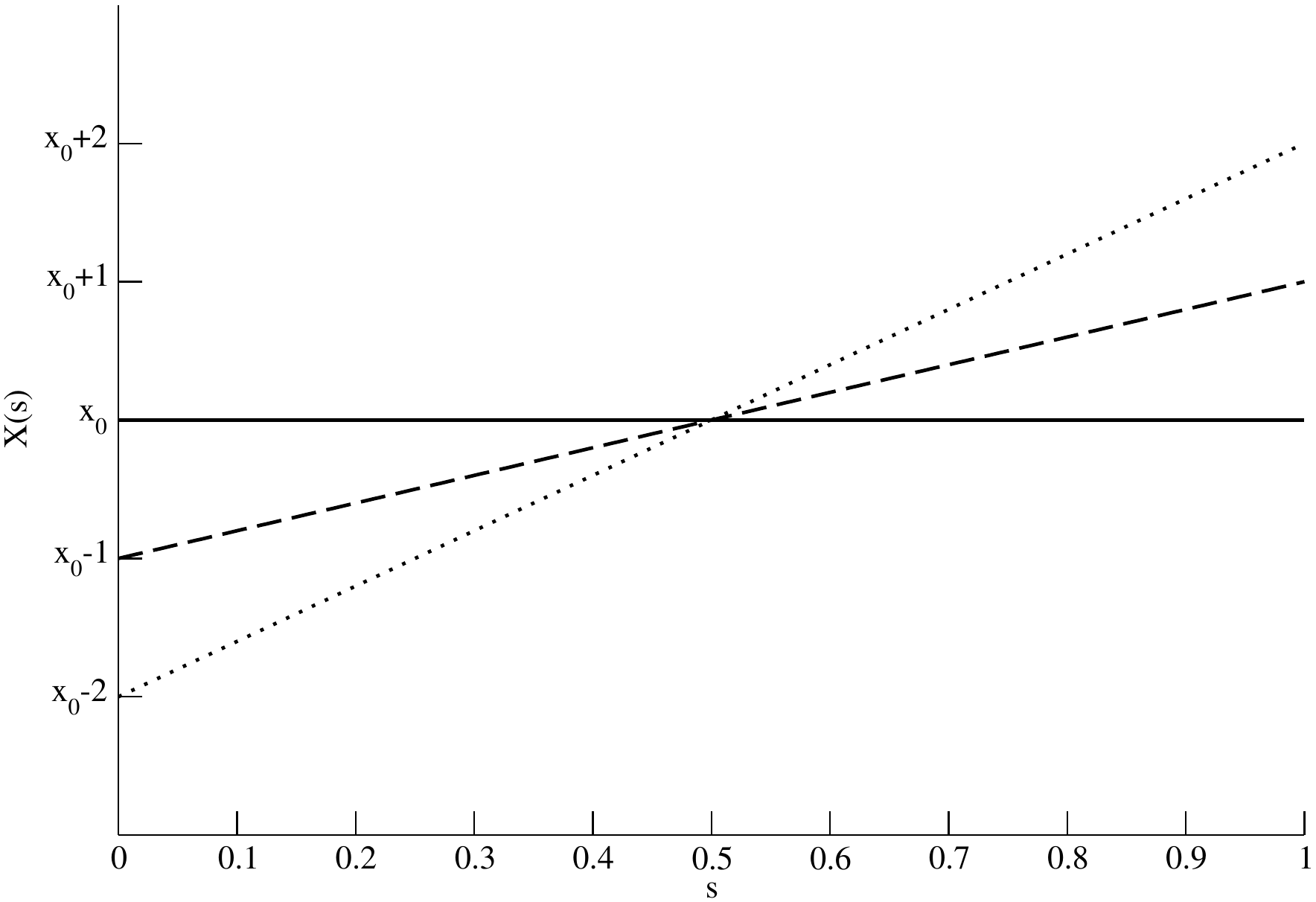}}  
\caption{Evolution from a Dirac delta ($\mu(x,0)=\delta_{x_0}$) for the potential $W(x)=-|x|$. The image plots
$\mu(x,t)$ and $X(s,t)$ at time $t = 0$ (bold line, initial data), $t = 1$ (dashed line) and $t = 2$ (dotted line).}
\end{figure}

\begin{figure}[p]
\centering
\subfloat[][\emph{Probability measure}.]
{\includegraphics[width=.80\columnwidth]{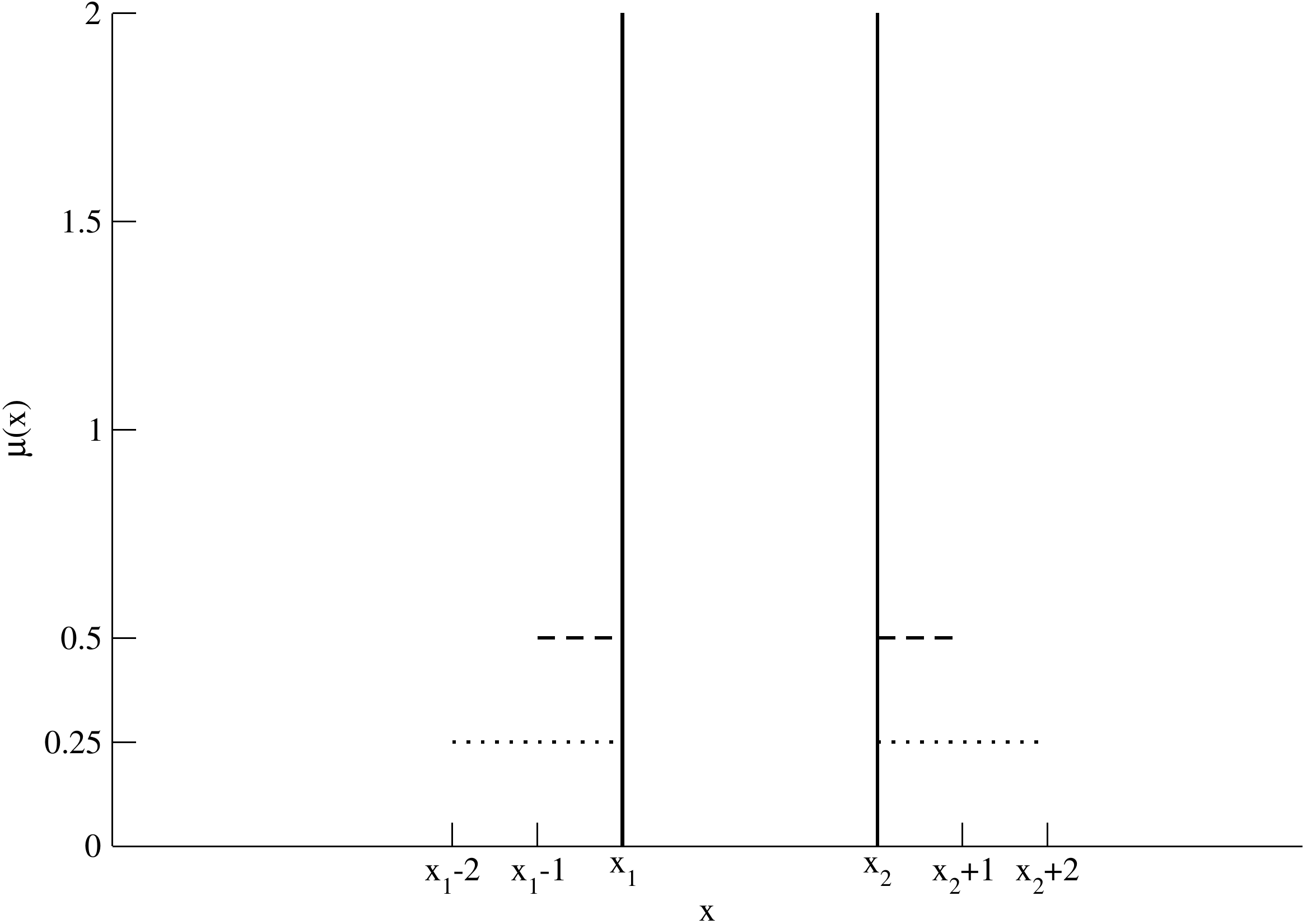}} \\
\subfloat[][\emph{Monotone rearrangement}.]
{\includegraphics[width=.80\columnwidth]{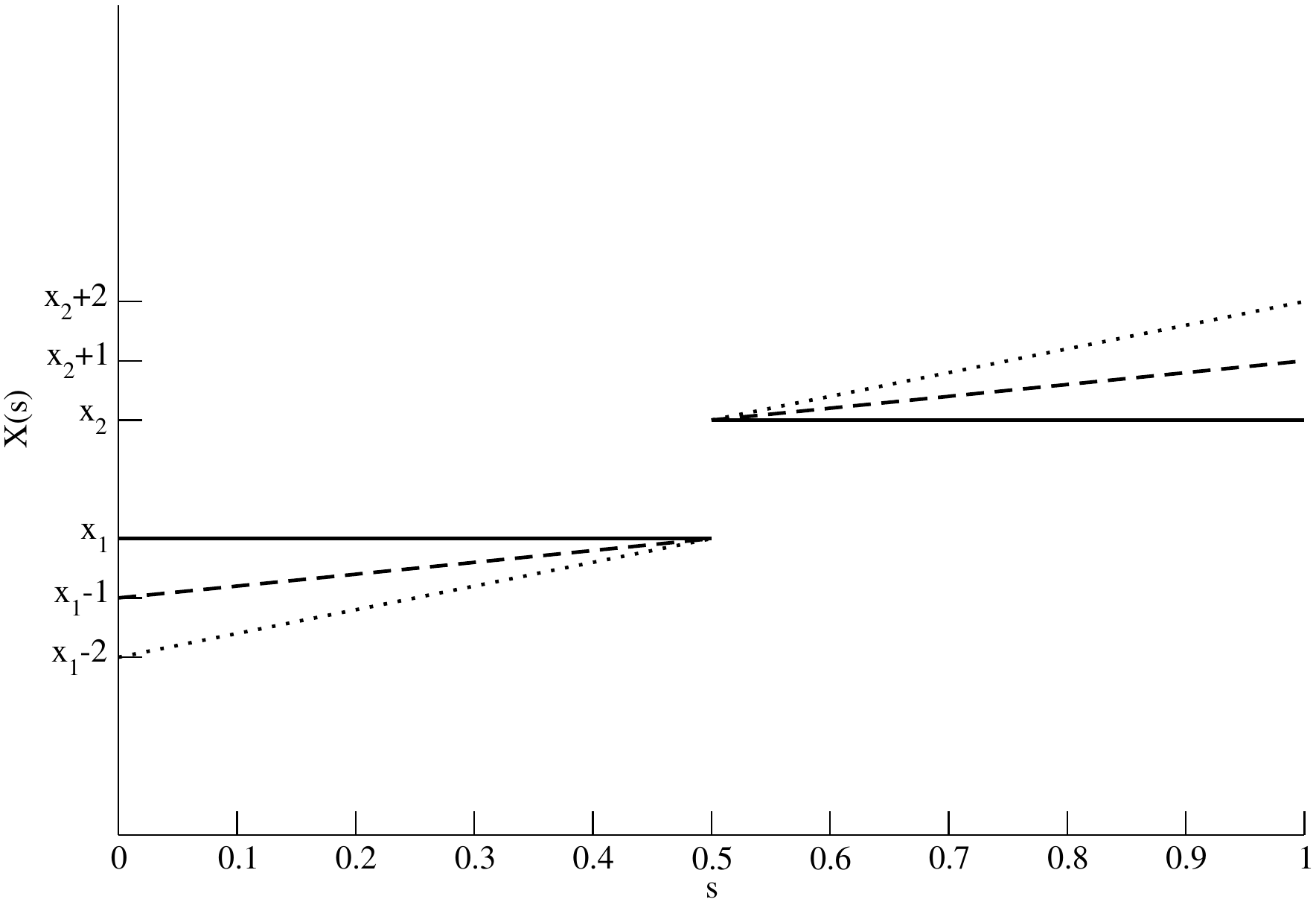}} 
\caption{Evolution from two Dirac delta ($\mu(x,0)=\frac{1}{2}\delta_{x_{1,0}}+\frac{1}{2}\delta_{x_{2,0}}$)for the potential $W(x)=-|x|$. The image plots
$\mu(x,t)$ and $X(s,t)$ at time $t = 0$ (bold line, initial data), $t = 1$ (dashed line) and $t = 2$ (dotted line).}
\end{figure}

\begin{figure}[p]
\centering
\subfloat[][\emph{Probability measure}.]
{\includegraphics[width=.80\columnwidth]{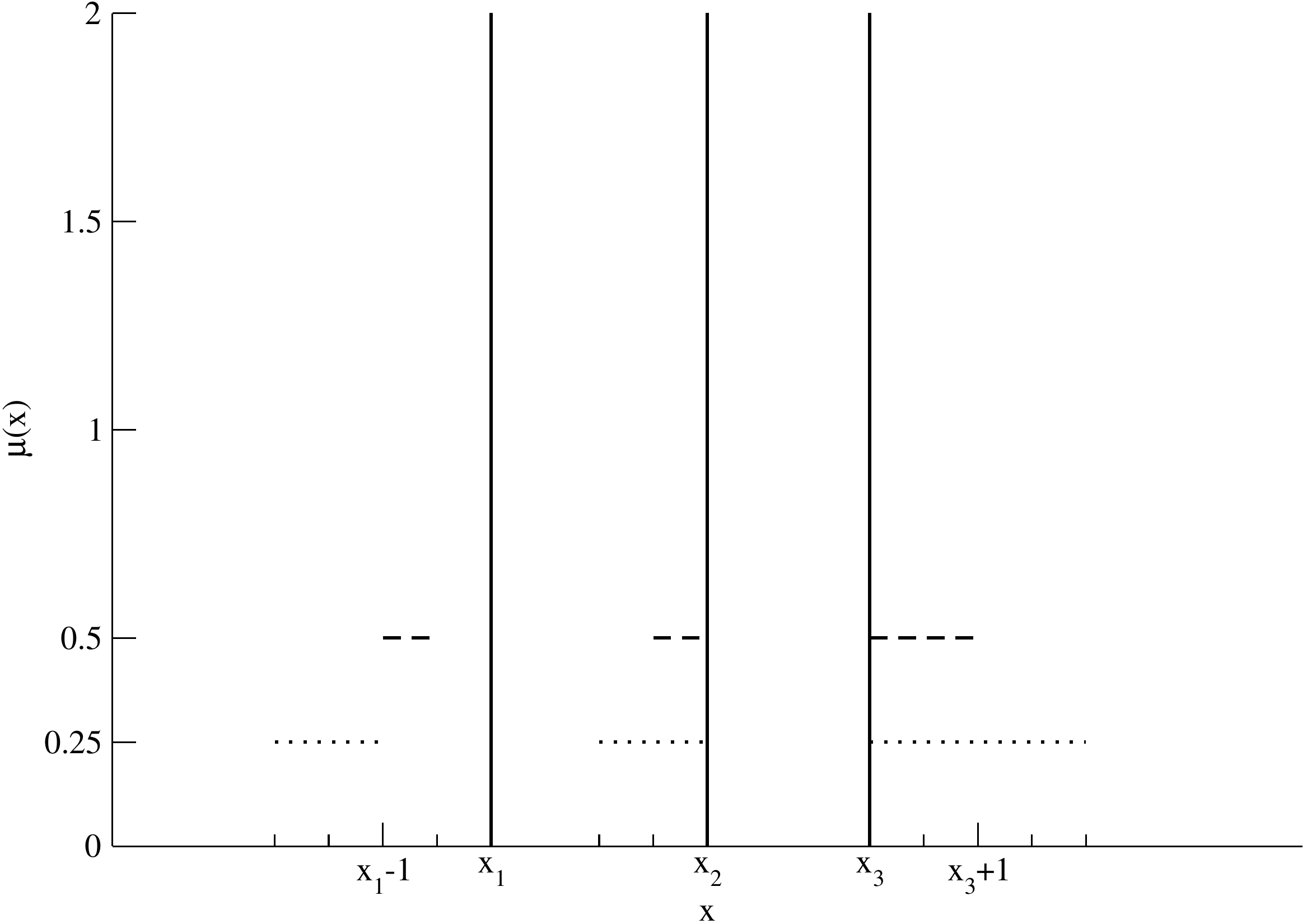}} \\
\subfloat[][\emph{Monotone rearrangement}.]
{\includegraphics[width=.80\columnwidth]{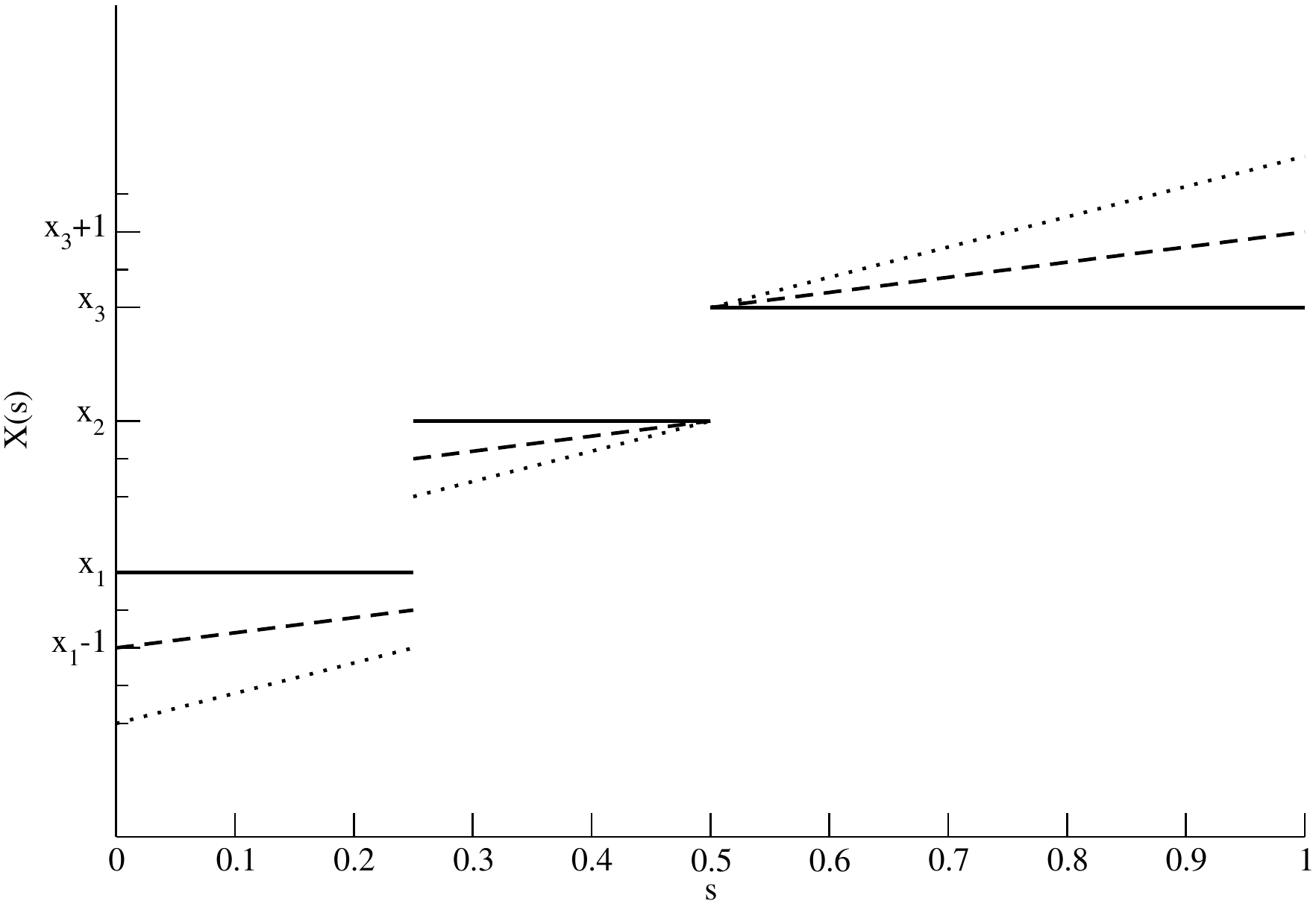}}  
\caption{Evolution from three Dirac delta ($\mu(0,x)=\frac{1}{4}\delta_{x_{1}}+\frac{1}{4}\delta_{x_{2}}+\frac{1}{2}\delta_{x_{3}}$)for 
the potential $W(x)=-|x|$. The image plots
$\mu(x,t)$ and $X(s,t)$ at time $t = 0$ (bold line, initial data), $t = 1$ (dashed line) and $t = 2$ (dotted line).}
\end{figure}

\clearpage

\section{Conclusions}

We showed how to study measure-valued solutions to the system of granular-flows in the presence of a repulsive interaction potential.

In dimension 1, in the framework of the theory of gradient flows, we showed how to select a distinguished solution, when $W$ is of the 
form $\overline{W}-\eta|x|$ with a smooth $\overline{W}$. This condition represent the case of a repulsive potential with a concave 
cusp in the origin. 

Our strategy is based on three fundamental steps: the well posedness of the minimizing movement scheme, the Wasserstein convexity of 
$\mathcal{W}$ and the characterization of the subdifferential. The MMS turns out to work well; the discrete solution can be defined and 
it converges to an absolutely continuous curve. The main problem is to see if it converges to a solution in the sense of curve of 
maximal slope (that is proved to be equivalent, in this case, to gradient flows solution). This can be achieved by characterizing 
the subdifferential and by showing the $\eta$-convexity of $\mathcal{W}$. In dimension $d>1$ two results 
give some hint for a possible characterization, but the fundamental property of regularity cannot be proved, a big problem at this 
level. Moreover, in the one dimensional framework, a change of variable formula allows us to get more insights and explicit formula in 
some relevant cases.

Simple examples show that the gradient flow solution select a diffusive evolution. It means, thinking at the time reversal, 
that after the aggregation of solutions there is a loss of information, it is impossible to determine the starting probability density 
that generated the aggregation, only one is selected among all the other. An explicit formula has been obtained 
for a general potential and it turns out that a Dirac delta cannot remain a Dirac delta at every time $t>0$. 

This thesis is a step towards a more detailed study of interaction potentials from a theoretical point of view, with 
a particular attention on possible applications. The general theory of gradient flows cannot be applied directly, but we showed how 
to extend to these potentials the methods and results of the theory, covering many interesting new cases. \\

\clearpage

\clearpage{\pagestyle{empty}\cleardoublepage}
	
\lhead[\thepage]{\small{\textsc{Bibliography}}}
\rhead[\small{\textsc{Bibliography} }] {\thepage}

\end{document}